\documentclass[11pt]{amsart}
\usepackage{amsmath,amssymb,amsthm}
\usepackage[shortlabels]{enumitem}
\usepackage[margin=1.2in]{geometry}
\usepackage[colorlinks=true,linkcolor=blue,citecolor=blue]{hyperref}

\newtheorem{theorem}{Theorem}[section]
\newtheorem{proposition}[theorem]{Proposition}
\newtheorem{lemma}[theorem]{Lemma}
\newtheorem{corollary}[theorem]{Corollary}
\newtheorem{conjecture}[theorem]{Conjecture}
\newtheorem{mainthm}{Theorem}

\theoremstyle{definition}
\newtheorem{definition}[theorem]{Definition}

\newtheorem{question}[theorem]{Question}
\newtheorem{convention}[theorem]{Convention}
\theoremstyle{remark}
\newtheorem{remark}[theorem]{Remark}
\numberwithin{equation}{section}

\newcommand{\Q}{\mathbb{Q}}
\newcommand{\R}{\mathbb{R}}
\newcommand{\Z}{\mathbb{Z}}
\newcommand{\cD}{\mathcal{D}}
\newcommand{\cS}{\mathcal{S}}
\newcommand{\cT}{\mathcal{T}}
\newcommand{\cU}{\mathcal{U}}
\newcommand{\cA}{\mathcal{A}}
\newcommand{\cH}{\mathcal{H}}
\newcommand{\cR}{\mathcal{R}}
\newcommand{\cF}{\mathcal{F}}
\newcommand{\cG}{\mathcal{G}}
\newcommand{\cQ}{\mathcal{Q}}
\newcommand{\fF}{\mathfrak{F}}
\newcommand{\fG}{\mathfrak{G}}
\newcommand{\fd}{\mathfrak{d}}
\newcommand{\rel}[2]{\langle #1,#2\rangle}
\newcommand{\Cm}{C^{m}}
\newcommand{\Csm}{C^{\mathrm{sm}}}
\newcommand{\Calt}{C^{\mathrm{sm,alt}}}

\newcommand{\GCRD}{GC^{R_\Delta}}
\newcommand{\GCE}{GC^{E}}
\newcommand{\GCG}{GC^{G}}
\newcommand{\GCcR}{GC^{\cR}}
\newcommand{\Calttop}{C^{\mathrm{alt}}}
\newcommand{\Bord}{\mathrm{Bord}}
\newcommand{\cut}{\mathrm{cut}}
\newcommand{\Sd}{\mathrm{Sd}}
\newcommand{\sgn}{\operatorname{sgn}}
\newcommand{\id}{\mathrm{id}}
\newcommand{\pr}{\mathrm{pr}}
\newcommand{\St}{\mathrm{St}}
\newcommand{\Lk}{\mathrm{Lk}}
\newcommand{\inter}{\operatorname{int}}

\begin{document}

\title[Gluing geometric chains]{Gluing geometric chains of manifolds and orbifolds with corners}
\author{Hao Yu}
\date{}

\begin{abstract}
Geometric chains are formal combinations of smooth maps from compact oriented manifolds or orbifolds with corners. In earlier work the author showed that, modulo a relation $D$ that kills degenerate chains, they compute rational singular homology. We show that $D$ can be replaced by a gluing relation: $[P,f]=\sum_i[Q_i,f\circ\iota_i]$ whenever $P$ is glued from finitely many pieces $Q_i$ along whole codimension-one faces. With this relation the comparison map from singular homology is an isomorphism in all degrees for compact polyhedra and smooth manifolds without boundary, with rational and with integer coefficients. Unlike $D$, the gluing relation is defined through the sources alone, and chain-level identities such as the equivalence of a degree $d$ covering with $d$ times its base, which fails modulo $D$, and the correction terms of the chain map from a Morse complex, which are degenerate modulo $D$, become gluing relations. The proof is axiomatic; injectivity comes from integration or dual cells, surjectivity from cutting along regular levels, where the cut of a gluing relation is again one, and at the point we use Whitehead triangulations of manifolds with corners subordinate to their corner strata. The gluing relation, read on underlying spaces, also suffices for orbifold sources, with integer coefficients modulo classes with an orientation reversing automorphism, and for effective orbifold targets.
\end{abstract}

\maketitle
\tableofcontents

\section{Introduction}\label{sec_intro}

\subsection{\texorpdfstring{Geometric chains and the relation $D$}{Geometric chains and the relation D}}
A geometric simplex of dimension $n$ in a space $Y$ is a pair $(Q,f)$, where $Q$ is a compact oriented $n$-manifold with corners and $f:Q\to Y$ is a map; the boundary of $(Q,f)$ is $(\partial Q,f\circ i_Q)$. Chains of this kind go back to Harrelson, Voronov and Z\'u\~niga \cite{HVZ}. A singular simplex $\sigma$ gives the geometric simplex $p(\sigma)=[\Delta^n,\sigma]$, and the question is in what sense the comparison map $p$ identifies geometric homology with singular homology.

In \cite{Yu1} the sources are orbifolds with corners, the maps are smooth, and the geometric chain complex is the span of isomorphism classes modulo three relations: orientation reversal, disjoint union, and the vanishing of \emph{degenerate} chains, i.e.\ combinations of thin simplices with thin boundary (a simplex is thin if its map factors smoothly through a lower-dimensional orbifold). Call the last relation $D$. The main theorem of \cite{Yu1} says that with $D$, $p_*$ is an isomorphism with rational coefficients for every smooth space with a good cover.

The role of $D$ is analysed in \cite[Remark on the role of degenerate chains]{Yu1}. Without $D$, the cutting (bisection) prisms, homotopy invariance and both halves of the Mayer--Vietoris argument still go through; $D$ is used essentially only to show that geometric chains in a point have no homology in positive degrees. Without $D$ this fails for orbifold sources: for a prime $p\ge5$ the signed count $\phi_p$ of isolated points with local group $\mu_p$ acting by scalars on $\mathbb{C}^2$ is a linear form on $4$-chains in a point that vanishes on boundaries and equals $\pm1$ on the weighted projective plane $\mathbb{P}(1,1,p)$. For manifold sources closed manifolds bound without $D$, and the question whether a third relation is needed at all is left open in \cite{Yu1}.

The first version of \cite{Yu1}, arXiv:2004.07698v1 \cite{Yu0}, used a different third relation, a \emph{gluing relation}: if $|P|$ is the gluing of $|Q_1|$ and $|Q_2|$ along their lower-dimensional faces, compatibly with the maps, then $(P,f_P)=(Q_1,f_{Q_1})+(Q_2,f_{Q_2})$ \cite[Definition~1.1]{Yu0}. The present version of \cite{Yu1} no longer imposes it, and notes that, being an equality in a quotient, it made the comparison with singular homology delicate, because cancellations could take place through relations involving chains outside a given open set. In this paper we return to the gluing relation, in the smooth setting, with manifold sources and finitely many pieces, and show that the difficulty can be avoided: for manifold sources the gluing relation suffices. In Section~\ref{sec_orb} we extend this to orbifold sources (Theorem~\ref{thm_orb}): if the pieces of a gluing are allowed to be manifolds, as in \cite[Definition~1.1]{Yu0}, where the gluing is of underlying spaces, the comparison theorem holds again, and the cycles $[\mathbb P(1,1,p),c]$ vanish at chain level. If the pieces are required to have the same local groups, they do not (Remark~\ref{rem_orb}). In Section~\ref{sec_int} we show that for manifold sources the comparison holds with integer coefficients (Theorem~\ref{thm_Z}), proving injectivity by intersecting with the dual cells of a smooth triangulation, and that for orbifold sources it holds over $\Z$ modulo the classes with an orientation reversing automorphism (Theorem~\ref{thm_Z_orb}). In Section~\ref{sec_orbtarget} we show that the gluing relation alone also suffices for effective orbifold targets, possibly with corners and not compact, with liftable maps or orbifold morphisms (Theorem~\ref{thm_orbtarget}); there the proof follows the good cover pattern of \cite{Yu1}, with injectivity from a de Rham pairing, and a finite covering of degree $d$ is equivalent to $d$ times its base (Proposition~\ref{prop_cover}). Section~\ref{sec_relprod} gives the relative version, excision through the comparison theorem, the cross product, pullbacks along fibre bundles (for orbifold sources with no condition on the faces of the fibre), and a chain map from the Morse complex, whose correction terms lie in $\operatorname{span}\cG$ at chain level, whereas in \cite{Yu1} they are degenerate (Propositions~\ref{prop_rel}, \ref{prop_cross}, \ref{prop_pull}, \ref{prop_orb_ops} and \ref{prop_morse_G}); the cutting proof of excision in \cite{Yu1} uses the locality of $D$, whose analogue for $G$ is not known (Remark~\ref{rem_locality}).

\subsection{Why the gluing relation}
The relation $D$ is defined through the maps: a chain is degenerate if it and its boundary are combinations of simplices whose maps factor smoothly through lower-dimensional orbifolds. The gluing relation $G$ itself is defined through the sources alone: it says that a chain does not change when its source is cut into pieces or reassembled from them. This is the form in which chain-level identities arise when a space, such as a compactified moduli space, is assembled from pieces, and with $G$ such identities are relations of the complex, rather than failing at chain level or holding only modulo thin chains. This has the following consequences:
\begin{enumerate}[(i)]
\item \emph{Finite coverings.} For an orbifold covering $\pi:P'\to P$ of degree $d$, $[P',\hat f\circ\hat\pi]=d\,[P,\hat f]$ modulo $\operatorname{span}\cG^{\mathrm{mor}}$ (Proposition~\ref{prop_cover}). Modulo $D$ this fails at chain level, for instance for $[S^1,z^2]-2[S^1,z]$; we expect it to be the identity needed for the connected form of the master equation for string vertices (cf.\ \cite[String vertices]{Yu1}, \cite{YuV}).
\item \emph{Morse complexes.} The map $\Phi$ from the Morse complex is a chain map into $\GCG$ whose correction terms, products with moduli spaces of positive dimension, are themselves sums of gluing relations (Lemma~\ref{lem_prod_thin}, Proposition~\ref{prop_morse_G}). In \cite{Yu1} they are degenerate, i.e.\ killed only by $D$.
\item \emph{Integer coefficients for orbifold sources.} With orbifold sources and the orbifold gluing relation, the comparison holds over $\Z$ for compact polyhedra and smooth manifolds without boundary, once the classes with an orientation reversing automorphism are set to zero (Theorem~\ref{thm_Z_orb}); \cite{Yu1} works over $\Q$, the natural choice with orbifold sources, and integrality for $D$ is not known. For manifold sources integer coefficients are the expected case: here they hold for $G$ (Theorem~\ref{thm_Z}), and for smooth manifolds without boundary and the relation of trivial and small rank degenerate chains, which differs from $D$, they are proved in \cite[Theorems 5.42, 5.45]{FMMS}.
\item \emph{An axiomatic criterion.} Any relation family with six properties, listed below, gives the comparison theorem, and the last of them, acyclicity at the point, is also necessary (Theorem~\ref{thm_abstract}). This separates what the comparison needs from the particular relation.
\end{enumerate}
The price is locality and scope. $D$ is local, a chain of an open set being degenerate there if it is degenerate in the whole space, and \cite{Yu1} uses this for excision by cutting, on arbitrary smooth spaces. For $G$ the corresponding statement is not known (Remark~\ref{rem_locality}); here excision follows from the comparison theorem (Proposition~\ref{prop_rel}), so only for the targets of our comparison theorems, and the smooth spaces of \cite{Yu1} beyond these are not treated (Remark~\ref{rem_rel_quotient}(2)). Separately, pullbacks along orbifold bundles over orbifold targets, which \cite{Yu1} uses for the BV operator, are not treated, since a pulled back gluing datum need not be one (Remark~\ref{rem_orbtarget}(3)).

\subsection{The gluing relation and the main theorem}
Throughout, sources are compact oriented smooth manifolds with corners with embedded faces (Definition~\ref{def_source}), targets are subsets $Y\subseteq\R^K$, and maps are smooth into $\R^K$ with image in $Y$ (Convention~\ref{conv_smooth}); a smooth manifold is treated as a closed submanifold of some $\R^K$. Coefficients are $\Q$.

A \emph{gluing decomposition} of a source $P$ is a finite family of sources $Q_i$ of the same dimension with smooth orientation-preserving embeddings $\iota_i:Q_i\to P$ of manifolds with corners, whose images cover $P$ and whose interiors are disjoint from the other images, such that every codimension-one face of every piece either lies in the boundary of $P$ or is glued, as a whole and by a diffeomorphism, to a codimension-one face of exactly one other piece (Definition~\ref{def_glue}). This is the literal smooth reading of ``$|P|$ is the gluing of the $|Q_i|$ along their faces''. The \emph{gluing relation} is
\[
[P,f]=\sum_i[Q_i,f\circ\iota_i],
\]
and $G(Y)$ is the subcomplex of geometric chains generated by these relations; $\GCG(Y)$ is the quotient complex. Triangulations of $P$ by smoothly embedded simplices are gluing decompositions, and so are bisections along regular levels.

\begin{mainthm}[Theorem~\ref{thm_A_G}]\label{thmA}
Let $Y$ be a compact polyhedron in $\R^K$, or a smooth manifold without boundary. Then
\[
p_*:H_*(Y;\Q)\longrightarrow H_*(\GCG(Y))
\]
is an isomorphism in all degrees.
\end{mainthm}

In particular the dimension axiom $H_d(\GCG(N))=0$ for $d>\dim N$ holds for every compact polyhedron $N$.

Other readings of ``gluing'' give the same result (Theorem~\ref{thm_readings}). If the pieces are only required to be smooth and injective, with the faces matched by homeomorphisms, or if they are topological embeddings such that the map composed with them is smooth, Theorem~\ref{thmA} still holds; in the topological reading injectivity comes from intersecting with dual cells instead of integration (Remark~\ref{rem_Z}). So does the larger \emph{decomposition relation} $E\supseteq G$, which allows non-immersive pieces and walls that do not match, and $H_*(\GCG)\to H_*(\GCE)$ is an isomorphism (Corollary~\ref{cor_GE}).

The relation of \cite{Yu0} has two pieces. For the two-piece version of Definition~\ref{def_glue}, $p_*$ is injective in all degrees and an isomorphism in degree $\le1$, and in degree $2$ modulo the classification of compact surfaces with corners, which we use without proof; for $m\ge3$ the comparison theorem in degrees $\le m$ is equivalent to the vanishing of $H_k$ of the two-piece complex of a point for $3\le k\le m$ (Proposition~\ref{prop_twopiece}). Whether this holds is open from degree $3$ on.

\subsection{Whitehead triangulations with corners}
At the point, the proof needs triangulations of sources whose simplices are themselves admissible pieces. We obtain them as a by-product.

\begin{mainthm}[Theorem~\ref{thm_whitehead}]\label{thmB}
Every compact smooth manifold with corners with embedded faces has a triangulation whose closed simplices are smooth embeddings and whose open simplices lie in single corner strata.
\end{mainthm}

We do not know a reference for this statement. The triangulations obtained by pushing an equivariant triangulation of the reflection cover down to $Q$ are not immersive at the mirrors; instead we triangulate a fundamental domain upstairs, by Illman's theorem \cite{Ill}, and identify it with $Q$ by a diffeomorphism built from compatible collars (Section~\ref{sec_point}). The weak reading of gluing and the decomposition relation need only the push-down (Corollary~\ref{thm_T7}), which uses Illman's theorem alone.

\subsection{\texorpdfstring{What replaces $D$}{What replaces D}}
The proof of Theorem~\ref{thmA} shows where $D$ was used in \cite{Yu1} and what replaces it.

At the point, triangulation relations, which are gluing relations, do what $D$ does. By Theorem~\ref{thmB} every source $Q$ of positive dimension has a Whitehead triangulation, so $[Q,c]$ is equivalent to a signed sum of constant simplices $[\Delta^n,c]$; each of these vanishes, because a transposition of two vertices is an orientation-reversing automorphism. So $\GCG_n(\mathrm{pt})=0$ for $n\ge1$ already at chain level.

In \cite{Yu1} bisection is realised by a prism: the bisection formula is an exact identity in the complex without $D$, the side prisms over faces on which the bisecting function has no zero vanishing by an orientation-reversing automorphism, and the cutting and prism operators must then be shown to preserve $D$ \cite[Lemma on cutting operators; Lemma on thin interfaces]{Yu1}. Here a bisection is itself a gluing relation, and the cut of a gluing relation along a generic level is again one (Lemma~\ref{lem_cutclosure}). This is the step that makes the Mayer--Vietoris argument close up without any triangulation of the pieces and without any locality statement.

For orbifold sources without any third relation the point genuinely fails, by the form $\phi_p$. That form would not survive a subdivision relation built from manifold simplices. The underlying space $|\mathbb P(1,1,p)|$ is not a topological manifold at the singular point, whose link is the lens space $L(p,1)$, but it can be triangulated; since simplices have no $\mu_p$-points, a relation $[\mathbb P(1,1,p),c]-\sum_i\varepsilon_i[\Delta^4,c]$ coming from such a triangulation would have $\phi_p=\phi_p(\mathbb P(1,1,p))\ne0$. So $\phi_p$ would not descend, and $[\mathbb P(1,1,p),c]$ would be killed by any such relation. Thus the $\mathbb{P}(1,1,p)$ obstruction is an obstruction to having no relation; it would not obstruct replacing $D$ by an orbifold version of a gluing relation. Section~\ref{sec_orb} develops such a version, reading the gluing on underlying spaces so that pieces may be manifolds; it needs no triangulation of orbifolds, only invariant simplicial fans in linear charts (Proposition~\ref{prop_orb_point}).

Theorem~\ref{thmA} does not answer the question of \cite{Yu1} whether manifold sources need a third relation at all, since $G$ is a genuine relation. It does locate the obstruction: any relation family that satisfies the six axioms of Section~\ref{sec_axioms} gives the comparison theorem, and without a third relation the splitting (A2) fails as an identity but can be replaced by the bisection prism of \cite{Yu1}, so that what remains is exactly (A6), the vanishing of the positive-degree homology of the point (Remark~\ref{rem_axiomatic}).

\subsection{Method}
The proof is axiomatic (Theorem~\ref{thm_abstract}). Consider a family of relations $[P,f]=\sum_i[Q_i,f\circ\iota_i]$ such that
\begin{enumerate}[(A1)]
\item membership does not depend on the target and is natural under smooth maps;
\item the cut relations $[Q,f]=[Q_-,f|]+[Q_+,f|]$ along regular levels belong to it;
\item the cut of a relation along a generic level, and its interface, are again relations;
\item the prism of a relation is a relation;
\item integrals of differential forms are additive over relations;
\item the quotient complex of a point has no homology in positive degrees.
\end{enumerate}
Then $p_*$ is an isomorphism in all degrees for compact polyhedra and for smooth manifolds without boundary, and (A6) is necessary.

Injectivity comes from (A5): Stokes' theorem for manifolds with corners handles boundaries of relations, de Rham's theorem gives injectivity for smooth manifolds, and a neighbourhood retraction extends it to compact polyhedra (Section~\ref{sec_inj}). Surjectivity is proved by cutting along a generic level $\{\ell=c\}$ of a linear function on $\R^K$ (Section~\ref{sec_cut}). Cutting along levels of functions pulled back from the target goes back to Lipyanskiy's creasing \cite{lip}, treated in detail by Friedman, Medina-Mardones and Sinha, who prove Mayer--Vietoris this way for manifolds \cite[\S4.5.2, \S5.2]{FMMS}; in \cite{Yu1} the cut is realised by a bisection prism. Here no crease or prism is needed, since a cut is itself a gluing relation (A2). The cut operators $\cut_\pm$ and the interface operator $\Lambda$ act on chains transverse to the level and satisfy $\partial\,\cut_-=\cut_-\partial+\Lambda$. By (A2) and (A3), a cycle modulo the relations splits into cycles modulo the relations on the two halves, up to a singular cycle built from the interface; this is a Mayer--Vietoris lemma in every degree (Lemma~\ref{lem_MV}). An induction over the vertices of a standard simplicial complex, with stars and links as the halves and the slice, proves the theorem for compact polyhedra, using (A4) and (A6) for the star-shaped half; manifolds follow through a tubular neighbourhood (Section~\ref{sec_proof}).

For $G$, axioms (A1), (A2), (A4) and (A5) are direct. Axiom (A3) holds because a cut is defined by pulling back one function: the pieces, their faces, the partners and the gluing diffeomorphisms are all cut by corresponding functions, so cuts of embeddings are embeddings and walls stay aligned (Lemma~\ref{lem_cutclosure}). Axiom (A6) is Theorem~\ref{thmB}. None of the extra flexibility of $E$ is used anywhere except at the point, where Theorem~\ref{thmB} replaces it.

\subsection{The simplex-only relation}
The relation $R_\Delta$ uses only triangulations: it imposes $[Q,f]=\sum_i\varepsilon_i[\Delta^n,f\circ\sigma_i]$ for every depth-subordinate triangulation of $Q$ by smooth, not necessarily immersive, simplices $\sigma_i$ with local degrees $\varepsilon_i$ (Section~\ref{sec_R}). It lies in the injective-only reading of $G$, and its Whitehead version lies in $G$. It makes $p$ surjective at chain level, and the question is concentrated in a kernel complex $K_*(Y)$. Here (A2) and (A3) fail as stated: the cut of a triangulation relation is a curvilinear tiling, which lies in $G^W$ but lies in $R_\Delta$ only if the triangulation can be subdivided coherently along the level; we write $\mathrm{RT}_j$ for this property in dimension $j$ (Definition~\ref{def_RT}).

\begin{mainthm}[Theorem~\ref{thm_R}, Proposition~\ref{prop_RG}]\label{thmC}
Let $Y$ be a compact polyhedron or a smooth manifold without boundary.
\begin{enumerate}[(a)]
\item $p_*:H_k(Y;\Q)\to H_k(\GCRD(Y))$ is injective for every $k$ and an isomorphism for $k\le3$.
\item If $\mathrm{RT}_j$ holds for $j\le m-1$, then $p_*$ is an isomorphism in degrees $\le m$.
\item The natural map $H_*(\GCRD(Y))\to H_*(GC^{G^W}(Y))$ to the injective-only gluing complex, whose homology is that of $\GCG(Y)$ through $p$, is onto in all degrees, with kernel $H_{*-1}(K(Y))$.
\end{enumerate}
\end{mainthm}
We expect that Conjecture~\ref{conj_reltri}, a corner version of the relative $C^1$ triangulation theorem, implies $\mathrm{RT}_j$ for all $j$, and hence that $p_*$ is an isomorphism in all degrees (Remark~\ref{cor_RT_all}); we give only sketches of the conjecture and of this reduction.

Part (c) is a consequence of Theorem~\ref{thm_readings}; it does not give the theorem for $R_\Delta$, whose remaining content is exactly the acyclicity of $K$.

\subsection{The continuous case}
For continuous maps, with topological triangulations, much less is known. The cutting methods of this paper do not apply, since continuous maps have no regular values. For Lipschitz maps the integration argument survives, and $p_*$ is injective in all degrees for compact polyhedra and Riemannian manifolds without boundary (Theorem~\ref{thm_lip_inj}), and an isomorphism in degrees $\le1$; in degree $2$ we expect but do not prove it, and in degrees $\ge3$ surjectivity is open. Section~\ref{sec_cont} sets up the continuous version and lists open problems.

\subsection{Inputs}
The following classical results are used as stated: Joyce's theory of manifolds with corners \cite{Jo}; Illman's smooth equivariant triangulation theorem \cite{Ill}, only in Section~\ref{sec_point}; the existence of compatible collars on manifolds with faces \cite{La}, only for Theorem~\ref{thmB}; Sard's theorem \cite{Sard}; Stokes' theorem on manifolds with corners, de Rham's theorem for smooth singular chains, the Whitney embedding and tubular neighbourhood theorems and the existence of smooth functions with prescribed zero set \cite{Lee,dR}; the small chains lemma and the fact that compact polyhedra are Euclidean neighbourhood retracts \cite{Hat}; acyclic models \cite{EM}, in Sections~\ref{sec_inj} and \ref{sec_R}; and Brown's collar theorem \cite{Bro}, for Lemma~\ref{lem_lip_goodpair}. Section~\ref{sec_orbtarget} also uses the good covers of orbifolds and the comparison for chains given by liftable maps and orbifold morphisms of \cite{Yu1}, and invariant partitions of unity on orbifolds. Section~\ref{sec_relprod} also uses the inward flow on manifolds with corners and the operations on geometric chains of \cite{Yu1}, the Eilenberg--Zilber shuffle map and excision in singular homology \cite{Hat}, and, for the Morse complex, de Rham's theorem for currents, the moduli spaces of \cite{qin}, the proof of \cite[Proposition on the Morse complex]{Yu1} and \cite{hl,hl2}. For Lipschitz maps (Subsection~\ref{ssec_lip}) we also use Rademacher's theorem \cite{Rad,EG}, that the almost-everywhere derivative of a Lipschitz map is its distributional derivative \cite{EG}, and, for the readings $G^W$ and $E$ only, the area formula \cite{Fed,EG}. For the sketch in Remark~\ref{cor_RT_all} we also use Lurie's notes on Whitehead triangulations \cite{Lur} and \cite{Mun,Whi}; Section~\ref{sec_int} uses smooth triangulations of manifolds without boundary \cite{Whi,Mun}, the isomorphism of simplicial and singular homology \cite{Hat}, and the parametric transversality theorem \cite{Lee}.

\section{Setting}\label{sec_setting}

\subsection{Manifolds with corners}
Manifolds with corners are in the sense of Joyce \cite{Jo}, as in \cite{Yu1}. For a manifold with corners $Q$ of dimension $n$, $\Bord_k(Q)$ is the set of points of depth $k$, a manifold of dimension $n-k$; its components are the \emph{depth strata}, and $\inter Q=\Bord_0(Q)$. The boundary $\partial Q$ is Joyce's boundary: the set of pairs $(x,\beta)$ with $\beta$ a local boundary component at $x$, a manifold with corners of dimension $n-1$, with the map $i_Q:\partial Q\to Q$, $(x,\beta)\mapsto x$. If $Q$ is oriented, $\partial Q$ carries the outward-normal-first orientation. We write $B_Q:=i_Q(\partial Q)$; as a topological space $Q$ is a compact $n$-manifold with boundary $B_Q$. We use:
\begin{enumerate}[(J1)]
\item a point of depth $k$ has exactly $k$ preimages under $i_Q$, each of depth $k-1$ in $\partial Q$;
\item every point of $\partial Q$ has a neighbourhood on which $i_Q$ is the inclusion of a boundary hyperplane $\{u_j=0\}$ of a chart $\R^n_k=[0,\infty)^k\times\R^{n-k}$; hence $i_Q$ maps each depth stratum of $\partial Q$ locally diffeomorphically onto a depth stratum of $Q$.
\end{enumerate}

\begin{convention}[smoothness]\label{conv_smooth}
A map from a manifold with corners $P$ to $\R^K$, or to another manifold with corners $Q$, is \emph{smooth} if it is continuous and, in charts, locally the restriction of a smooth map between open subsets of Euclidean spaces. This is Joyce's \emph{weak} smoothness \cite[Definition~3.1]{Jo}. We impose neither Joyce's smoothness, which asks that a boundary defining function of the target pull back to $0$ near the point or to a boundary defining function of the source \cite[Definition~3.1]{Jo}, nor the weaker condition of a b-map, which allows a product $u\cdot\prod_ib_i^{e_i}$ of boundary defining functions $b_i$ of the source with exponents $e_i\ge0$ integers and $u>0$ smooth \cite[Definition~3.2]{Jo}: the simplices of the triangulations below need not be smooth in that sense, and even the simplices of linear triangulations need not be. For instance, an affine $2$-simplex in the half-plane $\{y\ge0\}$ with exactly one vertex on $\{y=0\}$ pulls $y$ back, up to positive scalars, to $\lambda_1+\lambda_2$, where $\lambda_1,\lambda_2$ are the barycentric coordinates of the two other vertices; near the vertex on $\{y=0\}$ the boundary defining functions of $\Delta^2$ are $\lambda_1,\lambda_2$, and $\lambda_1+\lambda_2$ is not of the form $u\lambda_1^a\lambda_2^b$ with $u>0$, whatever exponents are allowed.
\end{convention}

With this convention, a smooth map into $Q$ is smooth into $\R^K$ after composition with any smooth map $Q\to\R^K$, and a map into $\R^K$ with image in a closed embedded submanifold $X$ is smooth into $X$ if and only if it is smooth into $\R^K$. A \emph{smooth embedding} of manifolds with corners is a smooth injective map whose differential is injective at every point, corners included.

\begin{definition}\label{def_source}
A \emph{source} is a compact oriented smooth manifold with corners $Q$, possibly disconnected or empty, with \emph{embedded faces}: $i_Q$ is injective on each component of $\partial Q$. The class of sources is denoted $\cS$. A \emph{face} of $Q$ is a component $F$ of $\partial Q$, identified with $i_Q(F)$.
\end{definition}

By (J1), if $Q$ has embedded faces, the $k$ local boundary components at a point of depth $k$ belong to $k$ distinct faces. Simplices, products of sources and boundaries of sources are sources. For boundaries: a component of $\partial F$, $F$ a face of $Q$, corresponds along its length to one second face of $Q$, and two preimages of one point would give two local components of that face. Pieces of sources cut along regular levels are sources (Lemma~\ref{lem_pieces}).

Embeddedness of faces is part of the definition of sources and is used throughout, for instance in (G3), in Lemma~\ref{lem_pieces}(b) and in Section~\ref{sec_int}. For admissible triangulations it enters in two places: the proof of Lemma~\ref{lem_adm}(b), which uses that $i_Q|_F$ is a homeomorphism onto its image, and the reflection cover of Section~\ref{sec_point}. For non-embedded faces the conclusion $t_\partial\in\cT(\partial Q)$ of Lemma~\ref{lem_adm}(b) still holds, by lifting simplices instead of faces. A simplex $\rho$ of $L$ with $t(\mathring\rho)\subseteq\Bord_k(Q)$ has exactly $k$ lifts to $\partial Q$: the map $i_Q^{-1}(t(\mathring\rho))\to t(\mathring\rho)$ is a $k$-sheeted covering by (J1) and (J2), with $k$ sections over the open simplex. A section extends to the closed simplex, because the preimages of a point of $\partial\rho$ have neighbourhoods with compact closures on which $i_Q$ is injective and, $i_Q$ being proper, the image of the rest of $\partial Q$ stays away from that point. Distinct lifts are disjoint, and the lifts form an admissible triangulation of $\partial Q$. Given this, everything below holds for sources with non-embedded faces, provided the existence statements of Section~\ref{sec_point} hold for them. Topological triangulations exist in general by \cite{Yu2}; smooth ones are not known.

\subsection{Chains}
Let $Y\subseteq\R^K$ be any subset. A \emph{geometric simplex} in $Y$ is a pair $(Q,f)$ with $Q\in\cS$ and $f:Q\to\R^K$ smooth with $f(Q)\subseteq Y$. Isomorphisms are orientation-preserving diffeomorphisms of manifolds with corners commuting with the maps. The complex $\Cm_n(Y)$ is the $\Q$-span of isomorphism classes $[Q,f]$ modulo $[-Q,f]=-[Q,f]$ and $[Q_1\sqcup Q_2,f_1\sqcup f_2]=[Q_1,f_1]+[Q_2,f_2]$, with $\partial[Q,f]=[\partial Q,f\circ i_Q]$. By \cite[Lemma on the basis of $C^g_*$]{Yu1} this is a chain complex, $[Q,f]=0$ whenever $(Q,f)$ has an orientation-reversing automorphism, and the connected classes without such an automorphism, each with a chosen orientation, form a basis. If $W\subseteq Y$, then $\Cm(W)\subseteq\Cm(Y)$ is spanned by a subset of this basis.

A map $\psi$ that is defined and smooth on an open neighbourhood of $Y$ in $\R^K$, with $\psi(Y)\subseteq Y'\subseteq\R^{K'}$, induces $\psi_*:\Cm(Y)\to\Cm(Y')$. A smooth manifold $X$ is always regarded as a closed submanifold of some $\R^K$ (Whitney's embedding theorem \cite{Lee}); then $\Cm(X)$ is the complex of smooth maps into $X$ in the intrinsic sense.

$\Csm_n(Y)$ is the complex of singular simplices $\sigma:\Delta^n\to\R^K$ that are smooth (Convention~\ref{conv_smooth}) with image in $Y$, and $p(\sigma)=[\Delta^n,\sigma]$. With the face signs of \cite{Yu1}, $p:\Csm_*(Y)\to\Cm_*(Y)$ is a chain map. For a smooth manifold $X$, $\Csm_*(X)$ computes $H_*(X;\Q)$ \cite{Lee}; for compact polyhedra see Proposition~\ref{prop_amb}.

\subsection{Triangulations}\label{ssec_adm}
\begin{definition}\label{def_adm}
An \emph{admissible triangulation} of a source $Q$ is a homeomorphism $t:|L|\to Q$ from a finite simplicial complex such that
\begin{enumerate}[(i)]
\item every closed simplex of $L$, with an affine parametrisation, is mapped by $t$ smoothly into $Q$;
\item $t$ is \emph{depth-subordinate}: each open simplex is mapped into a single depth stratum.
\end{enumerate}
The simplices need not be immersions. The set of these is $\cT(Q)$. If moreover every closed simplex is mapped by a smooth embedding, $t$ is a \emph{Whitehead triangulation}; the set of these is $\cT^{C^1}(Q)\subseteq\cT(Q)$.
\end{definition}

Fix a vertex order. Each top simplex $\rho_i$ gives a smooth injective map $\sigma_i:\Delta^n\to Q$, with local degree $\varepsilon_i=\pm1$ at interior points. Put
\[
\mu_t=\sum_i\varepsilon_i\sigma_i,\qquad\tau(Q,f,t)=f_*\mu_t=\sum_i\varepsilon_i\,f\circ\sigma_i\in\Csm_n(Y).
\]

\begin{lemma}\label{lem_adm}
Let $t\in\cT(Q)$, $\dim Q=n$.
\begin{enumerate}[(a)]
\item Every simplex of $L$ is a face of an $n$-simplex; open $n$-simplices lie in $\Bord_0(Q)$; an $(n-1)$-simplex either has its interior in $\Bord_0(Q)$ and is a face of exactly two $n$-simplices, or has its interior in $\Bord_1(Q)$ and is a face of exactly one.
\item For each face $F$ of $Q$, $i_Q(F)$ is the image of a subcomplex, and $(i_Q|_F)^{-1}\circ t$ restricted to it is an admissible triangulation $t_F$ of $F$. Together these give $t_\partial\in\cT(\partial Q)$.
\item The open $n$-simplices are open in $Q$, disjoint from all other closed simplices, and the closed $n$-simplices cover $Q$.
\end{enumerate}
\end{lemma}

\begin{proof}
(a) A maximal simplex of dimension $<n$ would be open in $Q$, contradicting invariance of domain; an open $j$-cell in $\Bord_k$ forces $j\le n-k$. The count of $n$-simplices at an $(n-1)$-simplex is read off from $H_n(Q,Q\setminus x)$, which is $\Q$ at interior points and $0$ on $B_Q$.
(b) $i_Q(F)$ is closed. Along an open simplex in $\Bord_k$ the set of faces containing a point is constant, by (J2) and connectedness, so $i_Q(F)$ is a union of open simplices, hence a subcomplex. Since $F$ is embedded, $i_Q|_F:F\to i_Q(F)$ is a homeomorphism, and by (J2) its inverse is smooth wherever it is applied to smooth maps into $i_Q(F)$; depth drops by one under lifting (J1). Different faces give disjoint lifted complexes, so $t_\partial$ is admissible.
(c) is standard for triangulations of manifolds.
\end{proof}

For any topological space $Z$ let $\Calttop_*(Z)$ be the complex of continuous singular chains modulo $\sigma\circ a=\sgn(a)\sigma$ for affine permutations $a$ of $\Delta^n$. It has a basis given by the orbits of singular simplices without an orientation-reversing self-symmetry, so $\Calttop(W)\to\Calttop(Z)$ is injective for $W\subseteq Z$, and so is $\Calt(Y)\to\Calttop(Y)$ for the smooth version $\Calt(Y)$ of $Y\subseteq\R^K$ (Section~\ref{sec_R}). We regard $\mu_t$ as a chain in $\Calttop_n(Q)$.

\begin{lemma}[boundary formula]\label{lem_bdry}
In $\Calttop_{n-1}(Q)$, $\partial\mu_t=(i_Q)_*\mu_{t_\partial}$ for the induced vertex order. Consequently $\partial r(Q,f,t)=r(\partial Q,f\circ i_Q,t_\partial)$ for the relations $r(Q,f,t)$ of Definition~\ref{def_rel}.
\end{lemma}

\begin{proof}
Interior $(n-1)$-simplices occur twice in $\partial\mu_t$, as $\sigma_id_j$ and $\sigma_kd_l=\sigma_id_j\circ a$ with $a$ an affine permutation, with opposite signs $\varepsilon_i(-1)^j=-\varepsilon_k(-1)^l\sgn a$: this is the coherence of the orientation of a triangulated oriented manifold, since $\mu_t$ represents the class of local degree $1$ in $H_n(Q,B_Q)$. So they cancel in $\Calttop$. A boundary $(n-1)$-simplex occurs once, as $\sigma_id_j=i_Q\circ\tilde\rho$ with $\tilde\rho$ its unique lift (Lemma~\ref{lem_adm}(b)). The boundary map $H_n(Q,B_Q)\to H_{n-1}(B_Q)$ sends the class of local degree $1$ to the class of local degree $1$ at depth-one points for the outward-normal-first orientation; this is a local statement near a depth-one point and is checked in one chart. Hence $\varepsilon_i(-1)^j$ is the local degree of $\tilde\rho$ in $\partial Q$. The second claim follows by applying $p$ and using $\partial[Q,f]=[\partial Q,f\circ i_Q]$.
\end{proof}

Existence: $\cT(Q)\ne\emptyset$ (Corollary~\ref{thm_T7}) and $\cT^{C^1}(Q)\neq\emptyset$ (Theorem~\ref{thm_whitehead}) for every source $Q$.

\emph{Staircase products.} If $t\in\cT(Q)$, the staircase triangulations of the prisms $\rho\times I$ for a fixed vertex order form $t\times I\in\cT(Q\times I)$: its simplices are $(\sigma\times\id)\circ(\text{affine})$, and every open simplex lies in an open face times $\{0\}$, $(0,1)$ or $\{1\}$, so depth adds. If $t$ is Whitehead, so is $t\times I$.

\subsection{Gluing decompositions and the relations}
\begin{definition}[gluing decompositions]\label{def_glue}
Let $P\in\cS$, $\dim P=n$. A \emph{gluing decomposition} $\cQ=(Q_i,\iota_i)_{i=1}^m$ of $P$ is a finite family with the following properties.
\begin{enumerate}[(G0)]
\item Each $Q_i\in\cS$ has dimension $n$ (possibly disconnected).
\item Each $\iota_i:Q_i\to P$ is a smooth embedding of manifolds with corners and preserves orientation.
\item $\bigcup_i\iota_i(Q_i)=P$, and $\iota_i(\inter Q_i)\cap\iota_j(Q_j)=\emptyset$ for $i\ne j$.
\item (Aligned gluing along faces.) For every $i$ and every face $F$ of $Q_i$, either (a) $\iota_i(F)\subseteq B_P$, and $F$ is \emph{outer}; or (b) there are $j\ne i$ and a face $F'$ of $Q_j$ with $\iota_i(F)=\iota_j(F')$ such that $\psi:=(\iota_j|_{F'})^{-1}\circ\iota_i|_F:F\to F'$ is a diffeomorphism of manifolds with corners, and $F$ is \emph{glued} to its \emph{partner} $F'$.
\end{enumerate}
The pair $\fd=(P;\cQ)$ is a \emph{gluing datum}, and $\fG$ denotes the class of gluing data. The \emph{two-piece} gluing data are those with $m\le2$.
\end{definition}

Lower-dimensional contacts, such as two triangles sharing only a vertex, arise automatically and need no condition. By Lemma~\ref{lem_glue}(b) below, (a) and (b) of (G3) are exclusive, and the relation is automatically compatible with the maps, since $f\iota_j\psi=f\iota_i$ on $F$. Some of the readings and generalisations in Section~\ref{ssec_readings} use the following weaker notion.

\begin{definition}[decompositions]\label{def_dec}
A \emph{decomposition} $(Q_i,\iota_i)_i$ of $P\in\cS$, $\dim P=n$, is a finite family with
\begin{enumerate}[(D0)]
\item $Q_i\in\cS$ of dimension $n$;
\item $\iota_i:Q_i\to P$ smooth into the manifold with corners $P$ in the sense of Convention~\ref{conv_smooth}, and injective;
\item $\bigcup_i\iota_i(Q_i)=P$, and $\iota_i(\inter Q_i)\cap\iota_j(Q_j)=\emptyset$ for $i\ne j$;
\item $\iota_i$ has local degree $+1$ at interior points (Lemma~\ref{lem_dec}(a)).
\end{enumerate}
No immersivity and no condition on depth are imposed, and walls need not match. Such data $\fd=(P;(Q_i,\iota_i))$ are \emph{decomposition data}; their class is denoted $\cD$.
\end{definition}

\begin{definition}\label{def_rel}
For a datum $\fd=(P;(Q_i,\iota_i)_i)$ of either kind and $f:P\to\R^K$ smooth with $f(P)\subseteq Y$ put
\[
\rel{\fd}{f}:=[P,f]-\textstyle\sum_i[Q_i,f\circ\iota_i]\in\Cm_n(Y).
\]
Let $\cG(Y)$ and $\cD(Y)$ be the sets of these relations for $\fd\in\fG$ and $\fd\in\cD$, and
\[
G(Y):=\operatorname{span}\cG(Y)+\partial\operatorname{span}\cG(Y),\qquad\GCG(Y)=\Cm(Y)/G(Y),
\]
and likewise $E(Y)$ and $\GCE(Y)$ for $\cD$. For $t\in\cT(Q)$ put $r(Q,f,t)=[Q,f]-p\,\tau(Q,f,t)=[Q,f]-\sum_i\varepsilon_i[\Delta^n,f\sigma_i]$. For a smooth $\chi:Q\to\R$ with $0$ a regular value on every depth stratum (a \emph{bisecting function} \cite{Yu1}), the \emph{bisection datum} is $(Q;(Q_{\chi\le0},\mathrm{incl}),(Q_{\chi\ge0},\mathrm{incl}))$.
\end{definition}

$G$ and $E$ are subcomplexes by construction, and natural under the maps $\psi_*$. For $E$, $\operatorname{span}\cD$ is not closed under $\partial$, because walls need not match. For $G$ one expects $\partial\operatorname{span}\cG\subseteq\operatorname{span}\cG$: glued faces cancel by Lemma~\ref{lem_glue}(b), and the outer faces lying in a face $H$ of $P$ form a gluing decomposition of $H$. We do not need this and define $G$ as the generated subcomplex.

Each condition of (G1) and (G2) is needed. If pieces could overlap, then $[I,f]=2[I,f]$, so every $1$-chain would vanish and $H_1(S^1)$ would be lost, and integrals would not be additive; $H_0$ would survive, since the boundaries of the killed chains are differences of points. Without covering, the empty family would give $[P,f]=0$. Without orientation preservation one would get $[S^1,\id]=-[S^1,\id]$.

\begin{lemma}[decompositions]\label{lem_dec}
\begin{enumerate}[(a)]
\item Let $C$ be a component of a piece $Q_i$ of a decomposition. Then $\iota_i(\inter C)\subseteq\inter P$, the restriction of $\iota_i$ to $\inter C$ is an open topological embedding, and its local degree is a constant $\pm1$.
\item For $t\in\cT(Q)$, the family $(\varepsilon_i\Delta^n,\sigma_i)_i$ is a decomposition of $Q$, with relation $r(Q,f,t)$.
\item $E_0(Y)=0$.
\end{enumerate}
\end{lemma}

\begin{proof}
(a) Let $x\in\inter C$. If $\iota_i(x)$ had positive depth, then in a chart $U\cong\R^n_k\subseteq\R^n$ of $P$ the image of a small ball about $x$ would be open in $\R^n$ by invariance of domain, contained in $\R^n_k$ and containing a point of $\partial\R^n_k$, which is impossible. So $\iota_i|_{\inter C}$ is an injective map between $n$-manifolds without boundary, hence an open embedding, and its local degree is locally constant on the connected set $\inter C$.

(b) Each $\sigma_i$ is smooth into $Q$ and injective; the closed top simplices cover $Q$ and an open top simplex misses all other closed simplices (Lemma~\ref{lem_adm}(c)); the local degree of $\sigma_i$ on $\varepsilon_i\Delta^n$ is $+1$; and $[\varepsilon_i\Delta^n,f\sigma_i]=\varepsilon_i[\Delta^n,f\sigma_i]$.

(c) A decomposition of a $0$-manifold is a bijection of points preserving signs, so $\cD_0=\{0\}$. For a decomposition of a $1$-manifold the images are arcs, or whole circle components. An interior point of $P$ that is an endpoint image is the terminal point of exactly one piece and the initial point of exactly one other, by (D2) and (D3), and an endpoint of $P$ is the image of exactly one piece endpoint with the same sign; so $\partial\rel{\fd}{f}=0$, and $E_0=0$.
\end{proof}

\begin{lemma}[gluing decompositions]\label{lem_glue}
Let $\cQ=(Q_i,\iota_i)$ be a gluing decomposition of $P$.
\begin{enumerate}[(a)]
\item (Depth.) An immersion between $n$-manifolds with corners does not increase depth. In particular $\iota_i(\inter Q_i)\subseteq\inter P$, and an outer face $F$ satisfies $\iota_i(\inter F)\subseteq\Bord_1(P)$.
\item (Interfaces.) If $F$ is glued to $F'$, then $\iota_i(\inter F)\subseteq\inter P$, and $\psi$ reverses the boundary orientations; hence $[F,f\iota_i|]+[F',f\iota_j|]=0$. In particular (a) and (b) of (G3) are exclusive.
\item $\cQ$ is a decomposition, so $\fG\subseteq\cD$, $\cG(Y)\subseteq\cD(Y)$, $G\subseteq E$ and $G_0=0$.
\item Bisection data are two-piece gluing data.
\item For $t\in\cT^{C^1}(Q)$, the data of Lemma~\ref{lem_dec}(b) form a gluing decomposition, so $r(Q,f,t)\in\cG(Y)$.
\end{enumerate}
\end{lemma}

\begin{proof}
(a) Let $x$ have depth $k$ in $Q_i$ and $\iota_i(x)$ depth $m$ in $P$. The isomorphism $d\iota_i(x)$ maps the tangent cone at $x$, a copy of $\R^n_k$ whose lineality space has dimension $n-k$, into the tangent cone at $\iota_i(x)$, whose lineality space has dimension $n-m$. A linear subspace contained in a convex cone lies in its lineality space, so $n-k\le n-m$. A depth-one point of an outer face maps into $B_P$, hence to depth exactly one.

(b) Let $x\in\inter F$ and $y=\psi(x)$. Near $x$ and $y$, $\iota_i$ and $\iota_j$ are diffeomorphisms of half-balls onto their images, which are separated by the common hypersurface $\iota_i(F)$ and lie on opposite sides of it, because by (G2) the images of interiors are disjoint. If $\iota_i(x)$ had depth one in $P$, both images would have the same half-space as tangent cone, so they would lie on the same side; hence $\iota_i(x)\in\inter P$, the union of the two images being a neighbourhood. The outward normal of $Q_i$ at $x$ is mapped to an inward normal of $Q_j$ at $y$, so with outward normal first and both $\iota$ orientation preserving, $\psi$ reverses orientation. Since $f\iota_j\psi=f\iota_i$, the two classes cancel. Exclusivity follows from (a).

(c) (G0)--(G2) give (D0)--(D2), and orientation preservation gives (D3). $G_0\subseteq E_0=0$ by Lemma~\ref{lem_dec}(c).

(d) The inclusions $Q_\pm\hookrightarrow Q$ are embeddings, by the slice charts of Lemma~\ref{lem_pieces}(a) below. The faces of $Q_\pm$ are the components of $(\partial Q)_\pm$, which are outer, and the components of $\{\chi=0\}$; each of the latter is a face of both pieces, glued by $\psi=\id$.

(e) The pieces $(\varepsilon_i\Delta^n,\sigma_i)$ are embeddings, and (G2) holds by Lemma~\ref{lem_adm}(c). The faces of $\Delta^n$ are its facets. By Lemma~\ref{lem_adm}(a), a facet with interior in $\Bord_1(Q)$ is a facet of exactly one top simplex and is outer, and a facet with interior in $\Bord_0(Q)$ is a facet of exactly two top simplices, glued by $\psi=\sigma_j^{-1}\sigma_i$ on the facet, which is affine.
\end{proof}

\section{Relation families and the abstract theorem}\label{sec_axioms}

\begin{definition}\label{def_family}
A \emph{decomposition-type datum} is $\fd=(P;(Q_i,\iota_i)_i)$ with $P$ and the $Q_i$ sources of the same dimension $n$ and $\iota_i:Q_i\to P$ continuous. A \emph{relation family} is a class $\fF$ of pairs $(\fd,f)$ with $f:P\to\R^K$ smooth for some $K$ and every $f\circ\iota_i$ smooth. For $Y\subseteq\R^K$ put
\[
\cF(Y)=\{\rel{\fd}{f}:(\fd,f)\in\fF,\ f(P)\subseteq Y\},\qquad\cR(Y):=\operatorname{span}\cF(Y)+\partial\operatorname{span}\cF(Y),
\]
and $\GCcR(Y)=\Cm(Y)/\cR(Y)$. For a datum $\fd$, a smooth $h:P\to\R$ and $c\in\R$, the pair $(h,c)$ is \emph{regular for $\fd$} if $c$ is a regular value of $h$ on every depth stratum of $P$ and of $h\circ\iota_i$ on every depth stratum of every $Q_i$. Then the \emph{cut data} are
\[
\fd_\pm=(P_\pm;(Q_{i,\pm},\iota_i|)_i),\qquad\fd_0=(L_P;(L_{Q_i},\iota_i|)_i),
\]
for the pieces of $P$ for $h$ and of $Q_i$ for $h\circ\iota_i$ (Lemma~\ref{lem_pieces}), with $L_P$ and $L_{Q_i}$ oriented as faces of $P_-$ and $Q_{i,-}$ and empty pieces omitted; since the functions correspond, $\iota_i^{-1}(P_\pm)=Q_{i,\pm}$ and $\iota_i^{-1}(L_P)=L_{Q_i}$. The \emph{prism datum} is $\fd\times I=(P\times I;(Q_i\times I,\iota_i\times\id)_i)$.
\end{definition}

The examples are $\fG\times\{\text{all }f\}$, giving $\cG$ and $G$, and $\cD\times\{\text{all }f\}$, giving $\cD$ and $E$; others appear in Section~\ref{ssec_readings}. The axioms are:
\begin{enumerate}[(A1)]
\item \emph{Naturality and source data.} Membership in $\fF$ does not depend on the target: $\cF(Y)$ consists of all $\rel{\fd}{f}$ with $(\fd,f)\in\fF$ and $f(P)\subseteq Y$, as above; and $(\fd,f)\in\fF$ implies $(\fd,\psi\circ f)\in\fF$ for every $\psi$ smooth on an open neighbourhood of $f(P)$.
\item \emph{Cut relations.} For every source $Q$, every smooth $f:Q\to\R^K$, linear $\ell$ and $c$ regular for $\ell f$ on every depth stratum of $Q$, the bisection datum for $\ell f-c$, with the map $f$, is in $\fF$.
\item \emph{Cut closure.} If $(\fd,f)\in\fF$ and $(\ell f,c)$ is regular for $\fd$, then $(\fd_\pm,f|)$ and $(\fd_0,f|)$ are in $\fF$.
\item \emph{Prism closure.} If $(\fd,f)\in\fF$ and $H$ is smooth on an open neighbourhood of $f(P)\times I$, then $(\fd\times I,H\circ(f\times\id))\in\fF$.
\item \emph{Additivity of integrals.} If $(\fd,f)\in\fF$ and $\omega$ is a smooth $n$-form on an open neighbourhood of $f(P)$, then $\int_Pf^*\omega=\sum_i\int_{Q_i}(f\iota_i)^*\omega$.
\item \emph{The point.} $H_k(\GCcR(\mathrm{pt}))=0$ for all $k\ge1$. We write (A6)$_m$ for the condition in degrees $1\le k\le m$.
\end{enumerate}

By (A1), $\cF(Y)\subseteq\cF(Y')$ for $Y\subseteq Y'$, $\psi_*\cF(Y)\subseteq\cF(Y')$ for maps $\psi$ as in Section~\ref{sec_setting}, and a relation of $\cF(Y')$ whose map takes values in $Y\subseteq Y'$ lies in $\cF(Y)$; the last point is what the proof for manifolds needs. For $\fG$ and $\cD$, (A1) holds because the data do not involve $f$ at all.

\begin{theorem}[abstract comparison theorem]\label{thm_abstract}
Let $\fF$ satisfy \textup{(A1)--(A5)}, and let $Y$ be a compact polyhedron in $\R^K$ or a smooth manifold without boundary.
\begin{enumerate}[(a)]
\item $p_*:H_k(Y;\Q)\to H_k(\GCcR(Y))$ is injective for every $k$ and an isomorphism for $k=0$.
\item If \textup{(A6)$_m$} holds, then $p_*$ is an isomorphism in degrees $\le m$. In particular \textup{(A1)--(A6)} imply that $p_*$ is an isomorphism in all degrees.
\item Conversely, if $p_*$ is onto in degree $k\ge1$ for $Y=\mathrm{pt}$, then $H_k(\GCcR(\mathrm{pt}))=0$.
\end{enumerate}
\end{theorem}

Part (a) is Propositions~\ref{prop_stokes} and \ref{prop_poly_inj} together with chain-level surjectivity in degree $0$: a $0$-chain is a signed sum of points, i.e.\ of $p$ of singular $0$-simplices. Part (b) is Propositions~\ref{prop_abs_poly} and \ref{prop_abs_mfd}. Part (c) holds because $H_k(\mathrm{pt};\Q)=0$ for $k\ge1$. Not used anywhere in the proof: any triangulation beyond what (A6) itself requires, locality of the relations, and relative triangulations.

\section{Injectivity}\label{sec_inj}

\subsection{Smooth singular chains of polyhedra}
A \emph{compact polyhedron} is the underlying set of a finite simplicial complex linearly embedded in some $\R^K$.

\begin{proposition}\label{prop_amb}
Let $Y\subseteq\R^K$ be a compact polyhedron. Then $\Csm_*(Y)\hookrightarrow C_*(Y;\Q)$ is a quasi-isomorphism.
\end{proposition}

\begin{proof}
Fix a linear triangulation $L$ of $Y$. For a simplex $s$ of $L$ let $\St(s)$ be the union of the open simplices having $s$ as a face; it is open in $Y$, and $\St(s)\cap\St(s')$ is $\St(s\vee s')$ if $s\cup s'$ spans a simplex and empty otherwise. The vertex stars cover $Y$.

\emph{Small chains.} Precomposition with affine maps preserves smooth simplices, so $\Csm$ is preserved by barycentric subdivision and by the chain homotopy between it and the identity used in the proof of excision \cite{Hat}, since both are given by fixed linear chains in $\Delta^n$. Hence for every open cover $\cU$ of an open $W\subseteq Y$ the inclusion of $\cU$-small smooth chains into $\Csm(W)$ is a chain homotopy equivalence, as for $C(W)$.

\emph{Stars.} Let $V=\St(s)$ with barycentre $b$ of $s$. The homotopy $h(y,u)=(1-u)y+ub$ is smooth on $\R^K\times I$ and maps $V\times I$ into $V$: for $y$ in an open simplex $\rho\ge s$ and $u<1$, $h(y,u)$ is a convex combination with positive weight on $y\in\mathring\rho$ and on $b\in\bar\rho$, hence lies in $\mathring\rho$. Its prism operator preserves $\Csm(V)$ and contracts it onto the constant simplices. So $\Csm(V)$ and $C(V)$ both have the homology of a point.

\emph{Induction.} For open $A,B\subseteq Y$ the sequences $0\to\Csm(A\cap B)\to\Csm(A)\oplus\Csm(B)\to C^{\mathrm{sm},\{A,B\}}(A\cup B)\to0$ are exact, where $C^{\mathrm{sm},\{A,B\}}(A\cup B)\subseteq\Csm(A\cup B)$ is the subcomplex of $\{A,B\}$-small chains, spanned by the simplices with image in $A$ or in $B$; all three are spanned by subsets of the set of simplices. The same holds for $C$. By induction on $k$ and the five lemma, the inclusion is a quasi-isomorphism on every union of $k$ open stars: for $O=A\cup B$ with $A$ a union of $k-1$ stars and $B$ a star, $A\cap B$ is a union of at most $k-1$ stars. Apply this to the union of the vertex stars.
\end{proof}

The only point specific to the smooth class is that the contracting homotopy is smooth on all of $\R^K\times I$, which a merely continuous or piecewise linear homotopy would not give.

\subsection{Integration}
Let $U\subseteq\R^K$ be open and $\omega$ an $n$-form on $U$. For $[Q,f]\in\Cm_n(U)$ put $\int_{[Q,f]}\omega:=\int_Qf^*\omega$; $f^*\omega$ is a smooth form on $Q$ by Convention~\ref{conv_smooth}. This is invariant under isomorphism, changes sign with orientation, is additive under disjoint union, and vanishes on classes with an orientation-reversing automorphism, so it is a well-defined linear form on $\Cm_n(U)$. Stokes' theorem for compact manifolds with corners \cite{Lee}, with Joyce's boundary and the outward-normal-first orientation (the points of depth $\ge2$ have measure zero in $\partial Q$), gives
\begin{equation}\label{eq_stokes}
\int_{\partial c}\omega=\int_cd\omega\qquad(c\in\Cm_{n+1}(U)).
\end{equation}
On singular chains, $\int_{p(\sigma)}\omega=\int_\sigma\omega$ is ordinary integration over smooth simplices.

\begin{proposition}[Stokes]\label{prop_stokes}
Let $\fF$ satisfy \textup{(A1)} and \textup{(A5)}, and let $U\subseteq\R^K$ be open. Integration of any form vanishes on $\cR(U)$. If $X$ is a smooth manifold without boundary, $p_*:H_n(X;\Q)\to H_n(\GCcR(X))$ is injective for every $n$.
\end{proposition}

\begin{proof}
For $\rel{\fd}{f}\in\cF_n(U)$, (A5) says $\int_{\rel{\fd}{f}}\omega=0$. For $d'\in\cF_{n+1}(U)$, \eqref{eq_stokes} gives $\int_{\partial d'}\omega=\int_{d'}d\omega=0$.

On a smooth manifold $X\subseteq\R^K$ with a tubular neighbourhood $O$ and smooth retraction $\varrho:O\to X$, a form $\omega$ on $X$ gives $\varrho^*\omega$ on $O$, with $f^*\varrho^*\omega=f^*\omega$ for maps into $X$; so integration of forms on $X$ also vanishes on $\cR(X)$. Let $s\in\Csm_n(X)$ be a cycle with $p(s)=\partial w+e$, $e\in\cR(X)$. For closed $\omega$, $\int_s\omega=\int_wd\omega+\int_e\omega=0$. By de Rham's theorem for smooth singular chains \cite{Lee,dR}, integration identifies $H^n_{dR}(X)$ with $\operatorname{Hom}(H_n(X;\R),\R)$, so $[s]=0$ in $H_n(X;\R)$, hence in $H_n(X;\Q)$.
\end{proof}

\begin{proposition}[compact polyhedra]\label{prop_poly_inj}
Let $\fF$ satisfy \textup{(A1)} and \textup{(A5)}, and let $Y\subseteq\R^K$ be a compact polyhedron. Then $p_*:H_n(Y;\Q)\to H_n(\GCcR(Y))$ is injective for every $n$.
\end{proposition}

\begin{proof}
$Y$ is a compact Euclidean neighbourhood retract \cite[Theorem~A.7, Corollary~A.10]{Hat}: there are an open $U\supseteq Y$ and a continuous retraction $\varrho:U\to Y$.
Let $s\in\Csm_n(Y)$ be a cycle with $p(s)=\partial w+e$, $e\in\cR(Y)$. All chains involved lie in $\Cm(U)$, and $\cR(Y)\subseteq\cR(U)$ by (A1). By Proposition~\ref{prop_stokes} applied in $U$, $\int_s\omega=0$ for every closed $n$-form $\omega$ on $U$, so $[s]=0$ in $H_n(U;\Q)$ by de Rham's theorem on the manifold $U$. Then $[s]=\varrho_*i_*[s]=0$ in $H_n(Y;\Q)$, and by Proposition~\ref{prop_amb}, $s\in\partial\Csm_{n+1}(Y)$.
\end{proof}

\providecommand{\Lip}{\mathrm{Lip}}
\providecommand{\CmL}{C^{m}_{\Lip}}
\providecommand{\CL}{C^{\Lip}}
\providecommand{\GCGL}{GC^{G}_{\Lip}}

\subsection{Lipschitz maps}\label{ssec_lip}
The injectivity half of Theorem~\ref{thm_abstract} uses only (A1) and (A5), and both survive when the maps are only Lipschitz. The cutting axioms (A2) and (A3) are not available, since Lipschitz maps need not have regular levels (Remark~\ref{rem_lipschitz}).

\begin{convention}[Lipschitz chains]\label{conv_lip}
Let $Q$ be a source. A map $f:Q\to\R^K$ is \emph{Lipschitz} if it is Lipschitz, on each component, for the length metric of a Riemannian metric on $Q$. This does not depend on the metric. Diffeomorphisms preserve it, and if $f$ is Lipschitz and $\iota$ is smooth, then $f\circ\iota$ is Lipschitz.

For $Y\subseteq\R^K$, $\CmL(Y)$ is defined as $\Cm(Y)$, with Lipschitz maps in place of smooth ones. The basis lemma holds verbatim, and $\partial$ preserves $\CmL$. Likewise $\CL_*(Y)$ is the complex of Lipschitz singular simplices $\Delta^n\to\R^K$ with image in $Y$, and $p(\sigma)=[\Delta^n,\sigma]$. Put
\[
\cG_{\Lip}(Y)=\{\rel{\fd}{f}:\fd\in\fG,\ f\text{ Lipschitz},\ f(P)\subseteq Y\},
\]
let $G_{\Lip}(Y)=\operatorname{span}\cG_{\Lip}(Y)+\partial\operatorname{span}\cG_{\Lip}(Y)$, and let $\GCGL(Y)=\CmL(Y)/G_{\Lip}(Y)$. Define $G^W_{\Lip}$, $E_{\Lip}$ and $R_{\Delta,\Lip}$ in the same way.
\end{convention}

Let $(X,g)$ be a Riemannian manifold without boundary, regarded as a closed submanifold of some $\R^K$. On compact subsets of one component of $X$, the distance $d_g$ and the Euclidean distance are comparable: near the diagonal both are comparable to coordinate distances, and away from it both are bounded above and below. Hence a map from a compact connected source is Lipschitz into $(X,d_g)$ if and only if it is Lipschitz into $\R^K$, and $\CmL(X)$ and $\CL(X)$ do not depend on $g$ or on the embedding.

\begin{proposition}\label{prop_lip_sing}
The inclusion $\CL_*(Y)\hookrightarrow C_*(Y;\Q)$ is a quasi-isomorphism in each of the following cases:
\begin{enumerate}[(a)]
\item $Y$ is a compact polyhedron;
\item $Y$ is open in $\R^K$;
\item there are an open $U\supseteq Y$ and a locally Lipschitz retraction $U\to Y$; in particular, $Y$ is a closed smooth submanifold without boundary.
\end{enumerate}
For a smooth manifold $M$ without boundary, $\Csm_*(M)\hookrightarrow\CL_*(M)$ is a quasi-isomorphism.
\end{proposition}

\begin{proof}
(a) The proof of Proposition~\ref{prop_amb} applies with ``smooth'' replaced by ``Lipschitz''. Barycentric subdivision and the excision homotopy are precompositions with fixed affine chains. The contraction $h(y,u)=(1-u)y+ub$ of a star is Lipschitz on bounded subsets of $\R^K\times I$, so its prism operator preserves Lipschitz simplices.

(b) A cycle, together with a continuous chain that it bounds, has compact support. That support lies in a finite union $P\subseteq Y$ of closed cubes of a fine grid, and $P$ is a compact polyhedron, so (a) applies in $P$.

(c) Composition with the retraction $\varrho$ preserves Lipschitz simplices, so $\varrho_*$ is a chain retraction of $\CL(Y)\hookrightarrow\CL(U)$. Let $z$ be a continuous cycle in $Y$. By (b) it is homologous in $U$ to a Lipschitz cycle $z'$, so $z=\varrho_*z\sim\varrho_*z'$. Let $z$ be a Lipschitz cycle that bounds in $C(Y)$. By (b) it bounds some $c'\in\CL(U)$, and then $z=\partial\varrho_*c'$. For a closed submanifold, take a tubular retraction.

The last claim follows because both complexes are quasi-isomorphic to $C(M)$: for $\Csm$ by \cite{Lee}, for $\CL$ by (c), after a closed embedding $M\subseteq\R^K$.
\end{proof}

This is folklore; for metric spaces satisfying cone inequalities compare Riedweg and Sch\"appi \cite{RiSch}. The proof above is self-contained, and the citation is not used.

\subsubsection*{Integration}
By Rademacher's theorem \cite{Rad,EG}, a Lipschitz map $f:Q\to U\subseteq\R^K$ is differentiable at almost every point of $\inter Q$. So for a continuous $n$-form $\omega$ on $U$, $f^*\omega$ is an $L^\infty$ form on $Q$, defined almost everywhere. Composition with a diffeomorphism $\varphi$ onto an open set satisfies $(f\varphi)^*\omega=\varphi^*f^*\omega$ almost everywhere, since $\varphi$ preserves null sets. Hence $\int_{[Q,f]}\omega:=\int_Qf^*\omega$ is a well-defined linear form on $C^{m}_{\Lip,n}(U)$, by the argument given above for smooth maps.

\begin{lemma}[weak continuity of minors]\label{lem_lip_minors}
Let $\Omega\subseteq\R^m$ be open, and let $g_j,g:\Omega\to\R^K$ be Lipschitz with $\sup_j\operatorname{Lip}(g_j)<\infty$ and $g_j\to g$ locally uniformly. Then for every $I=(i_1,\dots,i_k)$, $dg_j^{i_1}\wedge\dots\wedge dg_j^{i_k}\to dg^{i_1}\wedge\dots\wedge dg^{i_k}$ weakly-$*$ in $L^\infty(\Omega)$.
\end{lemma}

\begin{proof}
First, for Lipschitz $h$, a smooth compactly supported $(m-k)$-form $\psi$, $k\ge1$ and $I'=(i_2,\dots,i_k)$,
\begin{equation}\label{eq_lip_ibp}
\int_\Omega dh^I\wedge\psi=(-1)^k\int_\Omega h^{i_1}\,dh^{I'}\wedge d\psi.
\end{equation}
For smooth $h$ this is Stokes' theorem for a compactly supported exact form. In general, apply it to $h*\rho_\delta$. Then $d(h*\rho_\delta)=(dh)*\rho_\delta\to dh$ almost everywhere with a uniform bound, because the almost-everywhere derivative of a Lipschitz map is its distributional derivative \cite{EG}, and dominated convergence applies.

Now induct on $k$; the case $k=0$ is trivial. Since the $dg_j^I$ are bounded in $L^\infty$, it suffices to test against smooth compactly supported $\psi$. By \eqref{eq_lip_ibp} for $g_j$,
\[
\int dg_j^I\wedge\psi=(-1)^k\int(g_j^{i_1}-g^{i_1})\,dg_j^{I'}\wedge d\psi+(-1)^k\int g^{i_1}\,dg_j^{I'}\wedge d\psi.
\]
The first term tends to $0$. The second tends to $(-1)^k\int g^{i_1}dg^{I'}\wedge d\psi$ by induction, since $g^{i_1}d\psi$ is an $L^1$ test form. That limit is $\int dg^I\wedge\psi$, by \eqref{eq_lip_ibp} for $g$.
\end{proof}

This is the classical weak continuity of Jacobian minors; compare Ball \cite{Ball}.

\begin{proposition}[Stokes for Lipschitz maps]\label{prop_lip_stokes}
Let $Q$ be a compact oriented $n$-manifold with corners, $U\subseteq\R^K$ open, $f:Q\to U$ Lipschitz and $\eta$ a smooth $(n-1)$-form on $U$. Then $\int_Qf^*d\eta=\int_{\partial Q}(f\circ i_Q)^*\eta$. Consequently \eqref{eq_stokes} holds for $c\in C^{m}_{\Lip,n+1}(U)$.
\end{proposition}

\begin{proof}
\emph{Approximation.} Choose finitely many charts $\varphi_a$ of $Q$ that are restrictions of diffeomorphisms onto $\R^n_{k_a}\cap(-3,3)^n$, with $Q=\bigcup_a\varphi_a^{-1}((-1,1)^n)$, and a subordinate smooth partition of unity $\chi_a$. On the convex set $\R^n_{k_a}\cap[-2,2]^n$ the map $f\varphi_a^{-1}$ is Lipschitz. Extend it to $(-2,2)^n$ by reflection, $(u,w)\mapsto f\varphi_a^{-1}(|u_1|,\dots,|u_{k_a}|,w)$, which is still Lipschitz, and mollify to get smooth maps $g_{a,j}$ on $(-\tfrac32,\tfrac32)^n$ that converge uniformly, with bounded Lipschitz constants. Then
\[
f_j=\sum_a\chi_a\cdot(g_{a,j}\circ\varphi_a)
\]
is smooth in the sense of Convention~\ref{conv_smooth}, and $f_j\to f$ uniformly. Its Lipschitz constants are bounded, since $\operatorname{Lip}(\chi h)\le\|d\chi\|\sup|h|+\operatorname{Lip}h$ for $h=g_{a,j}\varphi_a-f$ on the chart domain; these local bounds are global bounds because the metric on $Q$ is a length metric. For large $j$, $f_j(Q)\subseteq U$, because $f_j(x)$ is a convex combination of points near $f(x)$.

\emph{Limit.} By \eqref{eq_stokes}, $\int_Qf_j^*d\eta=\int_{\partial Q}(f_ji_Q)^*\eta$. In each chart, the integrand on the left is $\tilde\chi_a\sum_I(d\eta)_I(g_j)\,dg_j^I$, with $g_j=f_j\varphi_a^{-1}$ on $\Omega_a=\inter(\R^n_{k_a})\cap(-2,2)^n$, which contains the support of $\tilde\chi_a=\chi_a\varphi_a^{-1}$ up to a null set and on which the $g_j$ have bounded Lipschitz constants. Here $(d\eta)_I(g_j)\to(d\eta)_I(g)$ uniformly, and $\tilde\chi_a(d\eta)_I(g)\in L^1(\Omega_a)$. By Lemma~\ref{lem_lip_minors}, applied with the $L^1$ test function $\tilde\chi_a(d\eta)_I(g)$, the left side converges to $\int_Qf^*d\eta$, the boundary of the quadrant being null. Lemma~\ref{lem_lip_minors} is only needed on open subsets of Euclidean space: the corner strata of $Q$ and of $\partial Q$ are null. The same argument on the compact manifold with corners $\partial Q$ applies to the maps $f_j\circ i_Q\to f\circ i_Q$, which converge uniformly with bounded Lipschitz constants. It gives convergence of the right side.
\end{proof}

Mollification alone controls $df_j$ only almost everywhere on $Q$, which says nothing on the null set $\partial Q$. Lemma~\ref{lem_lip_minors}, applied on $\partial Q$ itself, supplies the boundary term. Alternatively one can use that push-forward of integral currents by Lipschitz maps commutes with $\partial$ \cite{FF,Fed}.

\begin{proposition}[(A5) for Lipschitz maps]\label{prop_lip_add}
Let $\fd=(P;(Q_i,\iota_i))\in\fG$, $f:P\to U$ Lipschitz and $\omega$ a continuous $n$-form on $U$. Then $\int_Pf^*\omega=\sum_i\int_{Q_i}(f\iota_i)^*\omega$.
\end{proposition}

\begin{proof}
By Lemma~\ref{lem_glue}(a), $\iota_i$ maps $\inter Q_i$ diffeomorphically, preserving orientation, onto an open subset of $\inter P$. So $\int_{Q_i}(f\iota_i)^*\omega=\int_{\iota_i(\inter Q_i)}f^*\omega$. By (G2) these open sets are disjoint. Their complement in $P$ lies in $\bigcup_i\iota_i(B_{Q_i})$, which is null.
\end{proof}

The same holds for decompositions $\fd\in\cD$, and hence for $G^W$. On the set where $d\iota_i$ is invertible the argument is the same. On the critical set $Z_i$ one needs $(f\iota_i)^*\omega=0$ almost everywhere, which follows from the area formula for the Lipschitz map $f\iota_i$ \cite{Fed,EG}: its image $f(\iota_i(Z_i))$ is $\mathcal H^n$-null, by Sard's theorem and because Lipschitz maps preserve $\mathcal H^n$-null sets.

\begin{theorem}[injectivity for Lipschitz maps]\label{thm_lip_inj}
Let $Y\subseteq\R^K$ be a neighbourhood retract in $\R^K$ for which $\CL(Y)\hookrightarrow C(Y;\Q)$ is a quasi-isomorphism; by Proposition~\ref{prop_lip_sing} this includes compact polyhedra and all Riemannian manifolds without boundary. Then
\[
p_*:H_n(Y;\Q)\cong H_n(\CL(Y))\longrightarrow H_n(\GCGL(Y))
\]
is injective for every $n$. The same holds for $G^W_{\Lip}$, $E_{\Lip}$ and $R_{\Delta,\Lip}$.
\end{theorem}

\begin{proof}
Let $\varrho:U\to Y$ be a continuous retraction, $U\supseteq Y$ open. Let $s\in\CL_n(Y)$ be a cycle with $p(s)=\partial w+e$, where $w\in C^{m}_{\Lip,n+1}(Y)$ and $e=\sum_l\lambda_ld_l+\partial\sum_m\mu_md'_m$ with $d_l,d'_m\in\cG_{\Lip}(Y)\subseteq\cG_{\Lip}(U)$. For a closed $n$-form $\omega$ on $U$, Propositions~\ref{prop_lip_stokes} and \ref{prop_lip_add} give $\int_{\partial w}\omega=\int_wd\omega=0$, $\int_{d_l}\omega=0$ and $\int_{\partial d'_m}\omega=\int_{d'_m}d\omega=0$. Hence $\int_s\omega=0$.

By Proposition~\ref{prop_lip_sing}, $s=s'+\partial b$ with $s'\in\Csm_n(U)$ a cycle and $b\in\CL_{n+1}(U)$. By Proposition~\ref{prop_lip_stokes} on $\Delta^{n+1}$, $\int_{s'}\omega=\int_s\omega=0$. By de Rham's theorem \cite{Lee,dR}, $[s]=[s']=0$ in $H_n(U;\R)$, hence in $H_n(U;\Q)$. Therefore $[s]=\varrho_*[s]=0$ in $H_n(Y;\Q)$, where only the continuity of $\varrho$ is used, and $s\in\partial\CL_{n+1}(Y)$ by hypothesis.

For $G^W$ and $E$, use the remark after Proposition~\ref{prop_lip_add}. For $R_\Delta$, use $R_\Delta\subseteq G^W$, which holds with Lipschitz maps by the proof of Theorem~\ref{thm_readings}.
\end{proof}

\begin{proposition}\label{prop_lip_low}
Let $Y$ be as in Theorem~\ref{thm_lip_inj}.
\begin{enumerate}[(a)]
\item Axioms (A1), (A4), (A5) and (A6) hold in the Lipschitz class. Hence Lemmas~\ref{lem_prism} and \ref{lem_star} hold for locally Lipschitz maps $\psi$ and homotopies $H$, since gluing data do not involve the maps. In particular $H_k(\GCGL(V))=0$ for $k\ge1$ and every star-shaped $V$.
\item Let $\Calttop_{\Lip}(Y)$ be $\CL(Y)$ modulo $\sigma\circ a=\sgn(a)\sigma$ for affine permutations $a$ of $\Delta^n$. The map $\bar p:\Calttop_{\Lip}(Y)\to\GCGL(Y)$ from these alternating Lipschitz chains is onto at chain level. Let $K_G(Y)$ be its kernel. Then
\[
H_n(\GCGL(Y))\cong H_n(Y;\Q)\oplus H_{n-1}(K_G(Y))\qquad\text{for all }n,
\]
and $p_*$ is an isomorphism in degrees $\le1$.
\item The natural map $H_*(\GCG(Y))\to H_*(\GCGL(Y))$ is injective when $Y$ is a compact polyhedron or a smooth manifold without boundary. It is onto in degree $n$ if and only if $p_*$ is onto in degree $n$ for $\GCGL(Y)$.
\end{enumerate}
\end{proposition}

\begin{proof}
(a) Gluing data do not involve the map. Prisms of Lipschitz maps are Lipschitz. Use Proposition~\ref{prop_lip_add}. Maps to a point are constant, so $\GCGL(\mathrm{pt})=\GCG(\mathrm{pt})$, and Corollary~\ref{cor_point} applies.

(b) By Theorem~\ref{thm_whitehead} and Lemma~\ref{lem_glue}(e), $[Q,f]\equiv\sum_i\varepsilon_i[\Delta^n,f\sigma_i]$ modulo $\cG_{\Lip}$, and each $f\sigma_i$ is Lipschitz. By acyclic models \cite{EM}, $\CL\to\Calttop_{\Lip}$ is a quasi-isomorphism, as in the smooth case, since Lipschitz simplices are closed under affine precomposition. The long exact sequence of $0\to K_G\to\Calttop_{\Lip}\to\GCGL\to0$ splits by Theorem~\ref{thm_lip_inj}. Finally $K_{G,0}=0$, since $p$ is injective on $0$-chains and $G_{\Lip,0}=0$ by the proof of Lemma~\ref{lem_dec}(c), which does not involve the map.

(c) The composite of this map with $p_*$ for $\GCG$, an isomorphism by Theorem~\ref{thm_A_G}, is $p_*$ for $\GCGL$.
\end{proof}

\subsubsection*{Toward a second proof, and the support of kernel cycles}
\begin{lemma}\label{lem_lip_goodpair}
Let $P$ be a compact topological $n$-manifold with boundary, and let $V_1,\dots,V_m\subseteq P$ be compact topological $n$-manifolds with boundary, such that $\bigcup_iV_i=P$ and $\inter V_i\cap V_j=\emptyset$ for $i\ne j$. Put $W=\bigcup_i\partial V_i$. Then $W\cap V_i=\partial V_i$, $\partial P\subseteq W$, $(P,W)$ is a good pair, and the inclusions induce an isomorphism $\bigoplus_iH_n(V_i,\partial V_i)\cong H_n(P,W)$. In particular a class in $H_n(P,W)$ vanishes if and only if its local degrees at interior points of the $V_i$ vanish.
\end{lemma}

\begin{proof}
A point of $V_i\cap V_j$, $i\ne j$, lies neither in $\inter V_i$ nor in $\inter V_j$. Hence $V_i\cap V_j\subseteq\partial V_i\cap\partial V_j$ and $W\cap V_i=\partial V_i$. The sets $\inter V_i$ are open in $P$ by invariance of domain, so they miss $\partial P$.

If $W=\emptyset$, the $V_i$ are unions of components of $P$ and everything is clear; so let $W\ne\emptyset$, and read $V_i/\emptyset$ as $V_i$ with a disjoint base point. Choose collars $c_i:\partial V_i\times[0,1)\to V_i$ \cite{Bro}; for $V_i=\iota(Q)$ with $Q$ a manifold with corners one can take the image of a collar of $B_Q$. The set $N=\bigcup_ic_i(\partial V_i\times[0,\frac12))$ satisfies $N\cap V_i=c_i(\partial V_i\times[0,\frac12))$, so it is open in $P$. The deformations $c_i(x,s)\mapsto c_i(x,(1-t)s)$ agree on overlaps, where they are the identity, so they glue to a deformation retraction of $N$ onto $W$. Moreover $\bigsqcup_iV_i\to P/W$ induces a continuous bijection $\bigvee_iV_i/\partial V_i\to P/W$ from a compact space onto a Hausdorff space, hence a homeomorphism. So $H_n(P,W)\cong\tilde H_n(P/W)\cong\bigoplus_iH_n(V_i,\partial V_i)$ \cite{Hat}, naturally. For $x\in\inter V_i$ the map $H_n(V_j,\partial V_j)\to H_n(P,P\setminus x)$ is zero if $j\ne i$, and for $j=i$ it factors through the excision isomorphism $H_n(V_i,V_i\setminus x)\cong H_n(P,P\setminus x)$.
\end{proof}

No tameness of $W$ is used: pairwise intersections $\partial V_i\cap\partial V_j$ may be arbitrary compacta, and $H_*(W)$ may be large, but $H_n(P,W)$ only sees the quotient $P/W$. No orientability is needed either: for non-orientable components $H_n(V_i,\partial V_i;\Q)=0$.

\begin{remark}[the dimension argument for $G$]\label{rem_lip_nerve}
With Lemma~\ref{lem_lip_goodpair}, the argument sketched in Remark~\ref{rem_szpilrajn} for $R_\Delta$ would give a second proof of Theorem~\ref{thm_lip_inj} for $G$, using no differential forms; we only sketch it.

Fix topological admissible triangulations of all sources (Corollary~\ref{thm_T7}). Define $\Theta$ as in the proof of Proposition~\ref{prop_structure}, with values in continuous alternating chains modulo differences $\tau(Q,f,t)-\tau(Q,f,t')$ of topological triangulations. Let $d=[P,f]-\sum_i[Q_i,f\iota_i]$ be a gluing relation in any of the readings $G$, $G^W$, $G^T$ or $E$. Then $\Theta(d)$ is represented by $f_*\nu$, where $\nu=\mu_{t_P}-\sum_i\iota_{i*}\mu_{t_{Q_i}}$ is a relative cycle of $(P,W_d)$ with $W_d=\bigcup_i\iota_i(B_{Q_i})$. Its local degrees vanish, so by Lemma~\ref{lem_lip_goodpair} we have $\nu=\partial\beta+\gamma$ with $\gamma$ supported in $W_d$.

Applying $\Theta$ to $p(s)=\partial w+e$, one finds that $s$ is homologous to a cycle supported in a finite union $A$ of the sets $(f\iota_i)(B_{Q_i})$ and $f_j(B_{Q_j})$. The relations of degree $n+1$ contribute only boundaries. The set $A$ is a Lipschitz image of a compact $(n-1)$-manifold, so $\mathcal H^n(A)=0$ and $\dim A\le n-1$ by Szpilrajn's theorem \cite{Szp,HW}. A nerve of a fine cover of $A$ of order $\le n-1$, mapped back through a neighbourhood retraction, kills the class.

This would prove injectivity for Riemannian manifolds without boundary, with their own metric, and for compact polyhedra, once the support and nerve steps, also used for $R_\Delta$ in Remark~\ref{rem_szpilrajn}, are written out; we do not do so here, and Theorem~\ref{thm_lip_inj} does not depend on this remark.
\end{remark}

\begin{remark}[what is open for Lipschitz maps]\label{rem_lip_open}
By Proposition~\ref{prop_lip_low}(b), surjectivity of $p_*$ in degree $n$ for $\GCGL(Y)$ is equivalent to $H_{n-1}(K_G(Y))=0$, and, for compact polyhedra and smooth manifolds without boundary, by (c) it is equivalent to every class being represented by $p$ of a smooth singular cycle.

In degree $2$ we expect $H_1(K_G)=0$, so that $p_*$ is an isomorphism there. After subdivision, the pieces of a $1$-dimensional gluing relation assemble to a Whitehead triangulation, and this reduces the claim to the interpolation of triangulations of $Q\times I$ for $\dim Q=1$ in the class $\cT^{C^1}$. We have not written out that interpolation statement, so we do not claim the degree-$2$ isomorphism.
In degrees $\ge3$ surjectivity is open. The proof of Section~\ref{sec_proof} needs (A2) and (A3), i.e.\ cuts along regular levels, and Lipschitz maps need not have any level that cuts a source into manifolds with corners (Remark~\ref{rem_lipschitz}). Smoothing generator by generator does not help either, because it is not natural under the embeddings $\iota_i$ and $i_Q$ along which cancellations take place.
\end{remark}

\begin{remark}[continuous maps]\label{rem_lip_cont}
Neither argument extends to continuous maps.

The integration proof fails at the start. There is no pull-back $f^*\omega$, and without a bound on Lipschitz constants Lemma~\ref{lem_lip_minors} fails, so integration does not extend to continuous geometric chains by continuity. (A linear form extending integration, killing the relations and satisfying Stokes' theorem would prove injectivity; we do not know whether one exists.)

In the dimension argument of Remark~\ref{rem_lip_nerve}, everything up to the support statement is topological. The sketched support step would apply verbatim to continuous maps and topological triangulations, giving that every kernel cycle is homologous to a cycle supported in a finite union $A$ of continuous images of compact $(n-1)$-manifolds. The argument fails exactly at the estimate $\dim A\le n-1$, since continuous images of $(n-1)$-manifolds can have any dimension. For $R_{\rm top}$ (Section~\ref{sec_cont}) the analogous question is open.
\end{remark}

\section{Cutting and Mayer--Vietoris}\label{sec_cut}

\subsection{Pieces}
\begin{lemma}\label{lem_pieces}
Let $Q\in\cS$, $h:Q\to\R$ smooth, and $c$ a regular value of $h$ on every depth stratum of $Q$ (on $0$-dimensional strata this means $h\ne c$). Then $Q_-=\{h\le c\}$, $Q_+=\{h\ge c\}$ and $L_Q=\{h=c\}$ are sources (possibly empty), with the following properties.
\begin{enumerate}[(a)]
\item Near $y\in L_Q$ of depth $k$ there is a slice chart $(u_1,\dots,u_k,h-c,w')$ of $Q$, in which $Q_-$ becomes $\R^n_{k+1}$ and $L_Q$ becomes a boundary hyperplane. The depth of $y$ is $k+1$ in $Q_\pm$ and $k$ in $L_Q$.
\item As oriented manifolds with corners, $\partial(Q_-)=(\partial Q)_-\sqcup L_Q$ and $\partial(Q_+)=(\partial Q)_+\sqcup(-L_Q)$, where $(\partial Q)_\pm$ are the pieces of $\partial Q$ for $h\circ i_Q$, and $L_Q$ is oriented as a face of $Q_-$ (outward normal first). The faces of $Q_\pm$ are the components of $F_\pm$ for faces $F$ of $Q$ and the components of $L_Q$; the faces of $L_Q$ are the components of $F\cap L_Q$.
\item $c$ is a regular value of $h\circ i_Q$ on every depth stratum of $\partial Q$.
\end{enumerate}
\end{lemma}

\begin{proof}
(a) is \cite[Lemma on bisection pieces]{Yu1} with bisecting function $h-c$: regularity on the stratum $\{u=0\}$ means $dh\notin\operatorname{span}(du_1,\dots,du_k)$, so $(u,h-c)$ extend to coordinates. (c) follows from (J2): $i_Q$ maps strata of $\partial Q$ locally diffeomorphically onto strata of $Q$. (b) The local boundary components of $Q_-$ at $x$ are those of $Q$ if $h(x)<c$, and those of $Q$ together with $\{h=c\}$ if $h(x)=c$; along the first kind the outward normal of $Q_-$ is that of $Q$. On $Q_+$ the outward normal along $L_Q$ is reversed. Faces: a face of $Q_-$ is a component of $(\partial Q)_-$, embedded because $i_Q$ is injective on the faces of $Q$, or a component of $L_Q$, which is embedded. The faces of $L_Q$ are the components of $(\partial Q)\cap\{h\circ i_Q=c\}$, again embedded.
\end{proof}

\subsection{Cut operators}
Fix a linear function $\ell$ on $\R^K$ and $c\in\R$. For $Y\subseteq\R^K$ put $Y_\pm=Y\cap\{\pm(\ell-c)\ge0\}$ and $Y_0=Y\cap\{\ell=c\}$. A geometric simplex $(Q,f)$ is \emph{$c$-transverse} if $c$ is a regular value of $\ell\circ f$ on every depth stratum of $Q$. A relation $\rel{\fd}{f}$ is $c$-transverse if $(\ell f,c)$ is regular for $\fd$, i.e.\ if $[P,f]$ and every $[Q_i,f\iota_i]$ are; these are independent conditions in general. Let $C^\pitchfork_c(Y)\subseteq\Cm(Y)$ be the span of the $c$-transverse basis classes. By Lemma~\ref{lem_pieces}(c) it is a subcomplex.

On $C^\pitchfork_c(Y)$ define
\[
\cut_\pm[Q,f]=[Q_\pm,f|]\in\Cm(Y_\pm),\qquad\Lambda[Q,f]=[L_Q,f|]\in\Cm(Y_0),
\]
where the pieces are taken for $h=\ell\circ f$ and $L_Q$ is oriented as a face of $Q_-$.

\begin{lemma}\label{lem_cutops}
$\cut_\pm$ and $\Lambda$ are well-defined linear maps on $C^\pitchfork_c(Y)$, and for every $c$-transverse generator, even one that is zero in $\Cm$, they are given by the formulas above. They satisfy
\begin{equation}\label{eq_cut}
\partial\,\cut_-=\cut_-\partial+\Lambda,\qquad\partial\,\cut_+=\cut_+\partial-\Lambda,\qquad\Lambda\partial=-\partial\Lambda.
\end{equation}
\end{lemma}

\begin{proof}
The formulas respect isomorphisms, orientation reversal and disjoint unions. If $a$ is an orientation-reversing automorphism of $(Q,f)$, then $a$ commutes with $f$, hence preserves $Q_\pm$ and $L_Q$, so both sides vanish. Hence the maps are well defined, and since $C^\pitchfork_c$ is spanned by a subset of the basis, any two expressions of an element of $C^\pitchfork_c$ through transverse generators have the same image. The first two identities are Lemma~\ref{lem_pieces}(b), applied in degrees $k$ and $k-1$. Then $0=\partial^2\cut_-=\partial(\cut_-\partial+\Lambda)=\Lambda\partial+\partial\Lambda$.
\end{proof}

\begin{lemma}[genericity]\label{lem_generic}
Given finitely many geometric simplices $(Q,f)$ in $Y$, the set of $c\in\R$ for which all of them are $c$-transverse has full measure.
\end{lemma}

\begin{proof}
Each depth stratum is a second countable manifold without boundary, and by Sard's theorem \cite{Sard} the critical values of $\ell\circ f$ on it form a null set; on $0$-dimensional strata these are finitely many values. There are finitely many strata, up to countably many components.
\end{proof}

A relation $\rel{\fd}{f}$ involves finitely many simplices, so it is $c$-transverse for almost all $c$, whatever the rank of the $\iota_i$: transversality of $[Q_i,f\iota_i]$ concerns the single function $\ell f\iota_i$. The configurations that one might fear, such as $\iota_i$ mapping a face of $Q_i$ onto an open part of $L_P$, or $\operatorname{im}d\iota_i\subseteq\ker d(\ell f)$ at a point of the level, make $c$ a critical value of $\ell f\iota_i$ and are therefore excluded. For a gluing decomposition, a glued face of $Q_i$ is a stratum of $Q_i$ mapped into $\inter P$; its condition is imposed as part of the transversality of $[Q_i,f\iota_i]$.

\begin{lemma}[cuts of relations]\label{lem_cutrel}
Let $\fF$ satisfy \textup{(A1)--(A3)}. If $\rel{\fd}{f}\in\cF(Y)$ is $c$-transverse, then
\[
\cut_\pm\rel{\fd}{f}=\rel{\fd_\pm}{f|}\in\cF(Y_\pm),\qquad\Lambda\rel{\fd}{f}=\rel{\fd_0}{f|}\in\cF(Y_0),
\]
and for every $c$-transverse $[Q,f]$ in $Y$, $[Q,f]-\cut_-[Q,f]-\cut_+[Q,f]\in\cF(Y)$.
\end{lemma}

\begin{proof}
By linearity (Lemma~\ref{lem_cutops}), $\cut_-\rel{\fd}{f}=[P_-,f]-\sum_i[Q_{i,-},f\iota_i]$, the pieces of $Q_i$ being taken for $\ell f\iota_i$; this is $\rel{\fd_-}{f|}$, and likewise for $+$ and $\Lambda$, with the orientations of Definition~\ref{def_family}. Membership is (A3), and the maps take values in $Y_\pm$, $Y_0$, so (A1) applies. The last claim is (A2).
\end{proof}

\subsection{Mayer--Vietoris in every degree}
\begin{lemma}[$\mathrm{MV}_k$]\label{lem_MV}
Let $\fF$ satisfy \textup{(A1)--(A3)} and \textup{(A5)}. Let $Y\subseteq\R^K$ be a compact polyhedron, $\ell$ a linear function, $\mathcal I\subseteq\R$ an open interval and $k\ge1$. Suppose that for all $c$ in a subset of $\mathcal I$ of full measure, $Y_\pm$ and $Y_0$ are compact polyhedra, $p_*$ is onto in degree $k$ for $\GCcR(Y_\pm)$, and $p_*$ is onto in degree $k-1$ for $\GCcR(Y_0)$. Then $p_*:H_k(Y;\Q)\to H_k(\GCcR(Y))$ is onto.
\end{lemma}

\begin{proof}
Let $x\in\Cm_k(Y)$ with $\partial x\in\cR(Y)$, and write
\[
\partial x=\sum_l\lambda_ld_l+\partial\Big(\sum_m\mu_md'_m\Big),\qquad d_l\in\cF_{k-1}(Y),\ d'_m\in\cF_k(Y).
\]
\emph{Step 1: choice of the level.} By Lemma~\ref{lem_generic} choose $c$ in the given set such that the generators of $x$, of the $d_l$ and of the $d'_m$ (sources and pieces) are $c$-transverse. Their faces are then transverse too, so every chain written so far lies in $C^\pitchfork_c(Y)$, and by Lemma~\ref{lem_cutops} the two expressions of $\partial x$ have the same images under $\cut_\pm$ and $\Lambda$.

\emph{Step 2: the boundary terms.} Put $\theta_\pm:=\cut_\pm(\partial x)$. By \eqref{eq_cut},
\[
\theta_-=\sum_l\lambda_l\cut_-(d_l)+\sum_m\mu_m\big(\partial\,\cut_-(d'_m)-\Lambda(d'_m)\big),
\]
and each term lies in $\cR(Y_-)$ by Lemma~\ref{lem_cutrel}, since $\cF(Y_0)\subseteq\cF(Y_-)$ by (A1). Likewise $\theta_+\in\cR(Y_+)$, and
\[
\Lambda(\partial x)=\sum_l\lambda_l\Lambda(d_l)-\sum_m\mu_m\partial\Lambda(d'_m)\in\cR(Y_0).
\]

\emph{Step 3: cutting $x$.} Put $x_\pm=\cut_\pm x\in\Cm_k(Y_\pm)$ and $z=\Lambda x\in\Cm_{k-1}(Y_0)$. By Lemma~\ref{lem_cutrel}, $x\equiv x_-+x_+$ modulo $\operatorname{span}\cF(Y)$. By \eqref{eq_cut},
\[
\partial x_-=\theta_-+z,\qquad\partial x_+=\theta_+-z,\qquad\partial z=-\Lambda(\partial x)\in\cR(Y_0).
\]

\emph{Step 4: the slice.} $z$ is a cycle of $\GCcR(Y_0)$, so by hypothesis $z=p(s_0)+\partial u+e_0$ with $s_0\in\Csm_{k-1}(Y_0)$ a cycle, $u\in\Cm_k(Y_0)$ and $e_0\in\cR_{k-1}(Y_0)$.

\emph{Step 5: injectivity on the halves.} In $\Cm(Y_-)$,
\[
p(s_0)=\partial(x_--u)-\theta_--e_0,\qquad p(-s_0)=\partial(x_++u)-\theta_++e_0\ \text{ in }\Cm(Y_+).
\]
So $p_*[s_0]=0$ in $H_{k-1}(\GCcR(Y_-))$ and $p_*[-s_0]=0$ in $H_{k-1}(\GCcR(Y_+))$. By Proposition~\ref{prop_poly_inj}, $s_0=\partial s_-$ and $-s_0=\partial s_+$ with $s_\pm\in\Csm_k(Y_\pm)$.

\emph{Step 6: cycles on the halves.} Put $A_-=x_--u-p(s_-)$ and $A_+=x_++u-p(s_+)$. Then
\[
\partial A_-=\theta_-+e_0\in\cR(Y_-),\qquad\partial A_+=\theta_+-e_0\in\cR(Y_+),
\]
and $s=s_-+s_+\in\Csm_k(Y)$ is a cycle.

\emph{Step 7.} $x\equiv x_-+x_+=A_-+A_++p(s)$. By hypothesis $A_\pm\equiv p(a_\pm)$ modulo $\cR(Y_\pm)+\partial\Cm_{k+1}(Y_\pm)$, with $a_\pm\in\Csm_k(Y_\pm)$ cycles. Pushing into $Y$, $x\equiv p(a_-+a_++s)$ modulo $\cR(Y)+\partial\Cm_{k+1}(Y)$.
\end{proof}

The chains $s_0,u,e_0,s_\pm,a_\pm$ introduced after Step 1 are never cut, so no further genericity is needed. No triangulation, no small-simplex argument and no locality statement enters. The only inputs are Sard's theorem, (A1)--(A3) and injectivity on the halves; the identity $x\equiv x_-+x_+$ is a relation, generator by generator.

\section{Proof of the abstract theorem}\label{sec_proof}

\subsection{Homotopies and star-shaped sets}
\begin{lemma}\label{lem_prism}
Let $\fF$ satisfy \textup{(A1)} and \textup{(A4)}. Let $H:U\times I\to\R^{K'}$ be smooth, $U\supseteq Y$ open, $H(Y\times I)\subseteq Y'$, and $\cH[Q,f]=(-1)^n[Q\times I,H\circ(f\times\id)]$ for $\dim Q=n$. Then $\partial\cH+\cH\partial=H_{1*}-H_{0*}$ and $\cH(\cR(Y))\subseteq\cR(Y')$.
\end{lemma}

\begin{proof}
The homotopy formula is \cite[Proposition on homotopy invariance]{Yu1}. All terms of a relation have the same dimension, so $\cH\rel{\fd}{f}=(-1)^n\rel{\fd\times I}{H\circ(f\times\id)}$, which lies in $\cF(Y')$ by (A4) and (A1). Then $\cH\partial d=H_{1*}d-H_{0*}d-\partial\cH d\in\cR(Y')$ for $d\in\cF(Y)$.
\end{proof}

Note that $\cH(\cR)\subseteq\cR$ needs the relation subcomplex, not just $\cF$: $\cH$ of a generated boundary $\partial d$ is handled by the homotopy formula.

\begin{lemma}\label{lem_star}
Let $\fF$ satisfy \textup{(A1)}, \textup{(A4)} and \textup{(A5)}.
\begin{enumerate}[(a)]
\item Let $\psi_0,\psi_1$ be maps defined near $Y$, with $\psi_i(Y)\subseteq Y'$, joined by a homotopy $H$ that is smooth on $U\times I$ for an open $U\supseteq Y$ and maps $Y\times I$ into $Y'$. Then $\psi_{0*}$ and $\psi_{1*}$ are chain homotopic as maps $\GCcR(Y)\to\GCcR(Y')$, compatibly with $p$.
\item If $V\subseteq\R^K$ is star-shaped about $v$ (every segment $[v,y]$, $y\in V$, lies in $V$) and \textup{(A6)$_m$} holds, then $H_k(\GCcR(V))=0$ for $1\le k\le m$, $H_0=\Q$, and $p_*$ is an isomorphism for $V$ in degrees $\le m$.
\end{enumerate}
\end{lemma}

\begin{proof}
(a) is Lemma~\ref{lem_prism}; compatibility with $p$ is naturality, $\psi_*\circ p=p\circ\psi_*$. Only this is used to transfer surjectivity of $p_*$ along homotopy equivalences. (b) $h(y,s)=v+(1-s)(y-v)$ is smooth on $\R^K\times I$ and preserves $V$, so $\id\simeq(\mathrm{const}_v)_*$, which factors through $\GCcR(\mathrm{pt})$, and the inclusion of the point and the constant map are mutually inverse homotopy equivalences. By (A6)$_m$, $H_k=0$ for $1\le k\le m$. In degree $0$, $p_*$ is onto at chain level and injective by Proposition~\ref{prop_stokes} applied to $0$-forms, and $V$ is contractible.
\end{proof}

Only the homology of the point enters, not chain-level vanishing.

\subsection{Compact polyhedra}
Let $\cA$ be the class of images of standard simplicial complexes $N\subseteq\R^V$ (vertices the basis vectors, simplices their convex hulls) under injective affine maps $A:\R^V\to\R^K$. For such $A$ choose an affine $B$ with $BA=\id$; by (A1), $A_*$ and $B_*$ are inverse isomorphisms between $\GCcR(N)$ and $\GCcR(A(N))$, compatible with $p$.

\begin{proposition}\label{prop_abs_poly}
Let $\fF$ satisfy \textup{(A1)--(A5)} and \textup{(A6)$_m$}. For every compact polyhedron $Y\subseteq\R^K$, $p_*:H_k(Y;\Q)\to H_k(\GCcR(Y))$ is an isomorphism for $k\le m$.
\end{proposition}

\begin{proof}
Injectivity is Proposition~\ref{prop_poly_inj}. For surjectivity let $S(k)$ be the statement that $p_*$ is onto in degree $k$ for every member of $\cA$. In degree $0$, $p$ is onto at chain level. Assume $S(k-1)$, $1\le k\le m$, and prove surjectivity in degree $k$ for a standard $N\subseteq\R^V$ by induction on the number of vertices; a point is Lemma~\ref{lem_star}(b). For $y\in|N|$ the \emph{carrier} of $y$ is the smallest simplex of $N$ containing $y$; its vertices are the $w$ with $y_w>0$. For a vertex $v$, the \emph{link} $\Lk v$ is the closed link: the subcomplex of the simplices of $N$ that do not contain $v$ but span a simplex of $N$ together with $v$; it is again a standard complex. Fix a vertex $v$, let $\ell=x_v$ and $\mathcal I=(0,1)$. For $c\in(0,1)$:
\begin{itemize}
\item $N_+$ is star-shaped about $e_v$: if $y_v>0$ then $v$ is a vertex of the carrier of $y$, and $x_v$ increases along $[y,e_v]$, which lies in the carrier. So Lemma~\ref{lem_star}(b) applies.
\item Let $N_{-v}$ be the full subcomplex on $V\setminus\{v\}$ and $\pi(y)=(y-y_ve_v)/(1-y_v)$, smooth on $\{x_v<1\}$. The linear homotopy from $\id$ to $j\circ\pi$ ($j$ the inclusion) stays in $N_-$, since $\pi(y)$ lies in the face of the carrier opposite to $v$. For a cycle $w$ of $\GCcR(N_-)$, Lemma~\ref{lem_star}(a) gives $w\simeq j_*\pi_*w$, and $\pi_*w\equiv p(s)$ by the inner induction; so $w\equiv p(js)$.
\item $N_0=ce_v+(1-c)|\Lk v|$, since a point with $y_v=c$ is $ce_v+(1-c)z$ with $z$ in the face of its carrier opposite to $v$; so $N_0$ is an injective affine image of the standard complex $\Lk v$, and $S(k-1)$ applies.
\end{itemize}
Lemma~\ref{lem_MV} gives surjectivity for $N$, hence for all of $\cA$.

For a general $Y$, take a linear triangulation $L$ of $Y$ with vertex set $V$ and its standard model $|L|_{\rm std}\subseteq\R^V$; the linear map $A:e_w\mapsto w$ maps $|L|_{\rm std}$ onto $Y$. For each vertex $w$ the set $W_w=\R^K\setminus(Y\setminus\St(w))$ is open with $W_w\cap Y=\St(w)$; choose smooth $\lambda_w\ge0$ with $\{\lambda_w>0\}=W_w$ \cite{Lee}. Then $\sum_w\lambda_w>0$ on a neighbourhood $O$ of $Y$, and $\tilde\kappa=(\lambda_w/\sum\lambda)_w:O\to\R^V$ is smooth. For $y\in Y$, $\lambda_w(y)>0$ exactly for the vertices of the carrier of $y$, so $\tilde\kappa(y)\in|L|_{\rm std}$ and $A\tilde\kappa(y)$ lies in the closed carrier of $y$; the linear homotopy from $\id_Y$ to $A\tilde\kappa$ stays in $Y$ and is smooth on $O\times I$. For a cycle $w$ of $\GCcR(Y)$, $w\simeq A_*\tilde\kappa_*w$, and $\tilde\kappa_*w\equiv p(s)$ because $|L|_{\rm std}\in\cA$; so $w\equiv p(As)$.
\end{proof}

\subsection{Smooth manifolds}
\begin{proposition}\label{prop_abs_mfd}
Let $\fF$ satisfy \textup{(A1)--(A5)} and \textup{(A6)$_m$}. For every smooth manifold $X$ without boundary, $p_*:H_k(X;\Q)\to H_k(\GCcR(X))$ is an isomorphism for $k\le m$.
\end{proposition}

\begin{proof}
Injectivity is Proposition~\ref{prop_stokes}. Regard $X\subseteq\R^K$ as a closed submanifold with a tubular neighbourhood $O$ and a smooth retraction $\varrho:O\to X$ \cite{Lee}. Let $x$ be a cycle of $\GCcR(X)$, with a chosen expression of $\partial x$ as an element of $\cR(X)$. The images of all sources involved form a compact set $S\subseteq X$; choose a finite union $P$ of closed cubes of a fine grid with $S\subseteq P\subseteq O$. The relations in the chosen expression have maps with values in $S\subseteq P$, so by (A1) they lie in $\cF(P)$, and $i_*x$ is a cycle of $\GCcR(P)$. So $i_*x\equiv p(s)$ modulo $\cR(P)+\partial\Cm(P)$ by Proposition~\ref{prop_abs_poly}. Since $\varrho$ is smooth on $O\supseteq P$, $\varrho(P)\subseteq X$ and $\varrho\circ i=\id_X$, applying $\varrho_*$ gives $x\equiv p(\varrho\circ s)$ modulo $\cR(X)+\partial\Cm(X)$.
\end{proof}

\subsection{Remarks}
\begin{remark}[boundary and corners]\label{rem_corners}
For compact manifolds with corners as targets see Proposition~\ref{prop_corners} in Section~\ref{sec_relprod}: there is no smooth retraction of a neighbourhood of $X$ in $\R^K$ onto $X$ if $\partial X\ne\emptyset$, so the proof of Proposition~\ref{prop_abs_mfd} does not apply directly, and the inward flow is used instead.
\end{remark}

\begin{remark}[no locality]\label{rem_locality}
In \cite{Yu1} the gluing relation of \cite{Yu0} is said to have made the comparison delicate, because cancellations could take place through relations involving chains outside a given open set. The argument above never compares $\GCcR(U)$ with its image in $\GCcR(Y)$: it works with explicit presentations, cuts them, and closes up through integration. A posteriori, for $U$ and $Y$ targets of Theorem~\ref{thmA}, e.g.\ $U$ open in a manifold $Y$, $H_*(\GCG(U))\to H_*(\GCG(Y))$ is the map of singular homology, but chain-level locality of $G$ is neither known nor needed.
\end{remark}

\begin{remark}[what the axioms do]\label{rem_axiomatic}
Axiom (A5) gives injectivity, (A1)--(A3) give the Mayer--Vietoris step, and (A4), (A6) give the star-shaped half. For the zero family the splitting (A2) fails as an identity; in \cite{Yu1} it is replaced by the bisection prism, and by \cite[Remark on the role of degenerate chains]{Yu1} the comparison theorem without a third relation then amounts to $H_n(\Cm_*(\mathrm{pt}))=0$ for $n>0$, which is open for manifold sources. The gluing relation provides (A2) and (A3) at chain level.
\end{remark}

\begin{remark}[bisections only]\label{rem_narrow}
The family of bisection data satisfies (A1)--(A5): a generic cut of a bisection datum is the bisection datum of the cut piece, as one checks with Lemma~\ref{lem_pieces} in the same way as in Lemma~\ref{lem_cutclosure}. So the comparison theorem for the relation generated by bisections alone is equivalent to the vanishing of the positive-degree homology of its complex of a point. This is an intermediate version of the question in \cite[Remark on the role of degenerate chains]{Yu1}, and it is open.
\end{remark}

\section{The gluing relation}\label{sec_G}

\subsection{Integration over decompositions}
\begin{lemma}\label{lem_cov}
Let $\iota:Q\to P$ be smooth and injective, $\dim Q=\dim P=n$, with local degree $+1$ at interior points, and let $\alpha$ be an $n$-form on $P$. Then $\int_Q\iota^*\alpha=\int_{\iota(Q)}\alpha$.
\end{lemma}

\begin{proof}
Let $O\subseteq\inter Q$ be the open set of points where $d\iota$ is invertible. By Lemma~\ref{lem_dec}(a), $\iota|_O$ is an injective local diffeomorphism into $\inter P$, hence a diffeomorphism onto the open set $\iota(O)$, orientation preserving because at regular points the sign of the Jacobian equals the local degree. On $\inter Q\setminus O$, $\iota^*\alpha=0$ pointwise, since an $n$-form pulled back through a singular differential vanishes. The set $\iota(Q)\setminus\iota(O)$ lies in $\iota(\inter Q\setminus O)\cup\iota(B_Q)$; the first is null by Sard's theorem \cite{Sard}, the second is a smooth image of an $(n-1)$-dimensional manifold with corners, hence null. So $\int_Q\iota^*\alpha=\int_O\iota^*\alpha=\int_{\iota(O)}\alpha=\int_{\iota(Q)}\alpha$.
\end{proof}

For an embedding the critical set is empty and only the last step is needed.

\begin{proposition}\label{prop_additive}
Axiom \textup{(A5)} holds for decompositions, hence for gluing decompositions: for $\fd\in\cD$ and every smooth $n$-form $\alpha$ on $P$, $\int_P\alpha=\sum_i\int_{Q_i}\iota_i^*\alpha$.
\end{proposition}

\begin{proof}
By (D2), $\iota_i(Q_i)\cap\iota_j(Q_j)\subseteq\iota_i(B_{Q_i})\cup\iota_j(B_{Q_j})$ for $i\ne j$, a null set, and the images cover $P$. By Lemma~\ref{lem_cov}, $\int_P\alpha=\sum_i\int_{\iota_i(Q_i)}\alpha=\sum_i\int_{Q_i}\iota_i^*\alpha$. Apply this to $\alpha=f^*\omega$.
\end{proof}

\subsection{Cut closure}
\begin{lemma}\label{lem_cutclosure}
Let $\fd=(P;(Q_i,\iota_i))$ be a decomposition datum and $(h,c)$ regular for $\fd$.
\begin{enumerate}[(i)]
\item $\fd_\pm$ are decomposition data;
\item $\fd_0$ is a decomposition datum;
\item if $\fd$ is a gluing datum, then so are $\fd_\pm$ and $\fd_0$, and if $\fd$ is two-piece, so are they.
\end{enumerate}
In particular \textup{(A3)} holds for $\cD$ and for $\fG$.
\end{lemma}

\begin{proof}
Since $Q_{i,\pm}$ and $L_{Q_i}$ are defined by the function $h\iota_i$, we have $\iota_i^{-1}(P_\pm)=Q_{i,\pm}$ and $\iota_i^{-1}(L_P)=L_{Q_i}$ as sets.

(i) (D0): $Q_{i,-}\in\cS$ by Lemma~\ref{lem_pieces}. (D1): $\iota_i|$ is injective; it is smooth into $P_-$, because near $L_P$ the slice-chart coordinates $(u,h-c,w')$ of $P_-$ are smooth functions on an open set of $P$, and away from $L_P$ the charts of $P_-$ are those of $P$. (D2): $\bigcup_i\iota_i(Q_{i,-})=\bigcup_i\iota_i(Q_i)\cap P_-=P_-$; and $\inter Q_{i,-}\subseteq\inter Q_i$, while $\iota_i(\inter Q_i)$ misses $\iota_j(Q_j)\supseteq\iota_j(Q_{j,-})$. (D3): interior points of $Q_{i,-}$ are interior points of $Q_i$, with the same orientation. The case $+$ is the same.

(ii) (D0): $L_{Q_i}\in\cS$. (D1): $\iota_i|$ is injective and smooth into $L_P$, since $L_P$ is a face of $P_-$ with slice coordinates $(u,w')$. (D2): covering follows from $\iota_i^{-1}(L_P)=L_{Q_i}$. Interior points of $L_{Q_i}$ have depth $0$ in $Q_i$ (Lemma~\ref{lem_pieces}(a)), so $\iota_i(\inter L_{Q_i})\subseteq\iota_i(\inter Q_i)$, which misses $\iota_j(L_{Q_j})$.

(D3) needs an argument, since $\iota_i$ may be singular along $L_{Q_i}$. Let $y\in\inter L_{Q_i}$. Then $y\in\inter Q_i$, and by Lemma~\ref{lem_dec}(a) $\iota_i$ maps a neighbourhood $U$ of $y$ homeomorphically onto an open set $V\subseteq\inter P$, with $\iota_i(U\cap Q_{i,-})=V\cap P_-$. So $\iota_i$ restricts to a homeomorphism of pairs $(U\cap Q_{i,-},U\cap L_{Q_i})\to(V\cap P_-,V\cap L_P)$ of topological manifolds with boundary. The boundary orientation of a topological manifold with boundary can be defined by the connecting homomorphism of local homology; this is natural under homeomorphisms of pairs, and on smooth manifolds it differs from the outward-normal-first orientation by a sign that depends only on the dimension. That sign is the same for $Q_{i,-}$ and $P_-$. Since $\iota_i$ preserves the orientation of the interiors (D3), it preserves the outward-normal-first orientations of $L_{Q_i}$ and $L_P$ at $y$. No immersivity is used.

(iii) Let $\fd$ be a gluing datum. (G0) is (D0). (G1): $d(\iota_i|)=d\iota_i$, so the restrictions to $Q_{i,\pm}$ are immersive; the restriction to the submanifold $L_{Q_i}$ is immersive with image tangent to $L_P$. For $\fd_0$ the orientation is now direct: at $y\in\inter L_{Q_i}$, $\iota_i$ is a local diffeomorphism near $y$ (Lemma~\ref{lem_glue}(a)) mapping $(Q_{i,-},L_{Q_i})$ to $(P_-,L_P)$ and the outward normal $\{h\iota_i>c\}$ to the outward normal $\{h>c\}$, and preserves orientation, so it preserves the outward-normal-first boundary orientations. (G2) is (D2).

(G3) for $\fd_-$. By Lemma~\ref{lem_pieces}(b), the faces of $Q_{i,-}$ are the components of $F_-=F\cap\{h\iota_i\le c\}$ for faces $F$ of $Q_i$, and the components of $L_{Q_i}$.
\begin{itemize}
\item If $F$ is outer, $\iota_i(F_-)\subseteq B_P\cap P_-\subseteq B_{P_-}$, so every component is outer.
\item If $F$ is glued to $F'$ by $\psi$, then $h\iota_j\psi=h\iota_i$, so $\psi$ maps $F_-$ onto $F'_-$ and components to components. The restriction $\psi|:F_-\to F'_-$ is a diffeomorphism of manifolds with corners, being the restriction of a diffeomorphism to corresponding sublevel sets of corresponding functions that are regular at $c$ on the strata (the strata of $F$ are strata of $Q_i$). Moreover $\iota_i(F_-)=\iota_i(F)\cap P_-=\iota_j(F')\cap P_-=\iota_j(F'_-)$. So each component is glued.
\item The components of $L_{Q_i}$ satisfy $\iota_i(L_{Q_i})\subseteq L_P\subseteq B_{P_-}$, so they are outer.
\end{itemize}
The case $+$ is the same.

(G3) for $\fd_0$. The faces of $L_{Q_i}$ are the components of $F\cap L_{Q_i}$. If $F$ is outer, $\iota_i(F\cap L_{Q_i})\subseteq B_P\cap L_P=B_{L_P}$, since depth in $L_P$ equals depth in $P$ (Lemma~\ref{lem_pieces}(a)). If $F$ is glued by $\psi$, then $\psi$ maps $F\cap L_{Q_i}=\{h\iota_i|_F=c\}$ diffeomorphically onto $F'\cap L_{Q_j}$, matching components, and $\iota_i(F\cap L_{Q_i})=\iota_j(F'\cap L_{Q_j})$.

The number of pieces does not increase. No T-junctions arise: every new face is outer or glued as a whole.
\end{proof}

What makes this work is that a cut is defined by pulling back one function: pieces, faces, partners and the diffeomorphisms $\psi$ are all cut by corresponding functions. Immersivity and alignment are preserved, not merely tolerated. For $\fd=r(Q,f,t)$ with $t\in\cT(Q)$, part (i) produces the curvilinear tiling of $Q_-$ by the sets $\sigma_i(\Delta^n_{i,-})$, and part (iii) shows that for a Whitehead triangulation it is a gluing decomposition; whether it lies in the simplex-only relation is the question of Section~\ref{sec_R}.

\begin{lemma}[prisms]\label{lem_products}
If $\fd$ is a decomposition datum, a gluing datum or a two-piece gluing datum, so is $\fd\times I$. In particular \textup{(A4)} holds for $\cD$ and for $\fG$.
\end{lemma}

\begin{proof}
$Q_i\times I\in\cS$, $\iota_i\times\id$ is smooth, injective and, for embeddings $\iota_i$, an orientation-preserving embedding for the product orientations; (D2), (D3) and (G2) are clear. The faces of $Q_i\times I$ are $F\times I$ and $Q_i\times\{0\}$, $Q_i\times\{1\}$. The latter map into $P\times\{0,1\}\subseteq B_{P\times I}$ and are outer; $F\times I$ is outer if $F$ is, and if $F$ is glued by $\psi$, then $F\times I$ is glued to $F'\times I$ by $\psi\times\id$.
\end{proof}

\subsection{The main theorem}
\begin{theorem}[Theorem~\ref{thmA}]\label{thm_A_G}
For every compact polyhedron $Y\subseteq\R^K$ and every smooth manifold $Y$ without boundary, $p_*:H_*(Y;\Q)\to H_*(\GCG(Y))$ is an isomorphism in all degrees.
\end{theorem}

\begin{proof}
Theorem~\ref{thm_abstract}, with $\fF=\fG\times\{\text{all }f\}$ and the axioms verified as follows: (A1) because gluing data do not involve the map; (A2) by Lemma~\ref{lem_glue}(d); (A3) by Lemma~\ref{lem_cutclosure}(iii); (A4) by Lemma~\ref{lem_products}; (A5) by Proposition~\ref{prop_additive} and Lemma~\ref{lem_glue}(c); (A6) by Corollary~\ref{cor_point}, which rests on Theorem~\ref{thm_whitehead}.
\end{proof}

\begin{corollary}\label{cor_dim}
$H_d(\GCG(N))=0$ for every compact polyhedron $N$ and every $d>\dim N$.
\end{corollary}

\subsection{Other readings, and the decomposition relation}\label{ssec_readings}
Definition~\ref{def_glue} fixes four points that ``gluing in the usual sense'' leaves open: the regularity of the pieces, disjointness and covering, what ``along faces'' means, and orientations. Two other natural readings are the following.
\begin{itemize}
\item[(W)] \emph{Injective-only.} In (G1) the $\iota_i$ are only smooth and injective with local degree $+1$, as in (D1) and (D3), and in (G3)(b) $\psi$ is only a homeomorphism. Class $\fG^W$, relation $G^W$.
\item[(T)] \emph{Topological, as in \cite{Yu0}.} The $\iota_i$ are topological embeddings with $f\circ\iota_i$ smooth, of local degree $+1$, with (G2) and (G3), $\psi$ a homeomorphism. Pairs $\fG^T$, relation $G^T$. Here membership depends on $f$, but not on the target, and (A1) holds.
\end{itemize}
Each of the readings has a two-piece version. Further variants also satisfy the axioms by the same arguments: self-gluings ($j=i$ allowed in (G3)(b), injectivity only on interiors), and interfaces that are closed unions of corner strata instead of whole faces glued pairwise.

\begin{theorem}\label{thm_readings}
The conclusion of Theorem~\ref{thm_A_G} holds for $G^W$ and for $E$. It holds for $G^T$, with injectivity from Lemma~\ref{lem_dual} instead of (A5) (Remark~\ref{rem_Z}(2)); Remark~\ref{rem_T} sketches (A5) for $G^T$ directly. We have $G\subseteq G^W\subseteq E$, $G^W\subseteq G^T$, and $R_\Delta\subseteq G^W$, with $R_\Delta$ as in Section~\ref{sec_R}.
\end{theorem}

\begin{proof}
$G^W$: (A1) and (A2) are clear, and (A5) is Proposition~\ref{prop_additive}, since $\fG^W\subseteq\cD$. For (A3), Lemma~\ref{lem_cutclosure}(i),(ii) gives (D0)--(D3) for the cut data, and the bookkeeping of (G3) in the proof of (iii) uses only that $\psi$ is a homeomorphism commuting with the cutting functions; exclusivity of outer and glued faces follows from invariance of domain applied to the double of a neighbourhood instead of Lemma~\ref{lem_glue}(a),(b). (A4) is Lemma~\ref{lem_products}. For (A6), Corollary~\ref{thm_T7} and the argument of Lemma~\ref{lem_glue}(e) show that for $t\in\cT(Q)$ the pieces of Lemma~\ref{lem_dec}(b) form a datum in $\fG^W$, so $R_\Delta\subseteq G^W$ and $[Q,c]\in G^W(\mathrm{pt})$ as in Corollary~\ref{cor_point}.

$E$: (A1), (A2) by Lemma~\ref{lem_glue}(c),(d); (A3) by Lemma~\ref{lem_cutclosure}(i),(ii); (A4) by Lemma~\ref{lem_products}; (A5) by Proposition~\ref{prop_additive}; (A6) by $G^W\subseteq E$.

$G^T$: a compact smooth injective map is a topological embedding, so $G^W\subseteq G^T$ and (A6) follows without Theorem~\ref{thm_whitehead}. (A1)--(A4) hold by the same arguments, the cutting functions $\ell f\iota_i$ being smooth. Injectivity, including injectivity on the halves in Step 5 of Lemma~\ref{lem_MV}, is Remark~\ref{rem_Z}(2), and replaces (A5), which Remark~\ref{rem_T} sketches.
\end{proof}

\begin{remark}[integration in the topological reading]\label{rem_T}
Let $\iota:Q\to P$ be a topological embedding of degree $+1$ with $f$ and $f\iota$ smooth. By the area formula, $\int_Pf^*\omega$ is an integral over $f(P)$ against $n$-dimensional Hausdorff measure with signed multiplicities. By Sard's theorem the images of the critical sets of $f$ and of $f\iota$ are null. For almost every $y$, all preimages are regular for both maps, and there $\iota=f^{-1}\circ(f\iota)$ is locally a diffeomorphism whose Jacobian sign is its topological degree $+1$. So the multiplicities agree, and (A5) follows. This is a sketch.
\end{remark}

\begin{corollary}\label{cor_GE}
The maps $H_*(\GCG(Y))\to H_*(GC^{G^W}(Y))\to H_*(\GCE(Y))$ induced by the inclusions are isomorphisms for compact polyhedra and smooth manifolds without boundary.
\end{corollary}

\begin{proof}
They commute with $p$, which is an isomorphism onto each of them.
\end{proof}

So walls that do not match and non-immersive pieces add nothing in homology. In the proofs above, walls that do not match are never used: every decomposition that actually arises (triangulations, bisections, their cuts and prisms and the interfaces of these) is aligned. Non-immersive pieces are used only at the point, through Corollary~\ref{thm_T7}.

\subsection{Two pieces}
\begin{proposition}\label{prop_twopiece}
Let $G_2\subseteq G$ be the relation generated by two-piece gluing data.
\begin{enumerate}[(a)]
\item $G_2$ satisfies \textup{(A1)--(A5)}, so $p_*:H_k(Y;\Q)\to H_k(GC^{G_2}(Y))$ is injective for every $k$, for $Y$ as in Theorem~\ref{thm_A_G}.
\item $p_*$ is an isomorphism for $k\le2$; for $k=2$ this uses, without proof, the classification of compact surfaces with corners: after rounding the corners, the diffeomorphism type is determined by the genus, the number of boundary circles and the number of corners on each.
\item For $m\ge3$, $p_*$ is an isomorphism in degrees $\le m$ for all such $Y$ if and only if $H_k(GC^{G_2}(\mathrm{pt}))=0$ for $3\le k\le m$.
\end{enumerate}
\end{proposition}

\begin{proof}
(a) Bisection data are two-piece, and cuts and prisms of two-piece data are two-piece (Lemmas~\ref{lem_cutclosure}(iii), \ref{lem_products}); the pieces may be disconnected, which (G0) allows. (A1) and (A5) are as for $G$. Apply Theorem~\ref{thm_abstract}(a).

(b) By Theorem~\ref{thm_abstract}(b) it suffices that $\Cm_k(\mathrm{pt})=0$ for $k=1,2$, i.e.\ that every connected source of dimension $1$ or $2$ has an orientation-reversing diffeomorphism. In dimension $1$ the sources are intervals and circles. In dimension $2$, take a model of the surface that is symmetric under a reflection of $\R^3$, with each boundary circle meeting the mirror in two points and the corners on it placed symmetrically (one corner on the mirror if their number is odd); by the classification of compact surfaces with corners, every source is diffeomorphic to such a model.

(c) Theorem~\ref{thm_abstract}(b),(c), together with (b).
\end{proof}

From dimension $3$ on, chiral closed manifolds such as lens spaces exist, so some relation is needed at the point, and Theorem~\ref{thm_whitehead} does not help directly: unions of adjacent curvilinear simplices are in general not manifolds with corners, because of re-entrant angles, so a triangulation relation is not obviously a combination of two-piece relations.

\begin{question}\label{q_twopiece}
Is $H_k(GC^{G_2}(\mathrm{pt}))=0$ for $k\ge3$? Equivalently, does the gluing relation of \cite{Yu0}, in its smooth two-piece form, give singular homology in all degrees?
\end{question}

\section{The point: Whitehead triangulations with corners}\label{sec_point}

\subsection{The reflection cover and its fundamental domain}
Let $Q$ be a source of dimension $n$, $F_1,\dots,F_m$ its faces, $\Gamma=\Z_2^m$ with basis $e_a$, and for $x\in Q$ let $\Gamma_x$ be generated by the $e_a$ with $x\in F_a$; by embeddedness $\Gamma_x\cong\Z_2^{\operatorname{depth}x}$. Put
\[
M=(\Gamma\times Q)/\!\sim,\qquad(g,x)\sim(g',x)\iff g-g'\in\Gamma_x.
\]
This is Davis's basic construction for the right-angled reflection group, abelianised \cite{Dav}. Near $x$ take a chart of $Q$ with coordinates $(y,l)\in\R^n_k$, in which $\{y_i=0\}$ lies in the face $F_{a_i}$, and let $s(\eta,l)=(\eta_1^2,\dots,\eta_k^2,l)$. The maps $\eta\mapsto[g+\sum_{\eta_i<0}e_{a_i},\,s(\eta)]$ are bijections onto open subsets of $M$, since the $e_{a_i}$ are independent. Corner-preserving chart transitions lift to smooth maps of the doubled charts, up to sign changes of the $\eta$-coordinates between charts with different base points $g$ \cite[Doubling lemma]{Yu1}, and the face labels match because a transition maps $\{y_i=0\}$ into the hyperplane lying in the same face. So $M$ is a closed smooth $n$-manifold. It is Hausdorff: points over distinct points of $Q$ are separated by preimages of disjoint open sets of $Q$; and if $g-g'\notin\Gamma_x$, the chart images $\{[g+h,y]:h\in\Gamma_x,\ y\in U\}$ and $\{[g'+h,y]:h\in\Gamma_x,\ y\in U\}$, for a small neighbourhood $U$ of $x$, are disjoint, since $\Gamma_y\subseteq\Gamma_x$ for $y\in U$. Further, $\Gamma$ acts smoothly; the projection $\pi:M\to Q$ is $s$ in charts, hence smooth into $Q$; and the isotropy of $[g,x]$ is $\Gamma_x$.

Let $\hat Q\subseteq M$ be the image of $\{0\}\times Q$; it is closed, as the image of a compact set.

\begin{lemma}\label{lem_Qhat}
$\hat Q$ is a compact manifold with corners whose charts are restrictions of charts of $M$, and $\pi|_{\hat Q}:\hat Q\to Q$ is a smooth homeomorphism that preserves depth. Equivalently, $\hat Q$ is the set $|Q|$ with the atlas of the charts $(\sqrt{y_1},\dots,\sqrt{y_k},l)$, for corner charts $(y,l)$ of $Q$ in which $y_1,\dots,y_k$ vanish at the base point. On $\inter\hat Q$ this atlas is that of $\inter Q$, and we orient $\hat Q$ accordingly.
\end{lemma}

\begin{proof}
Take a chart based at $g=0$; this is allowed because $[g,x]=[0,x]$ if and only if $g\in\Gamma_x$. In it, $[\sum_{\eta_i<0}e_{a_i},s(\eta)]=[0,s(\eta)]$ if and only if no $\eta_i<0$, because the $e_{a_i}$ are independent. So $\hat Q=\{\eta_1,\dots,\eta_k\ge0\}$ locally, and the transitions preserve the signs of the $\eta$. In these charts $\pi$ is $(\eta,l)\mapsto(\eta^2,l)$, a smooth homeomorphism of $\R^n_k$ preserving depth.

For the second description: a corner chart transition has the form $y'_{\sigma(i)}=y_ia_i(y,l)$ with $a_i>0$, by Hadamard's lemma and the invertibility of the Jacobian. Then $\sqrt{y'_{\sigma(i)}}=\eta_i\sqrt{a_i(\eta^2,l)}$ is smooth in $(\eta,l)$, and so is $l'(\eta^2,l)$; coordinates $y'_j$ that are positive at the point have smooth square roots. The same holds in the inverse direction.
\end{proof}

\begin{lemma}[triangulating the fundamental domain]\label{lem_fund}
$\cT^{C^1}(\hat Q)\ne\emptyset$.
\end{lemma}

\begin{proof}
By Illman's theorem \cite[Theorem, p.~199]{Ill} there is a $\Gamma$-equivariant smooth triangulation $h:|J|\to M$. Smooth in Illman's sense means that $h$ has rank equal to the dimension on every closed simplex, i.e.\ it is a Whitehead triangulation, and $J$ is an equivariant simplicial complex in Illman's sense: if $g$ fixes a point of an open simplex $\mathring s$, it fixes $s$ pointwise \cite[pp.~201, 207]{Ill}. So isotropy is constant on open simplices, and no further subdivision is needed.

An open top simplex is open in $M$ and has constant isotropy; the non-free set is the nowhere dense union of the mirrors, so the isotropy is trivial. The free set is $\bigsqcup_{\gamma\in\Gamma}\gamma\cdot\inter\hat Q$, and an open top simplex is connected, so it lies in a single sheet $\gamma\cdot\inter\hat Q$. Let $J_0$ consist of the top simplices lying in $\inter\hat Q$ and all their faces. Then $h(|J_0|)$ is closed, lies in $\hat Q$ and contains the open top simplices in $\inter\hat Q$, which are dense there; hence $h(|J_0|)=\hat Q$, and $\hat t:=h|_{|J_0|}:|J_0|\to\hat Q$ is a homeomorphism. Each closed simplex is a smooth embedding into $M$ with image in $\hat Q$, hence a smooth embedding into $\hat Q$, whose charts are restrictions of charts of $M$ (Lemma~\ref{lem_Qhat}). The depth of $x\in\hat Q$ is the number of mirrors through $x$, i.e.\ $\log_2|\Gamma_x|$, which is constant on open simplices. So $\hat t\in\cT^{C^1}(\hat Q)$.
\end{proof}

\begin{corollary}[smooth admissible triangulations]\label{thm_T7}
Every source $Q$ has an admissible triangulation: $\cT(Q)\ne\emptyset$.
\end{corollary}

\begin{proof}
$\pi\circ\hat t\in\cT(Q)$, since $\pi|_{\hat Q}$ is a smooth homeomorphism preserving depth (Lemma~\ref{lem_Qhat}).
\end{proof}

This uses only Illman's theorem. The top simplices of $\pi\circ\hat t$ that meet $B_Q$ touch the faces with vanishing normal derivative ($\partial y_i/\partial\eta_i=0$), so they are not immersions there; this is why Definition~\ref{def_adm} does not ask for immersions, and why the gluing relation needs the following theorem instead.

\subsection{Straightening by collars}
\begin{theorem}[Theorem~\ref{thmB}]\label{thm_whitehead}
Every compact manifold with corners $Q$ with embedded faces, orientable or not, has a Whitehead triangulation: $\cT^{C^1}(Q)\ne\emptyset$.
\end{theorem}

\begin{proof}
It suffices to construct an orientation-preserving diffeomorphism $\Phi:Q\to\hat Q$ of manifolds with corners; then $\Phi^{-1}\circ\hat t\in\cT^{C^1}(Q)$ for $\hat t$ as in Lemma~\ref{lem_fund}, because diffeomorphisms preserve depth and embeddings. Since $\hat Q=|Q|$ as sets, $\Phi$ is a self-map of $|Q|$.

\emph{Collars.} A source is a manifold with faces, in the sense that each point of depth $k$ lies in $k$ distinct faces $F_a$ and $F_a\cap F_b$ is a union of components of $\partial F_a$, the labels being locally constant. With the faces $F_1,\dots,F_m$ of $Q$ as the $m$ faces, $Q$ is then an $\langle m\rangle$-manifold in the sense of \cite[Section~2.1]{La}, and the $\langle m\rangle$-diagram of collar embeddings $\R^m_+(a')\times Q(a)\hookrightarrow\R^m_+(b')\times Q(b)$ of \cite[Lemma~2.1.6]{La}, restricted to $[0,\varepsilon_0)$ in each collar coordinate, gives
open embeddings $c_a:F_a\times[0,\varepsilon_0)\to Q$ with $c_a(x,0)=x$, such that $c_a$ maps $(F_a\cap F_b)\times[0,\varepsilon_0)$ into $F_b$, and $c_a(c_b(x,s),t)=c_b(c_a(x,t),s)$ for $x\in F_a\cap F_b$ wherever defined. Let $t_a$ and $\pr_a$ be the collar coordinate and the projection to $F_a$ on $U_a(\varepsilon_0):=c_a(F_a\times[0,\varepsilon_0))$. For a set $A$ of faces with $F_A:=\bigcap_{a\in A}F_a\neq\emptyset$, iterating the collars gives a map $C_A:F_A\times[0,\varepsilon_0)^A\to Q$, independent of the order by compatibility; it is an open embedding onto a neighbourhood $W_A$ of $F_A$, by induction on $|A|$: $c_a^{-1}(F_b)=(F_a\cap F_b)\times[0,\varepsilon_0)$, so $c_a$ restricts to a collar of $F_a\cap F_b$ in $F_b$, and $C_A$ is a composite of open embeddings of manifolds with corners. Putting $c_{a_j}$ outermost shows $t_{a_j}\circ C_A(z,t)=t_j$. So the $(z,t)$ are corner charts of $Q$ near $F_A$, $z$ ranging over charts of $F_A$.

\emph{Shrinking.} We build the collars on the wider interval $[0,\varepsilon_0)$ and work with $U_a(\varepsilon):=\{t_a<\varepsilon\}\subseteq U_a(\varepsilon_0)$ for $\varepsilon\ll\varepsilon_0$. By compactness, $\varepsilon$ can be chosen so that for every set $A$ of faces, $\bigcap_{a\in A}\overline{U_a(\varepsilon)}\subseteq W_A$ if $F_A\ne\emptyset$, and $\bigcap_{a\in A}\overline{U_a(\varepsilon)}=\emptyset$ if $F_A=\emptyset$: the sets $\overline{U_a(\varepsilon)}$ shrink to $F_a$ as $\varepsilon\to0$, and there are finitely many $A$.

\emph{The maps $R_a$.} Fix a homeomorphism $\theta:[0,\infty)\to[0,\infty)$ with $\theta(t)=t^2$ near $0$, $\theta=\id$ on $[2\varepsilon/3,\infty)$, and $\theta$ smooth with $\theta'>0$ on $(0,\infty)$. Put $R_a(p)=c_a(\pr_a(p),\theta(t_a(p)))$ on $U_a(\varepsilon)$ and $R_a=\id$ elsewhere. This is a homeomorphism of $|Q|$, continuous because $\theta=\id$ near the far end of the collar. For $b\ne a$, $R_a$ maps $U_b(\varepsilon)$ onto itself and preserves $t_b$: if $p\in U_a(\varepsilon)\cap U_b(\varepsilon)$, then $p\in W_{\{a,b\}}$, where in the coordinates $(z,t_a,t_b)$ the map $R_a$ changes only $t_a$; and if $p\in U_a(\varepsilon)\setminus U_b(\varepsilon)$, the same argument applied to $R_a^{-1}$ shows $R_a(p)\notin U_b(\varepsilon)$. Hence the $R_a$ commute, and $\Phi:=\prod_aR_a$ is well defined, independently of the order.

\emph{Local form.} Let $p\in Q$ and $A=A(p)$ the set of $a$ with $p\in\overline{U_a(\varepsilon)}$; if $A=\emptyset$, $\Phi=\id$ near $p$. Otherwise $p\in W_A$, so $p=C_A(z,t)$, and $z$ lies in $\inter F_A$: $F_A$ is a manifold with corners whose boundary points lie in faces $F_b$ with $b\notin A$, since the depth in $Q$ of a point of $F_A$ is $|A|$ plus its depth in $F_A$; and if $z$ lay in a face $F_b$ with $b\notin A$, then $p\in F_b\subseteq U_b(\varepsilon)$ by compatibility, a contradiction. Faces not in $A$ have closed collars that miss a neighbourhood of $p$. So near $p$, in the corner chart $(z,t)$ of $Q$,
\[
\Phi(z,t)=\big(z,(\theta(t_a))_{a\in A}\big).
\]
Since $\theta(t)=0$ only for $t=0$, the coordinates $t_a$ that vanish at $\Phi(p)$ are those that vanish at $p$. Restrict the corner chart to $\{t_a>0\}$ for the coordinates with $t_a(p)>0$, which then become coordinates of the $l$-type, and shrink it so that $\theta(t)=t^2$ for the vanishing coordinates. In the chart of $\hat Q$ at $\Phi(p)$ given by Lemma~\ref{lem_Qhat}, which uses $\sqrt{t_a}$ for these coordinates and $t_a$ for the others, $\Phi$ reads $t_a\mapsto\sqrt{\theta(t_a)}=t_a$ near $0$ for the vanishing coordinates and $t_a\mapsto\theta(t_a)$ for the others, and $z\mapsto z$. Both are smooth with positive derivative, and so are the inverses ($t\mapsto t$ near $0$, and $\theta^{-1}$ away from $0$). So $\Phi$ is a local diffeomorphism $Q\to\hat Q$ in both directions, and a global homeomorphism, hence a diffeomorphism of manifolds with corners. It preserves orientation, since on $\inter Q=\inter\hat Q$ it is a diffeomorphism with positive derivatives in the collar coordinates. Orientation is used only here, so the theorem also holds for non-orientable $Q$, as stated in Theorem~\ref{thmB}.
\end{proof}

Something like $\Phi$ is indispensable: the set-theoretic identity $\hat Q\to Q$ is $\pi$, which is smooth but has zero normal derivative at the faces. Only the existence of an isomorphism $Q\cong\hat Q$ is used; triangulating $\hat Q$ by $\hat t$ and using $[Q,c]=[\hat Q,c]$ is equivalent.

\begin{corollary}[the point]\label{cor_point}
For every source $Q$ of dimension $k\ge1$, $[Q,c]\in\operatorname{span}\cG(\mathrm{pt})$. Hence $\GCG_k(\mathrm{pt})=0$ for $k\ge1$ at chain level, and \textup{(A6)} holds for $G$; and $GC^{G^W}_k(\mathrm{pt})=\GCE_k(\mathrm{pt})=0$ for $k\ge1$, by Corollary~\ref{thm_T7} alone.
\end{corollary}

\begin{proof}
By Theorem~\ref{thm_whitehead} and Lemma~\ref{lem_glue}(e), $[Q,c]-\sum_i\varepsilon_i[\Delta^k,c]\in\cG(\mathrm{pt})$. For $k\ge1$, $[\Delta^k,c]=0$ in $\Cm$, because a vertex transposition is an orientation-reversing automorphism commuting with $c$. For $G^W$ and $E$ use Corollary~\ref{thm_T7} and Theorem~\ref{thm_readings} instead.
\end{proof}

\begin{remark}[triangulation-free proofs]\label{rem_lowdim}
We do not know a proof of (A6) for $G$ that avoids triangulations. Products are harmless, since $[Q\times I,c]=0$ already in $\Cm$ by the flip of $I$, but this does not kill $[Q,c]$. Handle decompositions are not aligned: a handle is attached along part of a face of the previous stage. Cutting along regular levels of a Morse function gives aligned two-piece relations, but the pieces are elementary cobordisms, which need not have orientation-reversing automorphisms.
\end{remark}

\section{\texorpdfstring{The simplex-only relation $R_\Delta$}{The simplex-only relation R\_Delta}}\label{sec_R}

Let $R_\Delta(Y)$ be the span of the relations $r(Q,f,t)$ for $t\in\cT(Q)$; it is a subcomplex by Lemma~\ref{lem_bdry}, and $\GCRD(Y)=\Cm(Y)/R_\Delta(Y)$. By Theorem~\ref{thm_readings}, $R_\Delta\subseteq G^W\subseteq E$. The inclusion into $G$ itself holds for the Whitehead version $R_\Delta^{C^1}$, spanned by the $r(Q,f,t)$ with $t\in\cT^{C^1}(Q)$ (Lemma~\ref{lem_glue}(e)); whether $R_\Delta\subseteq G$ at chain level we do not know, because the simplices of an admissible triangulation need not be immersions.

\subsection{The structure theorem}
Let $\Calt_*(Y)$ be $\Csm_*(Y)$ modulo $\sigma\circ a=\sgn(a)\sigma$ for affine permutations $a$ of $\Delta^n$. By acyclic models \cite{EM}, $\Csm\to\Calt$ is a quasi-isomorphism, and $\Calt(W)\to\Calt(Y)$ is injective for $W\subseteq Y$, since both are spanned by orbits of simplices. Since $[\Delta^n,\sigma\circ a]=\sgn(a)[\Delta^n,\sigma]$, $p$ factors through $\bar p:\Calt\to\Cm$. Let
\[
K(Y)=\operatorname{span}\{\tau(Q,f,t)-\tau(Q,f,t'):t,t'\in\cT(Q)\}.
\]

\begin{proposition}\label{prop_structure}
$\bar p:\Calt(Y)\to\GCRD(Y)$ is surjective at chain level, with kernel $K(Y)$. Hence there is a long exact sequence
\[
\cdots\to H_n(K(Y))\to H_n(\Csm(Y))\xrightarrow{p_*}H_n(\GCRD(Y))\to H_{n-1}(K(Y))\to\cdots.
\]
\end{proposition}

\begin{proof}
Surjectivity: $[Q,f]\equiv\bar p\,\tau(Q,f,t)$ modulo $R_\Delta$, by Corollary~\ref{thm_T7}. $K$ is a subcomplex by Lemma~\ref{lem_bdry}. Define $\Theta:\Cm(Y)\to\Calt(Y)/K$ by $[Q,f]\mapsto\tau(Q,f,t)$. It is well defined: if $\sum_ac_a[Q_a,f_a]=0$ in $\Cm$, group the connected terms by isomorphism class $y$ up to orientation, with representatives $(Q_y,f_y)$ and signs $\delta_a$; if $y$ has no orientation-reversing automorphism then $\sum_{a\in y}c_a\delta_a=0$ and $\sum_{a\in y}c_a\tau_a=\sum c_a\delta_a(\tau(Q_y,f_y,h_at_a)-\tau(Q_y,f_y,t_y))\in K$; if it has one, $\rho$, then $\tau(Q_y,f_y,u)=\tfrac12(\tau(u)-\tau(\rho u))\in K$. $\Theta$ is a chain map by Lemma~\ref{lem_bdry}, $\Theta(r(Q,f,t))=0$, and $\Theta\circ\bar p=\id$, since $\tau(\Delta^n,\sigma,\id)=\sigma$. So $\ker\bar p\subseteq K$. Conversely $\bar p(K)\subseteq R_\Delta$ and $\bar p$ is a chain map.
\end{proof}

The relation family of triangulation data satisfies (A1) and (A5) (Lemma~\ref{lem_dec}(b) and Proposition~\ref{prop_additive}), so $p_*$ is injective for $\GCRD$ in all degrees (Propositions~\ref{prop_stokes}, \ref{prop_poly_inj}), and Proposition~\ref{prop_structure} gives
\[
H_n(\GCRD(Y))\cong H_n(Y;\Q)\oplus H_{n-1}(K(Y)).
\]
The prism of a triangulation datum is not a triangulation datum, but prisms preserve $R_\Delta$: with the staircase products of Subsection~\ref{ssec_adm},
\[
\cH r(Q,f,t)=(-1)^n\Big(r(Q\times I,H\circ(f\times\id),t\times I)-\sum_i\varepsilon_i\,r(\Delta^n\times I,H\circ(f\sigma_i\times\id),\id\times I)\Big),
\]
where $\id\times I$ is the staircase triangulation of the prism; so Lemma~\ref{lem_star}(a) holds for $R_\Delta$. At the point, $\GCRD_k(\mathrm{pt})=0$ for $k\ge1$ at chain level, by Corollary~\ref{thm_T7} as in Corollary~\ref{cor_point}. Hence the proof of Lemma~\ref{lem_star}(b) applies to $R_\Delta$ for every $m$: it uses only (a), the vanishing at the point and Proposition~\ref{prop_stokes} in degree $0$, which needs only (A1) and (A5). Axioms (A2) and (A3) fail as stated: a bisection is not a triangulation relation, and the cut of a triangulation relation is a curvilinear tiling, which lies in $G^W$ by the proof of (A3) in Theorem~\ref{thm_readings} but in $R_\Delta$ only if the triangulation can be subdivided coherently along the level.

\begin{proposition}\label{prop_RG}
For $Y$ a compact polyhedron or a smooth manifold without boundary, the map $H_n(\GCRD(Y))\to H_n(GC^{G^W}(Y))$ induced by $R_\Delta\subseteq G^W$ is onto for every $n$, with kernel $H_{n-1}(K(Y))$. The same holds for $R^{C^1}_\Delta\subseteq G$, with $K$ formed from Whitehead triangulations. So the theorem for $G$ implies the theorem for $R_\Delta$ exactly when $K$ is acyclic.
\end{proposition}

\begin{proof}
The composite with $p_*$ for $\GCRD$ is $p_*$ for $G^W$, an isomorphism by Theorem~\ref{thm_readings}; this gives surjectivity, and with the splitting above the kernel. For $R^{C^1}_\Delta$ the same argument applies, with Theorem~\ref{thm_whitehead} in place of Corollary~\ref{thm_T7} in Proposition~\ref{prop_structure}, and Theorem~\ref{thm_A_G}.
\end{proof}

\subsection{Coherent subdivisions and low degrees}
A \emph{smooth subdivision} of a finite simplicial complex $L$ is a homeomorphism $s:|L'|\to|L|$ from a simplicial complex such that each closed simplex of $L'$ is mapped smoothly into a closed simplex of $L$ and each open simplex into an open simplex of $L$.

\begin{definition}\label{def_RT}
Let $j\ge1$. Property $\mathrm{RT}_j$ holds if the following is true. Let $Q\in\cS$ with $\dim Q=j$, $t\in\cT(Q)$ with complex $L$ and simplices $\sigma_\rho:\rho\to Q$, $h:Q\to\R$ smooth, and $c$ a regular value of $h$ on every depth stratum of $Q$ and of $h\circ\sigma_\rho$ on every open face of every simplex $\rho$ of $L$. Then there is a smooth subdivision $s:|L'|\to|L|$ such that every open simplex of $L'$ lies in one of the subsets $\{h\circ t\circ s<c\}$, $\{h\circ t\circ s=c\}$, $\{h\circ t\circ s>c\}$ of $|L'|$. We set $\mathrm{RT}_0$ to be vacuous.
\end{definition}

Given such $s$, put $T=t\circ s$ and let $T_\pm$, $T_0$ be its restrictions to the subcomplexes $L'_\pm$, $L'_0$ over $Q_\pm$, $L_Q$ (for $h$ and $c$). Over each top simplex $\rho_i$, $s$ restricts to triangulations $S_{i,\pm}$, $S_{i,0}$ of the pieces $\Delta_{i,\pm}$ and $L_{\Delta_i}$ of $\Delta^j$ for $h\sigma_i$.

\begin{lemma}\label{lem_RT_cons}
In this situation, $T\in\cT(Q)$, $T_\pm\in\cT(Q_\pm)$, $T_0\in\cT(L_Q)$, $S_{i,\pm}\in\cT(\Delta_{i,\pm})$, $S_{i,0}\in\cT(L_{\Delta_i})$, and for every smooth $f$ on $Q$:
\begin{enumerate}[(a)]
\item $\tau(Q_\pm,f,T_\pm)=\sum_i\varepsilon_i\,\tau(\Delta_{i,\pm},f\sigma_i,S_{i,\pm})$ and $\tau(L_Q,f,T_0)=\sum_i\varepsilon_i\,\tau(L_{\Delta_i},f\sigma_i,S_{i,0})$ in $\Calt$;
\item if $h=\ell\circ f$, then $\cut_\pm r(Q,f,t)=r(Q_\pm,f,T_\pm)-\sum_i\varepsilon_ir(\Delta_{i,\pm},f\sigma_i,S_{i,\pm})$ and $\Lambda\,r(Q,f,t)=r(L_Q,f,T_0)-\sum_i\varepsilon_ir(L_{\Delta_i},f\sigma_i,S_{i,0})$;
\item if $Q=\Delta^j$, $t=\id$ and $h=\ell\circ g$, then $[\Delta^j,g]-\cut_-[\Delta^j,g]-\cut_+[\Delta^j,g]\in R_\Delta$.
\end{enumerate}
\end{lemma}

\begin{proof}
Admissibility. The closed simplices of $T_-$ are $\sigma_\rho\circ s|$, smooth into $Q$ with image in $Q_-$, hence smooth into $Q_-$ by the slice chart of Lemma~\ref{lem_pieces}(a); similarly for $L_Q$ and for the pieces of $\Delta^j$. An open simplex of $L'$ lies in an open face $\mathring\varphi$ of some $\rho$, and $\sigma_\rho(\mathring\varphi)$ lies in one stratum of $Q$, intersected with one of $\{h<c\}$, $\{h=c\}$, $\{h>c\}$; by Lemma~\ref{lem_pieces}(a) these intersections are unions of components of strata of $Q_\pm$ and of $L_Q$, and the image of the connected open simplex lies in one of these components. No immersivity of $\sigma_\rho$ and no $\sigma_\rho^{-1}$ is used.

(a) The top simplices of $T_-$ are the $\sigma_i\circ e$ with $e$ a top simplex of $S_{i,-}$, and local degrees multiply for homeomorphisms onto the image, which gives the first identity simplex by simplex. For the second we use a chain-level argument that needs no comparison of normals, so that it applies to non-immersive $\sigma_i$. Apply $\partial$ to the identity $\mu_{T_-}=\sum_i\varepsilon_i\sigma_{i*}\mu_{S_{i,-}}$ in $\Calttop(Q_-)$ and use Lemma~\ref{lem_bdry} on both sides. Then project onto the span of the orbits of simplices with image in $L_Q$; this is a well-defined linear map on $\Calttop(Q_-)$. On the left, a boundary $(j-1)$-simplex of $T_-$ over $(\partial Q)_-$ has its open interior in $\Bord_1(Q)\cap\{h<c\}$, since $\Bord_1(Q)\cap\{h=c\}$ has dimension $j-2$; so it is not in $L_Q$, and the projection is $\mu_{T_0}$. On the right, a boundary simplex of $S_{i,-}$ over $(\partial\Delta)_-$ has its interior in some $\mathring\varphi\cap\{h\sigma_i<c\}$, so its image lies in $\{h<c\}$ and is not contained in $L_Q$. So $\mu_{T_0}=\sum_i\varepsilon_i\sigma_{i*}\mu_{S_{i,0}}$, first in $\Calttop(Q_-)$ and then in $\Calttop(L_Q)$ by injectivity. Push forward by $f$.

(b) By linearity, $\cut_-r(Q,f,t)=[Q_-,f]-\sum_i\varepsilon_i[\Delta_{i,-},f\sigma_i]$; insert (a). The same for $\Lambda$, with the orientations of $L_{\Delta_i}$ and $L_Q$ as faces of $\Delta_{i,-}$ and $Q_-$, which is what Lemma~\ref{lem_bdry} uses.

(c) $[\Delta^j,g]-[\Delta_-,g]-[\Delta_+,g]=r(\Delta^j,g,T)-r(\Delta_-,g,T_-)-r(\Delta_+,g,T_+)$, since $\tau(T)=\tau(T_-)+\tau(T_+)$.
\end{proof}

\subsection{Dimensions one and two}
$\mathrm{RT}_1$ holds: in dimension one the simplices are smooth homeomorphisms of $[0,1]$ onto arcs, the level points on each open edge are finitely many and not vertices, and one subdivides affinely at them.

\begin{lemma}[strip]\label{lem_strip}
Let $e_1,\dots,e_k$ and $e'_1,\dots,e'_m$ be two triangulations of $[0,1]$ whose edges are smooth increasing homeomorphisms $e_i:[0,1]\to[a_{i-1},a_i]$, $e'_i:[0,1]\to[a'_{i-1},a'_i]$, possibly with vanishing derivatives. Then there is $S\in\cT([0,1]^2)$ that restricts to $(e_i)$ on $[0,1]\times0$ and to $(e'_i)$ on $[0,1]\times1$, and has the three affine segments $\{0,1\}\times[0,\tfrac13],[\tfrac13,\tfrac23],[\tfrac23,1]$ on the sides. The same holds for circles, in a coordinate $\R/\Z$, with at least three vertices.
\end{lemma}

\begin{proof}
Let $\ell_i$ be the affine parametrisation of $[a_{i-1},a_i]$ and
\[
\Phi_i(s,v)=\big((1-v)e_i(s)+v\ell_i(s),\ v/3\big),\qquad(s,v)\in[0,1]^2.
\]
For fixed $v$ the first coordinate is strictly increasing in $s$, so $\Phi_i$ is a smooth homeomorphism onto $[a_{i-1},a_i]\times[0,\tfrac13]$, affine on the vertical sides; cutting $[0,1]^2$ along a diagonal and composing with $\Phi_i$ gives two smooth closed triangles. Do the same at the top with $(e'_i)$, and fill $[0,1]\times[\tfrac13,\tfrac23]$ by the linear cone from its centre over the boundary polygon. The corners are vertices, open boundary edges lie in open sides, and everything else is interior. In the circle case the strips over $[0,\tfrac13]$ and $[\tfrac23,1]$ are built in the same way in $\R/\Z$, but the middle annulus $\R/\Z\times[\tfrac13,\tfrac23]$ has no centre, and cutting it open and coning would give two triangles with the same vertices. Triangulate it instead by an affine zipper: let $b_1,\dots,b_k$ and $a_1,\dots,a_m$, $k,m\ge3$, be the vertices on $\R/\Z\times\tfrac13$ and $\R/\Z\times\tfrac23$ in cyclic order, merge the two cyclic orders, and for each pair of consecutive vertices on one circle take the triangle with the vertex of the other circle that is current in the merged order. Each triangle has one edge on one circle and its third vertex on the other, each connecting edge lies in exactly two triangles, and since $k,m\ge3$ no two triangles have the same vertex set; the triangles are affine in the universal cover $\R\times\R$ with width less than $1$, hence embedded. So this is a simplicial complex, and the conclusion holds as before.
\end{proof}

\begin{lemma}[collar extension in dimension two]\label{lem_collar}
Let $P$ be a compact $2$-dimensional source and $\beta:|J_0|\to B_P$ a triangulation by a simplicial complex in which every corner of $P$ is a vertex and every closed edge is a smooth map into $P$ that is a homeomorphism onto an arc of a single face. Then there is $T\in\cT(P)$ whose restriction to $B_P$ is $\beta$.
\end{lemma}

\begin{proof}
Choose a smooth vector field $V$ on $P$ pointing strictly inward at every point of $B_P$; it exists by a partition of unity, since the inward cones are convex. Extend $V$ to an open manifold containing $P$ (locally by Seeley extension and a partition of unity); its flow $\varphi_s$ is defined on $P$ for $|s|\le s_0$, is injective with $\varphi_s\varphi_{s'}=\varphi_{s+s'}$, and for $s\ge0$ maps $P$ into $P$, since orbits can leave $P$ only through $B_P$, where $V$ points inward, as in the proof of \cite[Corollary on orbifolds with corners]{Yu1}. Put $\Psi(y,s)=\varphi_s(y)$ on $B_P\times[0,s_0]$.

$\Psi$ is injective: if $\varphi_s(y)=\varphi_{s'}(y')$ with $s<s'$, then $y=\varphi_{s'-s}(y')\in\inter P$. The set $P':=\varphi_{s_0}(P)$ is $P\setminus\Psi(B_P\times[0,s_0))$, since a point outside the collar flows backwards for time $s_0$ without leaving $P$. Take $T_0\in\cT(P)$ (Corollary~\ref{thm_T7}); then $T'=\varphi_{s_0}\circ T_0$ triangulates $P'\subseteq\inter P$, and all its open simplices have depth $0$ in $P$.

For each face arc $\gamma$ between consecutive corners, with a regular parametrisation $\gamma:[0,1]\to P$, the map $\Phi(u,s)=\varphi_{ss_0}(\gamma(u))$ is a smooth embedding of $[0,1]^2$: $\gamma'$ is tangent to the face, $V$ is transverse to it, and the flow preserves linear independence. Apply Lemma~\ref{lem_strip} with bottom $\gamma^{-1}\circ\beta|_\gamma$ and top $\gamma^{-1}\circ T_0|_\gamma$; both are triangulations of $[0,1]$ by smooth increasing homeomorphisms, since corners are vertices of both. Compose with $\Phi$; the top edges are those of $T'$ on $\partial P'$. Circle faces are treated with the circle case. Adjacent quadrilaterals share the flow line $\Psi(q,\cdot)$ from a corner $q$, subdivided at the same parameters. The only simplices of a strip with all vertices on its top are top edges and vertices, since the diagonals cross the layers, so the union with $T'$ is a simplicial complex along $\partial P'$; two strips sharing two corner flow lines (a face with two corners and two arcs) could only produce a repeated vertex set if both arcs had a single edge, which is excluded since each boundary circle of $P$ carries at least three vertices of $\beta$ and of $T_0$. Corners are vertices, bottom open edges lie in faces, and all other open simplices are interior.
\end{proof}

\begin{proposition}[$\mathrm{RT}_2$]\label{prop_RT2}
$\mathrm{RT}_2$ holds. More generally, finitely many functions $h_m$ with pairwise disjoint levels $\{h_m=c_m\}$ can be treated at once.
\end{proposition}

\begin{proof}
Vertices are off the level. Subdivide each edge of $L$ affinely at its finitely many level points. In a $2$-simplex $\rho$ the set $\{h\sigma_\rho=c\}$ is a compact $1$-submanifold, consisting of arcs that end transversally on open edges at the subdivision points, and of circles. Place at least two interior vertices on each arc and at least three on each circle, and parametrise the sub-arcs regularly; this avoids two edges with the same vertex pair, e.g.\ for an arc with both ends on one edge, or for a lens bounded by one arc and one edge segment. Each component $P$ of each piece of $\rho$ is a compact planar domain with corners; its faces are edge segments, level arcs and level circles, all embedded, so $P\in\cS$. Its boundary is triangulated as required in Lemma~\ref{lem_collar}, which gives $S_P$. The $S_P$ agree along shared level arcs, and each edge of $L$ was subdivided once, so the union is a smooth subdivision of $L$. For several disjoint levels, subdivide the edges at all level points; the pieces are bounded by disjoint level arcs and the argument is the same.
\end{proof}

\subsection{Small simplices}
\begin{lemma}\label{lem_small}
Let $Y\subseteq\R^K$ have a cover $\cU$ by sets open in $Y$ with $H_k(\GCRD(U))=0$ for all $U\in\cU$. For instance, for a compact polyhedron $Y$ with a linear triangulation, the open vertex stars $\St(w)$ (Proposition~\ref{prop_amb}) are open in $Y$, though not in $\R^K$, and $\St(w)$ is star-shaped about $w$, so Lemma~\ref{lem_star}(b) applies to it. Let $k\ge1$ and assume $\mathrm{RT}_{k-1}$. Let $x\in\Cm_k(Y)$.
\begin{enumerate}[(a)]
\item $x\equiv x':=p(\Sd^N\tau(x))$ modulo $R_\Delta(Y)$, where $\tau$ uses any triangulations of the sources, and $N$ can be chosen so that every simplex of $x'$ has image in a member of $\cU$.
\item If $c$ is such that all simplices of $x'$ are $c$-transverse, then $x'\equiv\cut_-x'+\cut_+x'$ modulo $R_\Delta(Y)+\partial\Cm_{k+1}(Y)$.
\end{enumerate}
\end{lemma}

\begin{proof}
(a) $[Q,f]-p\tau(Q,f,t)=r(Q,f,t)$, and $p(\sigma)-p(\Sd\sigma)=r(\Delta^k,\sigma,b)$ with $b$ the barycentric triangulation of $\Delta^k$. Small images follow from uniform continuity and a Lebesgue number.

(b) $x'-\cut_-x'-\cut_+x'=\sum_\sigma c_\sigma y_\sigma$ with $y_\sigma=[\Delta^k,\sigma]-[\Delta_-,\sigma]-[\Delta_+,\sigma]\in\Cm_k(U_\sigma)$, $\sigma(\Delta^k)\subseteq U_\sigma\in\cU$. By \eqref{eq_cut}, $\partial y_\sigma=\sum_j(-1)^jy_{\sigma d_j}-(\Lambda-\Lambda)[\Delta^k,\sigma]=\sum_j(-1)^jy_{\sigma d_j}$, and each $y_{\sigma d_j}$ lies in $R_\Delta(U_\sigma)$ by Lemma~\ref{lem_RT_cons}(c) for $j=k-1$ (for $k=1$ the faces are points off the level and $y_{\sigma d_j}=0$). So $y_\sigma$ is a cycle of $\GCRD(U_\sigma)$ and lies in $\partial\Cm_{k+1}(U_\sigma)+R_\Delta(U_\sigma)$.
\end{proof}

No homological input about $Y$ is used, and each simplex is cut once.

\subsection{\texorpdfstring{The theorem for $R_\Delta$}{The theorem for R\_Delta}}
\begin{lemma}[$\mathrm{MV}_k$ for $R_\Delta$]\label{lem_MV_R}
Lemma~\ref{lem_MV} holds with $\GCRD$ in place of $\GCcR$, provided $\mathrm{RT}_{k-1}$ holds.
\end{lemma}

\begin{proof}
Let $x\in\Cm_k(Y)$ with $\partial x\in R_\Delta(Y)$. Fix a linear triangulation of $Y$ and replace $x$ by $x'$ as in Lemma~\ref{lem_small}(a), for the cover of $Y$ by its open vertex stars. Since $x-x'\in R_\Delta(Y)$ at chain level and $R_\Delta$ is a subcomplex (Lemma~\ref{lem_bdry}), $\partial x'=\partial x-\partial(x-x')\in R_\Delta(Y)$; as $R_\Delta(Y)$ is the span of the $r(Q,f,t)$, we can write $\partial x'=\sum_l\lambda_lr_l$ with $r_l=r(Q_l,f_l,t_l)$ and $\dim Q_l=k-1$.

\emph{The case $k=1$.} Then $(R_\Delta)_0(Y)=0$: a $0$-dimensional admissible triangulation is a signed bijection of points, so $r(Q,f,t)=0$ for $\dim Q=0$. Hence $\partial x'=0$, the presentation is empty, and the boundary terms below vanish; Lemma~\ref{lem_RT_cons}, which is stated only for $j\ge1$, is not needed.

\emph{Choice of the level.} Choose $c$ in the given set such that the simplices of $x'$, the $(Q_l,f_l)$ and the $[\Delta^{k-1},f_l\sigma_{l,i}]$ are $c$-transverse (Lemma~\ref{lem_generic}). Faces of transverse generators are transverse, so every chain written so far lies in $C^\pitchfork_c(Y)$, and $\cut_\pm$, $\Lambda$ do not depend on the presentation (Lemma~\ref{lem_cutops}). For $k\ge2$ this is exactly the hypothesis of $\mathrm{RT}_{k-1}$ for $(Q_l,t_l,\ell f_l)$: by Lemma~\ref{lem_adm}(a) every simplex $\rho$ of the complex of $t_l$ is a face of a top simplex $\rho_i$, and $\sigma_\rho$ is the restriction of $\sigma_{l,i}$ to it, so $c$ is a regular value of $\ell f_l\sigma_\rho$ on every open face of every $\rho$; and $c$ is a regular value of $\ell f_l$ on every depth stratum of $Q_l$ because $(Q_l,f_l)$ is $c$-transverse.

\emph{The boundary terms} (Step~2 of Lemma~\ref{lem_MV}). By Lemma~\ref{lem_RT_cons}(b) with $j=k-1$, $\theta_\pm:=\cut_\pm\partial x'\in R_\Delta(Y_\pm)$ and $\Lambda\partial x'\in R_\Delta(Y_0)$.

\emph{The splitting} (Step~3 of Lemma~\ref{lem_MV}). In Lemma~\ref{lem_MV} the identity $x\equiv\cut_-x+\cut_+x$ is a relation, by Lemma~\ref{lem_cutrel}. For $R_\Delta$ it is replaced by Lemma~\ref{lem_small}(b):
\[
x'\equiv x'_-+x'_+\pmod{R_\Delta(Y)+\partial\Cm_{k+1}(Y)},\qquad x'_\pm:=\cut_\pm x',
\]
which holds only modulo $R_\Delta$ plus boundaries, not as a relation. The other identities of Step~3, $\partial x'_\pm=\theta_\pm\pm z$ with $z=\Lambda x'$ and $\partial z=-\Lambda\partial x'\in R_\Delta(Y_0)$, are \eqref{eq_cut}.

\emph{Steps 4--7} of Lemma~\ref{lem_MV} then apply with $R_\Delta$ in place of $\cR$ and $x'$ in place of $x$; injectivity on $Y_\pm$ in degree $k-1$ is Proposition~\ref{prop_poly_inj} for the family of triangulation data, which satisfies (A1) and (A5). Step~7 gives $x'\equiv A_-+A_++p(s)\equiv p(a_-+a_++s)$ modulo $R_\Delta(Y)+\partial\Cm_{k+1}(Y)$. This suffices, because Step~7 concludes modulo $R_\Delta(Y)+\partial\Cm_{k+1}(Y)$ anyway, so the weaker form of the splitting costs nothing. Finally $x\equiv x'$ modulo $R_\Delta(Y)$.
\end{proof}

\begin{theorem}\label{thm_R}
Let $Y$ be a compact polyhedron or a smooth manifold without boundary, and $m\ge1$. If $\mathrm{RT}_j$ holds for $j\le m-1$, then $p_*:H_k(Y;\Q)\to H_k(\GCRD(Y))$ is an isomorphism for $k\le m$; it is injective for all $k$. Moreover $\mathrm{RT}_1$ and $\mathrm{RT}_2$ hold, so $p_*$ is an isomorphism for $k\le3$.
\end{theorem}

\begin{proof}
Injectivity is Propositions~\ref{prop_stokes} and \ref{prop_poly_inj}. For surjectivity in degrees $\le m$ the proofs of Propositions~\ref{prop_abs_poly} and \ref{prop_abs_mfd} apply verbatim, with Lemma~\ref{lem_MV_R} in place of Lemma~\ref{lem_MV}. In the induction $S(k)$, for $1\le k\le m$: the star-shaped half is handled by Lemma~\ref{lem_star}(b), which applies to $R_\Delta$ as shown above; the other half by the radial projection, Lemma~\ref{lem_star}(a) for $R_\Delta$, and the inner induction on vertices; the slice by $S(k-1)$; and $\mathrm{MV}_k$ for $R_\Delta$ needs $\mathrm{RT}_{k-1}$, which holds since $k-1\le m-1$. General polyhedra and manifolds follow as in Propositions~\ref{prop_abs_poly} and \ref{prop_abs_mfd}, by naturality of $R_\Delta$. No homological statement about a space is used before it is proved, and $H_1(K)=0$ is not used. The last claim uses $\mathrm{RT}_1$ and Proposition~\ref{prop_RT2}.
\end{proof}

By Proposition~\ref{prop_structure}, Theorem~\ref{thm_R} in degrees $\le3$ says that $H_1(K(Y))=H_2(K(Y))=0$. The cases $k\le2$ also follow from the structure theorem and $H_1(K)=0$, which is proved generator by generator with the prism over a strip (Lemma~\ref{lem_strip}): for a $1$-dimensional $Q$ and $t,t'\in\cT(Q)$, the chain $\tau(Q\times I,f\circ\pr,S)$ lies in $K$ by the flip of the interval, and has boundary $\pm(\tau(t)-\tau(t'))$ up to constant side terms, which vanish in $\Calt$.

\subsection{Relative triangulations with corners}
$\mathrm{RT}_j$ for $j\ge3$ would follow from a relative triangulation theorem in the class of $C^1$ (Whitehead) triangulations. The collar construction of $\mathrm{RT}_2$ does not lift: it would need interpolation between two smooth triangulations of $(\partial P)\times I$ with $\dim\partial P\ge2$, which is not available in the literature as far as we know; it would follow from Conjecture~\ref{conj_reltri}. Note that Theorem~\ref{thm_whitehead} is absolute, not relative, and does not help here.

A \emph{PD boundary triangulation} of a compact manifold with corners $P$ is a homeomorphism $\beta:|J_0|\to B_P$ from a finite simplicial complex such that each closed simplex is a smooth immersion into $P$ contained in a single face, and each open simplex lies in one depth stratum.

\begin{conjecture}[relative $C^1$ triangulation with corners]\label{conj_reltri}
Let $P$ be a compact $j$-manifold with corners and $\beta$ a PD boundary triangulation. Then there are a finite simplicial complex $J\supseteq J_0$ and a homeomorphism $T:|J|\to P$ with $T|_{|J_0|}=\beta$, with the same simplicial structure on $J_0$, such that each closed simplex of $J$ is a smooth immersion into $P$ and open simplices not in $J_0$ lie in $\inter P$. In particular $T\in\cT(P)$; here $\cT(P)$ is defined by Definition~\ref{def_adm} for any compact manifold with corners, with embedded faces or not.
\end{conjecture}

We state this as a conjecture because the argument we have is only a sketch and has not been checked. It runs as follows: thicken $\beta$ by an inward collar, which is a piecewise differentiable embedding of $J_0\times[0,s_0]$ with its staircase triangulation; make it and finitely many interior cube charts compatible by successive $C^1$-small modifications that never touch the collar, as in the existence proof of Whitehead triangulations in \cite[Lectures 3 and 4]{Lur}, which is the corner version of \cite[Lecture 4, Variant 7]{Lur}; and finally retriangulate the collar prisms by cones so that $J_0\times0$ is not subdivided. For manifolds with faces, Massuyeau and Turaev state that the standard methods of \cite[Section~10.6]{Mun} give smooth triangulations in which all faces are subcomplexes \cite[Section~3.1.1]{MT}.

We know no reference for the corner-relative statement with a prescribed boundary triangulation. Equivariant and stratified triangulation theorems (Illman \cite{Ill}, Goresky \cite{Gor}, Verona \cite{Ver}) are relative to subsets, not to prescribed triangulations of them, and so do not give coherence across simplices.

\begin{remark}[expected consequence]\label{cor_RT_all}
Assuming Conjecture~\ref{conj_reltri}, we expect that $\mathrm{RT}_j$ holds for all $j$, and $p_*:H_*(Y;\Q)\to H_*(\GCRD(Y))$ is an isomorphism in all degrees for compact polyhedra and smooth manifolds without boundary. In particular $R_\Delta$, $G$ and $E$ would give the same homology.
\end{remark}

The reduction inducts over the skeleta of $L$ inside the affine simplices, applying Conjecture~\ref{conj_reltri} first to the level set in each simplex, relative to its already triangulated boundary, and then to each piece, relative to its boundary. It has not been verified: only its logical form has been checked, together with the fact that the hypotheses of Conjecture~\ref{conj_reltri} on the boundary triangulations are met inductively; the immersivity bookkeeping at the corners of the pieces has not been written out.

\subsection{Remarks on other classes}
\begin{remark}[a second route to injectivity]\label{rem_szpilrajn}
For $R_\Delta$ we sketch an argument for injectivity that does not use differential forms and would work for Lipschitz maps into Riemannian manifolds and polyhedra; its support and nerve steps are not written out, and nothing else in the paper depends on it. By Proposition~\ref{prop_structure}, one must show that every cycle $k=\sum_jc_j(\tau(Q_j,f_j,t_j)-\tau(Q_j,f_j,t'_j))$ of $K$ bounds in $Y$. Since $\mu_{t_j}$ and $\mu_{t'_j}$ represent the same fundamental class in $H_n(Q_j,B_{Q_j})$, $\mu_{t_j}-\mu_{t'_j}=\partial\beta_j+\gamma_j$ with $\gamma_j$ supported in $B_{Q_j}$, even for topological triangulations; so $[k]$ is represented by a cycle supported in $A=\bigcup_jf_j(B_{Q_j})$. This is a Lipschitz image of an $(n-1)$-dimensional set, so its $n$-dimensional Hausdorff measure vanishes and $\dim A\le n-1$ by Szpilrajn's theorem \cite{Szp,HW}; a nerve of a fine cover of $A$ of order $\le n-1$, mapped back into $Y$ through a neighbourhood retraction, then kills $[k]$. For $G$, $G^W$, $G^T$ and $E$ the correction term is supported in $\bigcup_i\iota_i(B_{Q_i})$, whose pairwise intersections can be arbitrary compacta; by Lemma~\ref{lem_lip_goodpair} the relative class is nevertheless detected by local degrees, so the argument also applies to these relations; see Remark~\ref{rem_lip_nerve}.
\end{remark}

\begin{remark}[Lipschitz maps]\label{rem_lipschitz}
The surjectivity arguments of this paper need regular values, and fail for Lipschitz maps. Injectivity holds in all degrees for Lipschitz maps (Theorem~\ref{thm_lip_inj}). For instance, let $W:[0,1]\to[0,1]$ be convex, $1$-Lipschitz and non-differentiable on a dense set, such as a normalised $\sum_n2^{-n}|x-q_n|$, and $f(x,y)=\tfrac12(y+W(x))$ on $[0,1]^2$, composed with an affine coordinate of an edge of a graph; every level in an open interval is the graph of a function that is nowhere $C^1$, so no level cuts $[0,1]^2$ into manifolds with corners. In the Lipschitz class $p_*$ is an isomorphism in degrees $\le1$ (Proposition~\ref{prop_lip_low}); we do not know whether it is onto in degree $2$, where we expect it, or in degrees $\ge3$. Already the case of a wedge of circles is open there in degrees $\ge2$, although its smooth analogue holds by Theorem~\ref{thm_A_G} and Corollary~\ref{cor_dim}.
\end{remark}

\section{Orbifold sources}\label{sec_orb}
In this section sources are allowed to be orbifolds. We show that the gluing relation, read on underlying spaces as in \cite{Yu0}, then also replaces $D$: pieces may be manifolds even when the glued source is an orbifold, and this is what kills the cycles of \cite[Remark on the role of degenerate chains]{Yu1}.

\subsection{Setting}
Let $\cS^{\mathrm{orb}}$ be the class of compact oriented effective orbifolds with corners in the sense of \cite[Definitions of orbifolds with corners]{Yu1}, possibly disconnected or empty. Faces need not be embedded, and local groups may permute boundary hyperplanes; the boundary $i_Q:\partial Q\to Q$ is that of \cite[Convention on the boundary]{Yu1}. A map $f:|Q|\to\R^K$ is \emph{smooth} if it lifts on charts to invariant smooth maps; since the target is a manifold, this agrees with liftability, and an orbifold morphism to a manifold is the same as an invariant smooth map on charts \cite[Lemma on orbifold morphisms, (a)]{Yu1}. The complex $C^{\mathrm{orb}}_*(Y)$ is defined as $\Cm_*(Y)$ in Section~\ref{sec_setting}, with sources in $\cS^{\mathrm{orb}}$; the basis lemma \cite[Lemma on the basis of $C^g_*$]{Yu1} is stated there for orbifolds. A manifold with corners, with embedded faces or not, is an element of $\cS^{\mathrm{orb}}$, and $\Cm_*(Y)\subseteq C^{\mathrm{orb}}_*(Y)$ is spanned by a subset of the basis. All constructions below (boundaries, cut pieces and level sets, products, the pieces of Proposition~\ref{prop_orb_point}) stay in $\cS^{\mathrm{orb}}$, so the weights $1/|K|$ of non-effective sources in \cite{Yu1} never occur.

\begin{definition}[orbifold gluing decompositions]\label{def_glue_orb}
Let $P\in\cS^{\mathrm{orb}}$, $\dim P=n$. An \emph{orbifold gluing decomposition} $(Q_i,\iota_i)_{i=1}^m$ of $P$ satisfies (G2) of Definition~\ref{def_glue} and:
\begin{enumerate}[(G0$'$)]
\item[(G0$'$)] each $Q_i\in\cS^{\mathrm{orb}}$ has dimension $n$;
\item[(G1$'$)] each $\iota_i:|Q_i|\to|P|$ is a \emph{liftable smooth embedding}: it is injective, and near every point it has a lift to charts which is a smooth embedding of manifolds with corners; it preserves orientation;
\item[(G3$'$)] for every $i$ and every face $F$ of $Q_i$, either (a) $\iota_i\circ i_{Q_i}(F)\subseteq B_P$, or (b) there are $j\ne i$, a face $F'$ of $Q_j$ and a diffeomorphism of orbifolds with corners $\psi:F\to F'$ with $\iota_j\circ i_{Q_j}\circ\psi=\iota_i\circ i_{Q_i}|_F$.
\end{enumerate}
Here $B_P:=i_P(|\partial P|)$ is the set of points of positive depth, (G2) is read in $|P|$, and orientations are those of the lifts. No condition is imposed on the local groups of the pieces. We write $\fG^{\mathrm{orb}}$ for the class of these data, $\cG^{\mathrm{orb}}(Y)$ for the corresponding relations, and $GC^{G}_{\mathrm{orb}}(Y)$ for the quotient complex, as in Definition~\ref{def_family}.
\end{definition}

(G3$'$) is (G3) written without inverting $\iota_j\circ i_{Q_j}|_{F'}$, which need not be injective when faces are not embedded. Every gluing datum of Definition~\ref{def_glue} is an orbifold gluing datum, so $\fG\subseteq\fG^{\mathrm{orb}}$. The pieces may have smaller local groups than $P$, but only on their boundaries:

\begin{lemma}\label{lem_orb_interior}
Let $\iota:Q\to P$ satisfy \textup{(G1$'$)}, $q\in Q$ and $x=\iota(q)$, with local groups $H_q$ and $G_x$.
\begin{enumerate}[(a)]
\item $|H_q|\le|G_x|$.
\item If $q\in\inter Q$, then $x\in\inter P$, $H_q\cong G_x$, and lifts of $\iota$ near $q$ are equivariant for such an isomorphism.
\end{enumerate}
\end{lemma}

\begin{proof}
Composing with embeddings of charts and shrinking, as in the proof of \cite[Lemma on orbifold morphisms, (b)]{Yu1}, we may take a lift $\tilde\iota:\tilde U\to\tilde V$ near $q$ which is an embedding, with $\tilde U$ an $H_q$-invariant neighbourhood of a lift $\tilde q$ in a chart centred at $\tilde q$, and $\tilde V$ a chart centred at $\tilde x=\tilde\iota(\tilde q)$. Points of $\tilde V$ with the same image in $|P|$ lie in one $G_x$-orbit.

(a) $H_q$ acts effectively, so points of $\tilde U$ with trivial isotropy accumulate at $\tilde q$, also when $q$ has positive depth. A free $H_q$-orbit gives $|H_q|$ points with distinct images under the injective $\tilde\iota$, all with the same image in $|P|$, hence in one $G_x$-orbit; so $|H_q|\le|G_x|$.

(b) Since $\tilde q$ is interior, $\tilde\iota(\tilde U)$ contains a neighbourhood of $\tilde x$ homeomorphic to $\R^n$, so $\tilde x$ has depth $0$ by invariance of domain, $x\in\inter P$, and $\tilde\iota$ is a diffeomorphism onto an open neighbourhood of $\tilde x$, which contains a $G_x$-invariant ball $B$. A point of $B$ with trivial isotropy has $|G_x|$ translates in $B$, all mapped to the same point of $|P|$; their preimages under $\tilde\iota$ map to the same point of $|Q|$ by injectivity of $\iota$, so they form part of one $H_q$-orbit, and $|G_x|\le|H_q|$. Applying \cite[Lemma on orbifold morphisms, (d)]{Yu1} to the restriction of $\iota$ to an open suborbifold of $\inter Q$ around $q$, on which it is a local diffeomorphism, gives a homomorphism $\lambda:H_q\to G_x$ with $\tilde\iota\circ\gamma=\lambda(\gamma)\circ\tilde\iota$. It is injective, since $\lambda(\gamma)=1$ gives $\tilde\iota\gamma=\tilde\iota$, hence $\gamma=1$ by injectivity of $\tilde\iota$ and effectiveness. By (a) and $|G_x|\le|H_q|$ it is an isomorphism.
\end{proof}

For instance, a smoothly embedded simplex $\Delta^4\to\mathbb C^2/\mu_p$ with a vertex at the singular point satisfies (G1$'$): its lift is the inclusion of a narrow simplicial cone $\sigma$ with $\sigma\cap g\sigma=\{0\}$ for $g\ne1$, which exists since $\mu_p$ acts freely on the unit sphere.

If $P$ and all $Q_i$ are manifolds with embedded faces, (G3$'$) implies (G3), since $\iota_j\circ i_{Q_j}|_{F'}$ is then injective; so on the data of Definition~\ref{def_glue}, $\fG^{\mathrm{orb}}$ restricts to $\fG$. If $P$ is a manifold, the pieces are manifolds by Lemma~\ref{lem_orb_interior}(a), but their faces need not be embedded.

\subsection{The abstract theorem for orbifold sources}
\begin{proposition}\label{prop_abstract_orb}
Theorem~\ref{thm_abstract} holds with $\cS$ replaced by $\cS^{\mathrm{orb}}$ and $\Cm$ by $C^{\mathrm{orb}}$, for relation families of decomposition-type data with sources in $\cS^{\mathrm{orb}}$.
\end{proposition}

\begin{proof}
The proof of Theorem~\ref{thm_abstract} in Sections~\ref{sec_inj}, \ref{sec_cut} and \ref{sec_proof} uses the following properties of sources, all of which hold in $\cS^{\mathrm{orb}}$ when read in charts. Since the depth strata of $|Q|$ need not be manifolds, regularity is read in charts throughout: ``$(h,c)$ regular for $\fd$'' in Definition~\ref{def_family}, ``$c$-transverse'', and the regular values in Lemmas~\ref{lem_pieces} and \ref{lem_generic} refer to the invariant lifts on all corner strata of all charts. By the chart form of (J2) below, $C^\pitchfork_c$ is then a subcomplex, as in Lemma~\ref{lem_cutops}.
\begin{itemize}
\item \emph{Pieces} (Lemma~\ref{lem_pieces}). For $h:|Q|\to\R$ smooth and $c$ a regular value of the invariant lifts of $h$ on all corner strata of all charts, $Q_\pm$ and $L_Q$ are compact orbifolds with corners, with $\partial Q_-=(\partial Q)_-\sqcup L_Q$, by \cite[Lemma on bisection pieces]{Yu1}. They are effective: a local group element fixing a neighbourhood in the hypersurface $\{\tilde h=c\}$ preserves the covector $d\tilde h$, so its differential is unipotent of finite order, hence trivial \cite[proof of the Lemma on orbifold morphisms, (c)]{Yu1}. The statement in Lemma~\ref{lem_pieces} that faces are embedded is not used in the abstract theorem. The property (J2) of the boundary is used in its chart form: in a linear chart, $i_Q$ lifts to the inclusions of boundary hyperplanes, which map points of depth $j$ of a hyperplane to points of depth $j+1$; on underlying spaces this fails, see the example $v\mapsto|v|$ in \cite[proof of the Lemma on the bisection prism]{Yu1}.
\item \emph{Genericity} (Lemma~\ref{lem_generic}). Sard's theorem applied to the invariant lifts in finitely many charts covering $Q$.
\item \emph{Prisms and stars} (Lemmas~\ref{lem_prism}, \ref{lem_star}). $Q\times I\in\cS^{\mathrm{orb}}$, with the boundary formula of \cite[Section on operations]{Yu1}.
\item \emph{Integration} (Propositions~\ref{prop_stokes}, \ref{prop_poly_inj}). For smooth $f$, $f^*\omega$ is an invariant form on charts, and Stokes' theorem holds on orbifolds with corners with the weights $1/|G|$ \cite[Remark on injectivity via Stokes]{Yu1}.
\item \emph{Basis and orientation reversal} (Section~\ref{sec_setting}). \cite[Lemma on the basis of $C^g_*$]{Yu1}.
\item \emph{Degree $0$} (Theorem~\ref{thm_abstract}(a)). A $0$-dimensional effective orbifold is a finite set of points.
\end{itemize}
The remaining steps (Lemmas~\ref{lem_cutops}, \ref{lem_cutrel}, \ref{lem_MV} and Propositions~\ref{prop_abs_poly}, \ref{prop_abs_mfd}) are formal consequences of these and of the axioms. Triangulations of sources, Lemma~\ref{lem_adm}, the topological properties of $B_Q$ used in Lemma~\ref{lem_bdry}, Lemma~\ref{lem_dec}, which enters Proposition~\ref{prop_additive} through Lemma~\ref{lem_cov}, and Section~\ref{sec_point} enter only through (A3), (A5) and (A6) for the specific families, which are proved separately below.
\end{proof}

\subsection{The axioms for orbifold gluing}
\begin{proposition}\label{prop_orb_A15}
$\fG^{\mathrm{orb}}\times\{\text{all }f\}$ satisfies \textup{(A1)}--\textup{(A5)}. Moreover, for every $Q\in\cS^{\mathrm{orb}}$, every smooth $h:|Q|\to\R$ (lifting to invariant smooth functions on charts) and every regular value $c$ in the chart sense of Proposition~\ref{prop_abstract_orb}, the bisection datum, with pieces $Q_\pm$ and their inclusions, lies in $\fG^{\mathrm{orb}}$.
\end{proposition}

\begin{proof}
(A1) holds because the data do not involve $f$, and $f\circ\iota_i$ is smooth since compositions of liftable maps are liftable \cite[Lemma on orbifold morphisms, (b)]{Yu1}. For bisections, the inclusions of $Q_\pm$ preserve local groups and lift to inclusions of half-spaces in charts, so (G1$'$) holds; (G2) is clear; the faces of $Q_\pm$ are the components of $(\partial Q)_\pm$, which are outer, and of $L_Q$, which are glued to themselves in the other piece by $\psi=\mathrm{id}$. In particular (A2) holds.

(A3): let $(P;(Q_i,\iota_i))\in\fG^{\mathrm{orb}}$ and $(h,c)$ regular for it. We go through the proof of Lemma~\ref{lem_cutclosure} for the cut data $\fd_\pm$ and $\fd_0$. (G0$'$): $Q_{i,\pm}$ and $L_{Q_i}$ lie in $\cS^{\mathrm{orb}}$ by the item \emph{Pieces} of Proposition~\ref{prop_abstract_orb}. (G1$'$): near a point, a lift $\tilde\iota$ which is an embedding restricts to embeddings $\tilde U\cap\{\tilde h\tilde\iota\le c\}\to\tilde V\cap\{\tilde h\le c\}$ and $\tilde U\cap\{\tilde h\tilde\iota=c\}\to\tilde V\cap\{\tilde h=c\}$ into slice charts of $P_-$ and $L_P$, and similarly for $P_+$; injectivity is inherited. For $\fd_0$, $L_{Q_i}$ and $L_P$ are oriented with the outward normal first, and at a point of $\inter L_{Q_i}$, which is interior in $Q_i$, $\tilde\iota$ is a local diffeomorphism by Lemma~\ref{lem_orb_interior}(b) mapping outward normals to outward normals, so orientations are preserved; this replaces the use of Lemma~\ref{lem_dec}(a). (G2): $\inter Q_{i,\pm}\subseteq\inter Q_i$ and $\iota_i^{-1}(L_P)=L_{Q_i}$. (G3$'$): an outer face $F$ of $Q_i$ gives faces of $Q_{i,\pm}$ mapping into $B_{P_\pm}$ and faces of $L_{Q_i}$ mapping into $B_{L_P}$, since depth in $L_P$ equals depth in $P$; the new faces $L_{Q_i}$ of $Q_{i,\pm}$ map into $L_P\subseteq B_{P_\pm}$. For a glued face, $\iota_j\circ i_{Q_j}\circ\psi=\iota_i\circ i_{Q_i}$ gives $h\iota_j\psi=h\iota_i$ on $F$, so $\psi$ maps $F_\pm$ onto $F'_\pm$ and $F\cap L$ onto $F'\cap L$, component by component, and restricts to diffeomorphisms of orbifolds with corners, since $c$ is a regular value on the chart strata of $F$ and $F'$ by the chart form of (J2). Embeddedness of faces is not used.

(A4): $\partial(Q\times I)=\partial Q\times I\sqcup Q\times\{0,1\}$; the faces $F\times I$ are glued by $\psi\times\mathrm{id}$ and the faces $Q_i\times\{0,1\}$ are outer, as in Lemma~\ref{lem_products}.

(A5): Let $\alpha$ be an $n$-form on $P$. In a chart $(\tilde V,G)$ the points with trivial isotropy form an open set $\tilde V_{\mathrm{free}}$ with null complement, and $\tilde V_{\mathrm{free}}\to\tilde V_{\mathrm{free}}/G$ is a $|G|$-sheeted covering; so the orbifold integral $\int_P\alpha$, with the weight $1/|G|$, equals the integral of $\alpha$ over the manifold $|P|_{\mathrm{reg}}$ of interior points with trivial isotropy, whose complement in $|P|$ is null. By Lemma~\ref{lem_orb_interior}(b), $\iota_i$ maps $\inter Q_i\cap|Q_i|_{\mathrm{reg}}$ diffeomorphically and orientation preservingly onto an open subset of $|P|_{\mathrm{reg}}$, and $\iota_i(\inter Q_i\setminus|Q_i|_{\mathrm{reg}})$ is null, since it lies in the singular set of $P$ by Lemma~\ref{lem_orb_interior}(b). The image $\iota_i(B_{Q_i})$ is null, being locally a smooth image of finitely many $(n-1)$-dimensional manifolds with corners. So $\int_{Q_i}\iota_i^*\alpha=\int_{\iota_i(Q_i)}\alpha$, as in Lemma~\ref{lem_cov}, and Proposition~\ref{prop_additive} follows as before from (G2).
\end{proof}

\subsection{The point}
\begin{proposition}\label{prop_orb_point}
For every $Q\in\cS^{\mathrm{orb}}$ of dimension $k\ge1$, $[Q,c]\in\operatorname{span}\cG^{\mathrm{orb}}(\mathrm{pt})$. Hence $GC^{G}_{\mathrm{orb},k}(\mathrm{pt})=0$ for $k\ge1$ at chain level, and \textup{(A6)} holds for $\fG^{\mathrm{orb}}$.
\end{proposition}

\begin{proof}
Write $n=k$.

\textit{Step 1: invariant simplicial fans.} Let $G$ be a finite group acting on $\R^n=\R^d\times\R^{n-d}$ by $g(y,l)=(P_gy,A_gl)$ with permutation matrices $P_g$ and orthogonal $A_g$, so orthogonally for the standard inner product, preserving the orientation. Choose a finite $G$-invariant set of hyperplanes through $0$ containing these coordinate hyperplanes, with normals spanning $\R^n$. Its arrangement is a $G$-invariant complete fan of pointed polyhedral cones, and $\R^n_d:=\{y_1,\dots,y_d\ge0\}$ is a union of its cones. Subdivide it once barycentrically, the ray of a cone $\tau$ being spanned by the sum of the unit vectors on its extremal rays; this choice is $G$-equivariant, since $G$ is orthogonal. The result is a $G$-invariant complete simplicial fan $\Sigma$ whose cones correspond to chains of cones of the arrangement. If the rays $b(s)$ and $b(gs)$ span a common cone of $\Sigma$, then $s$ and $gs$ are comparable in such a chain and have the same dimension, so $s=gs$. Hence, if $x$ and $gx$ lie in a cone $\sigma\in\Sigma$, every ray of the carrier cone of $x$ is fixed by $g$, and $gx=x$: \emph{every closed cone of $\Sigma$ maps injectively to $\R^n/G$}. Moreover, if $\sigma\ne\sigma'=g\sigma''$ share a facet $\tau$ and $\sigma''=\sigma$, then $\tau$ and $g^{-1}\tau$ both lie in $\sigma$, so $g$ fixes $\tau$, hence the hyperplane $\operatorname{span}\tau$, pointwise; as $g\ne1$ is orthogonal, it is the reflection in this hyperplane, which reverses the orientation. So \emph{a facet of a cone $\sigma$ is never shared with a translate $g\sigma\ne\sigma$}. Finally, if $g\sigma=\sigma$ for an $n$-dimensional cone, then $g$ fixes $\sigma$ pointwise by injectivity, so $g=1$: \emph{$G$ acts freely on the $n$-dimensional cones}.

\textit{Step 2: cutting into linear charts.} By \cite[Doubling lemma]{Yu1}(1), with the $l$-coordinates orthonormal for a $G_j$-invariant inner product, and compactness choose finitely many linear charts $(\tilde U_j,G_j,\varphi_j)$, with $G_j$ acting as in Step 1 on $\tilde U_j\subseteq \R^n_{d_j}$, whose images $U_j$ cover $|Q|$, compact sets $K_j\subseteq U_j$ covering $|Q|$, and smooth invariant functions $\chi_j:|Q|\to[0,1]$, obtained by averaging over $G_j$, with $\chi_j=1$ on $K_j$ and closed support in $U_j$. Fix fans $\Sigma_j$ as in Step 1. Choose $t_1,t_2,\dots\in(0,1)$ successively such that $t_j$ is a regular value of the lifts $\tilde\chi_j$ of $\chi_j$ to the charts $\tilde U_a$ on every manifold $\tau^\circ\cap\{\tilde\chi_i=t_i,\ i\in S\}$, where $S\subseteq\{1,\dots,j-1\}$, $\tau^\circ$ runs over the open cones of the fans $\Sigma_a$ (including $\{0\}$) in the charts $\tilde U_a$, and over the open corner strata of the charts; by induction these are finitely many manifolds, and Sard's theorem applies. Bisect $Q$ along $\chi_1=t_1$, then the piece $\{\chi_1\le t_1\}$ along $\chi_2=t_2$, and so on. By Proposition~\ref{prop_orb_A15} these are orbifold gluing relations, and they give $[Q,c]\equiv\sum_a[Q_a,c]$ modulo $\operatorname{span}\cG^{\mathrm{orb}}(\mathrm{pt})$, where $Q_1=\{\chi_1\ge t_1\}\subseteq U_1$, $Q_2=\{\chi_1\le t_1,\ \chi_2\ge t_2\}\subseteq U_2$, and so on; the last piece lies in $K_N\subseteq U_N$, since $\chi_j=1>t_j$ on $K_j$. Each $Q_a$ is the quotient $G_a\ltimes\tilde V_a$ of the full preimage $\tilde V_a\subseteq\tilde U_a$, a $G_a$-invariant region, compact since the quotient map is proper, cut out by finitely many of the functions $\tilde\chi_i\le t_i$, $\tilde\chi_i\ge t_i$.

\textit{Step 3: manifold pieces.} Fix $a$, write $G=G_a$, and choose one representative $\sigma$ in each $G$-orbit of $n$-dimensional cones of $\Sigma_a$ contained in $\R^n_{d_a}$. Put $Q_{a,\sigma}:=\sigma\cap\tilde V_a$, possibly empty, with the orientation of the chart $\tilde U_a$, and let $\iota_{a,\sigma}$ be the inclusion followed by the quotient map. At a point of an open cone face $\tau^\circ$ of dimension $l$, $\sigma$ is locally given by $n-l$ independent linear inequalities, and the active cut functions have independent differentials on $\operatorname{span}\tau$ by Step 2; and $0\notin\{\tilde\chi_i=t_i\}$. So $Q_{a,\sigma}$ is a compact manifold with corners, and its faces are embedded, since the local boundary components at a point lie on distinct labelled hypersurfaces (facet hyperplanes of $\sigma$, cut hypersurfaces) and labels are constant on faces. By Step 1, $\iota_{a,\sigma}$ is injective; its lift is an inclusion, so (G1$'$) holds. The images cover $|Q_a|$, and interiors are disjoint, as distinct representatives lie in distinct orbits, so (G2) holds. For (G3$'$), let $F$ be a face of $Q_{a,\sigma}$. If $F$ lies in a cut hypersurface or in a coordinate hyperplane $\{y_i=0\}$, it maps to $B_{Q_a}$. Otherwise $F$ lies in a facet $\tau$ of $\sigma$ whose relative interior meets the interior of $\R^n_{d_a}$, since $\Sigma_a$ refines the coordinate hyperplanes; then $\tau$ is a facet of another $n$-dimensional cone $\sigma'$, which lies in $\R^n_{d_a}$ because its interior meets the interior of $\R^n_{d_a}$ and $\Sigma_a$ refines the coordinate hyperplanes, so $\sigma'=g\sigma''$ with $\sigma''$ a representative, $\sigma''\ne\sigma$ by Step 1, and $\psi=g^{-1}$ maps $F$ diffeomorphically onto a face of $Q_{a,\sigma''}$, as $g^{-1}$ is linear and preserves $\tilde V_a$ and the cut hypersurfaces. So $(Q_a;(Q_{a,\sigma},\iota_{a,\sigma})_\sigma)\in\fG^{\mathrm{orb}}$, and $[Q_a,c]\equiv\sum_\sigma[Q_{a,\sigma},c]$. As a check of the weights: each orbit contains $|G|$ cones by Step 1, so $\int_{Q_a}\alpha=\frac1{|G|}\int_{\tilde V_a}\tilde\alpha=\sum_\sigma\int_{Q_{a,\sigma}}\tilde\alpha$ for every $n$-form $\alpha$, as (A5) requires.

\textit{Step 4.} Each $Q_{a,\sigma}$ is a source in the sense of Definition~\ref{def_source}. By Corollary~\ref{cor_point}, $[Q_{a,\sigma},c]\in\operatorname{span}\cG(\mathrm{pt})\subseteq\operatorname{span}\cG^{\mathrm{orb}}(\mathrm{pt})$. Hence $[Q,c]\in\operatorname{span}\cG^{\mathrm{orb}}(\mathrm{pt})$.
\end{proof}

\begin{theorem}\label{thm_orb}
For every compact polyhedron $Y\subseteq\R^K$ and every smooth manifold $Y$ without boundary, $p_*:H_*(Y;\Q)\to H_*(GC^{G}_{\mathrm{orb}}(Y))$ is an isomorphism in all degrees.
\end{theorem}

\begin{proof}
Proposition~\ref{prop_abstract_orb}, with the axioms verified in Propositions~\ref{prop_orb_A15} and \ref{prop_orb_point}.
\end{proof}

\begin{remark}\label{rem_orb}
(1) For orbifold sources, $D$ cannot be dropped without a replacement \cite[Remark on the role of degenerate chains]{Yu1}, but Theorem~\ref{thm_orb} shows that the gluing relation replaces it, for the targets of this paper and smooth maps. The cycles $[\mathbb P(1,1,p),c]$ of \cite{Yu1}, which are not boundaries without a third relation, vanish in $GC^G_{\mathrm{orb}}(\mathrm{pt})$ at chain level, by Proposition~\ref{prop_orb_point}: near the singular point $\mathbb P(1,1,p)$ is glued from manifold pieces. Orbifold targets are treated in Section~\ref{sec_orbtarget}; we do not treat the gluing relation for the more general smooth spaces of \cite{Yu1}, nor for orbifold sources with continuous maps.

(2) It is essential that pieces may have smaller local groups than $P$. If (G1$'$) is strengthened to require that the lifts of $\iota_i$ be equivariant for isomorphisms of local groups at all points, then every $p$-point of $P$ (in the sense of \cite[Remark on the role of degenerate chains]{Yu1}) lies in the interior of exactly one piece, since at a point of positive depth the local group fixes a non-zero vector and the local group $\mu_p$ of a $p$-point does not; interior $p$-points map to $p$-points with the same sign, since a real linear map normalising the scalar action of $\mu_p$ is complex linear or antilinear and has positive determinant. Hence the signed count $\phi_p$ of $p$-points is additive on such relations, and, since $\phi_p\circ\partial=0$ for $p\ge5$ \cite[Remark on the role of degenerate chains]{Yu1}, $[\mathbb P(1,1,p),c]$ remains a non-zero class in degree $4$ of the point for $p\ge5$. By Lemma~\ref{lem_orb_interior}, in Definition~\ref{def_glue_orb} the local groups can only become smaller on the boundaries of the pieces; accordingly the $p$-points lie in corners of the pieces of Step 3.
\end{remark}

\section{Integer coefficients}\label{sec_int}
In this section coefficients are $\Z$. We show that the gluing relation computes integral singular homology. Apart from the verification of (A6) for $G$ and one argument in the proof of Lemma~\ref{lem_cutops}, which used that classes with an orientation reversing automorphism vanish, the only step of the proof of Theorem~\ref{thm_abstract} that is not integral is injectivity, which came from integration of forms; we replace it by intersecting with the dual cells of a smooth triangulation, which gives a chain map to simplicial chains.

\subsection{Setting}
Let $\Cm_*(Y;\Z)$ be defined as in Section~\ref{sec_setting} with integer coefficients: the free abelian group on the classes of connected geometric simplices modulo $[-Q,f]=-[Q,f]$. A class with an orientation reversing automorphism now satisfies $2[Q,f]=0$, but need not vanish. Since sources have embedded faces, $\partial^2=0$ still holds: the components of $\partial^2Q$ are the components of $F_a\cap F_b$, labelled by ordered pairs $(a,b)$ with $a\ne b$, and the involution of \cite[Convention on the boundary]{Yu1} maps such a component to the same component of $F_b\cap F_a$ with the label $(b,a)$, reversing orientation, so it fixes no component and the terms cancel in pairs. The relations $\cG(Y)$, the complex $\GCG(Y;\Z)=\Cm_*(Y;\Z)/G(Y;\Z)$ with $G(Y;\Z)=\operatorname{span}_\Z\cG(Y)+\partial\operatorname{span}_\Z\cG(Y)$, and the map $p$ from smooth singular chains $\Csm_*(Y;\Z)$ are defined as before.

\begin{theorem}\label{thm_Z}
For every compact polyhedron $Y\subseteq\R^K$ and every smooth manifold $Y$ without boundary, $p_*:H_*(Y;\Z)\to H_*(\GCG(Y;\Z))$ is an isomorphism in all degrees.
\end{theorem}

Since $\GCG(Y;\Q)=\GCG(Y;\Z)\otimes\Q$ by right exactness and $\Q$ is flat, Theorem~\ref{thm_Z} implies Theorem~\ref{thm_A_G}, with a proof that does not use de Rham's theorem.

\subsection{The dual cell map}
Let $X$ be a smooth manifold without boundary, of dimension $n$, possibly non-compact and non-orientable. Choose a smooth triangulation $K$ of $X$ \cite{Whi,Mun}, locally finite, and let $K'$ be its barycentric subdivision, whose simplices $[b(\tau_0),\dots,b(\tau_r)]$, $\tau_0<\dots<\tau_r$, are again smoothly embedded. For a simplex $\tau$ of $K$ the \emph{dual cell} $D(\tau)$ is the union of the simplices of $K'$ with $\tau_0=\tau$; it has dimension $n-\dim\tau$, and $\partial D(\tau)=\bigcup D(\tau')$, the union over $\tau'>\tau$ with $\dim\tau'=\dim\tau+1$. Moreover $D(\tau)$ is a pseudomanifold with boundary: every simplex of codimension one of $D(\tau)$ not in $\partial D(\tau)$ is a face of exactly two top-dimensional simplices of $D(\tau)$, and every simplex of codimension one in $\partial D(\tau)$ is a face of exactly one, obtained by adding the vertex $b(\tau)$.

\emph{Coorientation.} Every top-dimensional simplex of $D(\tau)$ lies in an $n$-simplex $\sigma\supseteq\tau$ of $K$, and in the affine coordinates of $\sigma$ it lies in the affine subspace $A_\sigma(\tau)$ on which the barycentric coordinates of the vertices of $\tau$ are equal. The linear span $W_\tau$ of the edge vectors of $\tau$ is a complement of the direction of $A_\sigma(\tau)$, and we coorient $D(\tau)$ by the orientation of $W_\tau$ given by the orientation of $\tau$; no orientation of $X$ is used. This is consistent across codimension one faces in the following sense. Let $T,T'$ be top simplices of $D(\tau)$ with a common face $S$ not in $\partial D(\tau)$. If they lie in the same $n$-simplex $\sigma$, they lie in $A_\sigma(\tau)$ and are cooriented by the same $W_\tau$. Otherwise $T\subseteq\sigma$, $T'\subseteq\sigma'$, the last vertex of $S$ is the barycentre of the $(n-1)$-simplex $\rho'=\sigma\cap\sigma'\supseteq\tau$, so the open simplex $S^\circ$ lies in the open simplex $\rho'^\circ$, and $W_\tau\subseteq T\rho'$ is the same subspace in both charts, since the parametrizations of $\sigma$ and $\sigma'$ restrict to that of $\rho'$; at $x\in S^\circ$ it is a complement of $T_xS$ in $T_x\rho'$. Now let $f:Q\to X$, $\dim Q=\dim\tau+1$, be transverse to $S^\circ,T^\circ,T'^\circ$, with $f(q)=x\in S^\circ$. Then $f$ is transverse to $\rho'^\circ$ at $q$, so $H=f^{-1}(\rho')$ is a hypersurface near $q$, and an oriented basis of $W_\tau$ modulo $T_xS$ lifts to a basis $\tilde w$ of $T_qH$. Near $q$, $f^{-1}(T\cup T')$ is an arc through $q$ whose two halves lie on the two sides of $H$, since they map into $\sigma$ and $\sigma'$. Orient each half by $a$ with $(a,\tilde w)$ positive in $T_qQ$, the arc direction followed by the pulled back coorientation, which does not depend on the lift. Then both directions lie on the same side of $H$, so the orientation is continuous along the arc.

\emph{Transversality.} A geometric simplex $[Q,f]$ in $X$ is \emph{$K$-transverse} if $f$ restricted to every open corner stratum of $Q$ is transverse to every open simplex of $K'$. Let $C^\pitchfork_*(X)\subseteq\Cm_*(X;\Z)$ be the span of the $K$-transverse classes, a subcomplex, and call a gluing datum $(\fd,f)$ \emph{$K$-transverse} if $f$ and all $f\circ\iota_i$ are; let $\cR^\pitchfork(X)$ be the span of these relations and their boundaries.

For a $K$-transverse $[Q,f]$ of dimension $k$ and $\dim\tau=k$, the set $f^{-1}(D(\tau))$ is finite and consists of interior points of $Q$ mapped to open top-dimensional simplices of $D(\tau)$, by dimension. Let $\langle f,D(\tau)\rangle\in\Z$ be its signed count, a point $q$ counting $+1$ if $df_q$ maps the orientation of $T_qQ$ to the coorientation of $D(\tau)$. Put
\begin{equation*}
I[Q,f]:=\sum_{\dim\tau=k}\langle f,D(\tau)\rangle\,\tau\in C^{\mathrm{simp}}_k(K;\Z),
\end{equation*}
a finite sum, since $f(Q)$ is compact and $K$ is locally finite.

\begin{lemma}\label{lem_dual}
\begin{enumerate}[(a)]
\item $I$ is well defined on $C^\pitchfork_*(X)$, and vanishes on classes with an orientation reversing automorphism.
\item $\partial\circ I=I\circ\partial$.
\item $I(\cR^\pitchfork(X))=0$.
\item For every finite simplicial chain $c$ of $K$, regarded as a smooth singular chain, there is $\varepsilon>0$ such that $(\mathbf F_\beta)_*p(c)$ is $K$-transverse and $I((\mathbf F_\beta)_*p(c))=c$ for almost all $\beta$ with $|\beta|<\varepsilon$, where $\mathbf F:X\times B\to X$, $B$ a ball around $0$ in some $\R^m$, is a \emph{submersion family}: $\mathbf F(x,0)=x$, and $\beta\mapsto\mathbf F(x,\beta)$ is a submersion at $\beta=0$ for every $x$, hence near $\beta=0$ uniformly on compact sets; for instance $\mathbf F(x,\beta)=r(x+\beta)$ for a tubular retraction $r$ of $X\subseteq\R^K$.
\end{enumerate}
\end{lemma}

\begin{proof}
(a) The count is invariant under isomorphism and changes sign with the orientation. If $a$ is an orientation reversing automorphism with $f\circ a=f$, then $a$ permutes $f^{-1}(D(\tau))$ and reverses the sign of each point, so the count $S$ satisfies $S=-S$ in $\Z$, and $S=0$.

(b) Fix $\rho$ with $\dim\rho=k-1$, $k=\dim Q$. By transversality, $f^{-1}(D(\rho))$ is a compact one-dimensional manifold with boundary, piecewise smooth: the preimages of the open top simplices of $D(\rho)$ are open arcs; at preimages of interior codimension one simplices of $D(\rho)$ two arcs meet, since $D(\rho)$ is a pseudomanifold with boundary, and their orientations defined below match by the consistency of the coorientation; strata of codimension $\ge2$ are not met. Its ends are the points of $(f\circ i_Q)^{-1}(D(\rho))$ and the points of $f^{-1}(D(\tau))$, $\tau>\rho$, $\dim\tau=k$. With the outward normal first for $\partial Q$ and the incidence numbers $[\tau:\rho]$ of $K$, also outward normal first, orient the arcs so that at a point of $f^{-1}(D(\rho)^\circ_{\mathrm{top}})$ the orientation of $T_qQ$ is the arc direction followed by the pullback of the coorientation of $D(\rho)$; then, counting final points $+1$ and initial points $-1$, an end on $\partial Q$ counts as in $\langle f\circ i_Q,D(\rho)\rangle$, and an end on $D(\tau)$ counts $-[\tau:\rho]$ times its sign in $\langle f,D(\tau)\rangle$, since there the arc enters the top simplex of $D(\rho)$ adjacent to $D(\tau)$ in the direction $b(\rho)-b(\tau)$ modulo the tangent space of $D(\tau)$, which lies in $W_\tau$ and points like the outward normal of $\rho$ in $\tau$. The signed count of the ends of a compact oriented one-manifold is zero, which gives $\langle f\circ i_Q,D(\rho)\rangle=\sum_\tau[\tau:\rho]\langle f,D(\tau)\rangle$, i.e. the coefficient of $\rho$ in $I(\partial x)$ and in $\partial I(x)$ agree.

(c) Let $(P;(Q_i,\iota_i))$ be a $K$-transverse gluing datum with map $f$. The interfaces $f\iota_i(B_{Q_i})$ have dimension $k-1$ and are transverse to $D(\tau)$, of codimension $k$, hence disjoint from it. So every point of $f^{-1}(D(\tau))$ lies in $\iota_i(\inter Q_i)$ for exactly one $i$, where $\iota_i$ is an orientation preserving local diffeomorphism, with the same sign. Hence $I[P,f]=\sum_iI[Q_i,f\iota_i]$, and $I(\partial\cdot)=\partial I(\cdot)$ by (b).

(d) Let $\tau'$ be a $k$-simplex of $c$, with parametrization $\sigma_{\tau'}$. If $\tau'\neq\tau$, $\dim\tau=k$, then $\tau'\cap D(\tau)=\emptyset$, since a simplex $[b(\tau_0),\dots]$ of $D(\tau)$ meets $\tau'$ only in the face spanned by those $b(\tau_i)$ with $\tau_i\le\tau'$, and $\tau=\tau_0\le\tau'$ forces $\tau=\tau'$; by compactness this persists for small $\beta$. For $\tau'=\tau$ we have $\tau\cap\partial D(\tau)=\emptyset$ and $\partial\tau\cap D(\tau)=\emptyset$, so there is a $C^0$-neighbourhood $\mathcal N$ of $\sigma_\tau$ such that $g(\Delta^k)\cap\partial D(\tau)=\emptyset$ and $g(\partial\Delta^k)\cap D(\tau)=\emptyset$ for $g\in\mathcal N$. For two smooth $g_0,g_1\in\mathcal N$ transverse to the open simplices of $D(\tau)$, $\langle g_0,D(\tau)\rangle=\langle g_1,D(\tau)\rangle$: join them by the straight homotopy in $\R^K$ followed by $r$, constant near the ends, which stays in $\mathcal N$ if $g_0,g_1$ are close enough to $\sigma_\tau$ (we shrink $\mathcal N$ by a factor depending on the Lipschitz constant of $r$), perturb it by $\mathbf F(\cdot,\varphi(s)\gamma)$ with small generic $\gamma$ and a bump function $\varphi$ that vanishes only near $s=0,1$, where the homotopy equals $g_0$, $g_1$ and is already transverse, which makes it transverse on all corner strata of $\Delta^k\times I$, and apply the argument of (b), which now has no ends on $\partial\Delta^k\times I$ or on $\partial D(\tau)$. To compute the count, choose an $n$-simplex $\sigma\supseteq\tau$, a generic point $z\in\inter\sigma$ and $g_\delta(y)=(1-\delta)y+\delta z$ in the affine coordinates of $\sigma$. For small $\delta$, $g_\delta(\Delta^k)$ meets $A_\sigma(\tau)$ in a single point of an open top simplex of $D(\tau)$, transversally, with tangent space $W_\tau$, so with sign $+1$ by the choice of coorientation; it misses $\partial D(\tau)$, and $g_\delta(\partial\Delta^k)$ misses $D(\tau)$. Finally, $K$-transversality of $(\mathbf F_\beta)_*p(c)$ for almost all $\beta$ is parametric transversality for the countably many open simplices of $K'$.
\end{proof}

\begin{proposition}\label{prop_Z_inj}
For every smooth manifold $X$ without boundary and every compact polyhedron $Y\subseteq\R^K$, $p_*:H_k(\,\cdot\,;\Z)\to H_k(\GCG(\,\cdot\,;\Z))$ is injective.
\end{proposition}

\begin{proof}
Let $X$ be a manifold. Fix a total order of the vertices of $K$, so that simplicial chains are smooth singular chains. Every class of $H_k(X;\Z)$ is represented by a finite simplicial cycle $c$ of $K$, via the isomorphism induced by $C^{\mathrm{simp}}_*(K;\Z)\to\Csm_*(X;\Z)$, which uses that the simplices of $K$ are smooth and that $\Csm_*\to C_*$ is a quasi-isomorphism \cite{Hat,Lee}. Suppose $p(c)=\partial w+r$ with $w\in\Cm_{k+1}(X;\Z)$ and $r\in\cR(X)$, written with finitely many relation data. Only finitely many chains are involved, with images in a compact set on which $\mathbf F_\beta$ is defined for small $\beta$ (for non-compact $X$ use $\mathbf F(x,\beta)=r(x+\eta(x)\beta)$ with $\eta>0$ small). By (A1), $(\mathbf F_\beta)_*$ is a chain map preserving $\cR$, so $p(c)_\beta=\partial w_\beta+r_\beta$. For almost all small $\beta$, $w_\beta$, the moved relation data and $p(c)_\beta$ are $K$-transverse, by parametric transversality, and by Lemma~\ref{lem_dual}
\begin{equation*}
c=I(p(c)_\beta)=\partial I(w_\beta)+I(r_\beta)=\partial I(w_\beta).
\end{equation*}
So $c$ is a simplicial boundary.

Let $Y\subseteq\R^K$ be a compact polyhedron. It is a Euclidean neighbourhood retract, so it has an open neighbourhood $U$ with a continuous retraction onto $Y$, and $i_*:H_*(Y;\Z)\to H_*(U;\Z)$ is injective. By naturality, $p_U\circ i_*=i_*\circ p_Y$ on homology. If $p_Y(a)=0$, then $p_U(i_*a)=0$, hence $i_*a=0$ by the manifold case, and $a=0$.
\end{proof}

\subsection{Proof of Theorem~\ref{thm_Z}}
Injectivity is Proposition~\ref{prop_Z_inj}. Surjectivity follows from the proofs of Propositions~\ref{prop_abs_poly} and \ref{prop_abs_mfd} with integer coefficients, with three changes.
\begin{itemize}
\item In Step 5 of Lemma~\ref{lem_MV}, injectivity on the halves $Y_\pm$, which are compact polyhedra, is Proposition~\ref{prop_Z_inj} instead of Proposition~\ref{prop_poly_inj}.
\item In Lemma~\ref{lem_star}(b), injectivity in degree $0$ follows from the augmentation, the signed count of points, which vanishes on boundaries of $1$-chains and on the (trivial) gluing relations in degree $0$.
\item (A6) over $\Z$: by the proof of Corollary~\ref{cor_point}, $[Q,c]\equiv\sum_i\varepsilon_i[\Delta^k,c]$ modulo $\operatorname{span}\cG(\mathrm{pt})$. For $k\ge1$, $[\Delta^k,c]\in\operatorname{span}\cG(\mathrm{pt})$: the hyperplane through $v_2,\dots,v_k$ and the midpoint of $[v_0,v_1]$ divides $\Delta^k$ into two simplices $\sigma_1,\sigma_2$, with the orientations induced from $\Delta^k$, which form a Whitehead triangulation, so $[\Delta^k,c]-[\sigma_1,c]-[\sigma_2,c]\in\cG(\mathrm{pt})$ by Lemma~\ref{lem_glue}(e); the transposition of $v_0,v_1$ maps $\sigma_1$ onto $\sigma_2$ reversing orientation and commutes with $c$, so $[\sigma_2,c]=-[\sigma_1,c]$.
\end{itemize}
All other steps (Lemmas~\ref{lem_cutrel}, \ref{lem_prism}, Proposition~\ref{prop_amb}) are integral. In Lemma~\ref{lem_cutops} the formulas are defined on isomorphism classes and are compatible with $(-Q)_\pm=-(Q_\pm)$, so over $\Z$ they are well defined without the argument for classes with an orientation reversing automorphism; for such a class both sides are $2$-torsion, consistently with $2[Q,f]=0$. \qed

\begin{remark}\label{rem_Z}
(1) Lemma~\ref{lem_dual} gives a proof of injectivity that does not use de Rham's theorem, also over $\Q$.

(2) The decomposition relation and the readings $G^W$, $G^T$ of Section~\ref{ssec_readings} satisfy Lemma~\ref{lem_dual}(c) as well, since their pieces are injective with local degree $+1$ at interior points; and (A6) over $\Z$ holds for them because $[Q,c]\in\operatorname{span}\cG(\mathrm{pt})$ at chain level and $\cG\subseteq\cG^W\subseteq\cD$, $\cG^W\subseteq\cG^T$; so Corollary~\ref{cor_GE} holds over $\Z$. In particular, for $G^T$ injectivity, and injectivity on the halves in Step 5 of Lemma~\ref{lem_MV}, follow from Lemma~\ref{lem_dual}(c), which needs only (A1) for $G^T$, without the area formula sketch of Remark~\ref{rem_T}; with (A2)--(A4) and (A6), verified in the proof of Theorem~\ref{thm_readings}, the surjectivity argument above applies, so the comparison theorem holds for $G^T$ over $\Z$ and over $\Q$.

(3) Not covered over $\Z$: Remark~\ref{rem_lowdim} ($[Q\times I,c]$ is only $2$-torsion there), Proposition~\ref{prop_twopiece}(b), which uses $\Cm_k(\mathrm{pt})=0$, the division by $2$ in Proposition~\ref{prop_structure}, and the Lipschitz maps of Section~\ref{ssec_lip}, which cannot be made transverse.

(4) Orbifold sources are treated in Subsection~\ref{ssec_Z_orb}.
\end{remark}

\subsection{Orbifold sources}\label{ssec_Z_orb}
For the orbifold gluing relation of Section~\ref{sec_orb} faces need not be embedded, and the components of $\partial^2Q$ fixed by the involution of \cite[Convention on the boundary]{Yu1} give terms that are only $2$-torsion. Following \cite[Remark on stratified sources, (4)]{Yu1}, which proposes to impose the vanishing of classes with an orientation reversing automorphism as a separate relation over $\Z$, as in \cite{lip}, we divide by these classes.

\begin{definition}\label{def_Z_orb}
Let $C^{\mathrm{orb}}_*(Y;\Z)$ be the abelian group generated by the isomorphism classes of geometric simplices $(Q,f)$ in $Y$ with $Q\in\cS^{\mathrm{orb}}$, modulo $[-Q,f]=-[Q,f]$ and $[Q_1\sqcup Q_2,f_1\sqcup f_2]=[Q_1,f_1]+[Q_2,f_2]$, with the homomorphism $\partial[Q,f]=[\partial Q,f\circ i_Q]$ ($\partial^2\ne0$ in general, see Lemma~\ref{lem_Z_orb_complex}). Let $R_{\mathrm{or}}(Y)$ be the span of the classes, possibly disconnected, that have an orientation reversing automorphism, and put
\[
\bar C_*(Y):=C^{\mathrm{orb}}_*(Y;\Z)/R_{\mathrm{or}}(Y),\qquad \overline{GC}{}^{G}_{\mathrm{orb}}(Y):=C^{\mathrm{orb}}_*(Y;\Z)\big/\big(R_{\mathrm{or}}(Y)+G^{\mathrm{orb}}(Y;\Z)\big),
\]
where $G^{\mathrm{orb}}(Y;\Z)=\operatorname{span}_\Z\cG^{\mathrm{orb}}(Y)+\partial\operatorname{span}_\Z\cG^{\mathrm{orb}}(Y)$. The map $p$ sends a smooth singular simplex $\sigma$ to the class of $[\Delta^n,\sigma]$.
\end{definition}

\begin{lemma}\label{lem_Z_orb_complex}
\begin{enumerate}[(a)]
\item $C^{\mathrm{orb}}_*(Y;\Z)$ is the direct sum of a copy of $\Z$ for each connected class without an orientation reversing automorphism, with a chosen orientation, and a copy of $\Z/2$ for each connected class with one. $R_{\mathrm{or}}(Y)$ is its torsion subgroup, so $\bar C_*(Y)$ is free on the connected classes without such an automorphism, $\bar C_*(W)\subseteq\bar C_*(Y)$ is spanned by a subset of this basis for $W\subseteq Y$, and $\bar C_*(Y)\otimes\Q=C^{\mathrm{orb}}_*(Y)$.
\item $\partial R_{\mathrm{or}}\subseteq R_{\mathrm{or}}$ and $\partial^2C^{\mathrm{orb}}_*(Y;\Z)\subseteq R_{\mathrm{or}}$. Hence $\bar C_*(Y)$ and $\overline{GC}{}^{G}_{\mathrm{orb}}(Y)$ are chain complexes, and $p$ is a chain map to both.
\end{enumerate}
\end{lemma}

\begin{proof}
(a) The relations identify $[-Q]$ with $-[Q]$ for connected $Q$ and are additive over components; a connected class with $Q\cong-Q$ thus generates $\Z/2$, and the others $\Z$. A possibly disconnected class with an orientation reversing automorphism $a$ lies in the torsion subgroup: $a$ permutes the components, on an orbit of odd length $l$ the power $a^l$ is an orientation reversing automorphism of each component, and on an orbit of even length consecutive components are isomorphic with opposite orientations and cancel in pairs. Conversely the torsion subgroup is spanned by the connected classes with such an automorphism. The rest follows, since over $\Q$ these classes vanish \cite[Lemma on the basis of $C^g_*$]{Yu1}.

(b) $\partial$ is a homomorphism, so it preserves torsion. For a generator $(Q,f)$, $\partial^2Q$ is the boundary of $\partial Q$, and the involution $\sigma$ of $\partial^2Q$ reverses orientation and commutes with the map. Components exchanged by $\sigma$ cancel; on a component $C$ with $\sigma(C)=C$, $\sigma|_C$ is an orientation reversing automorphism, so $[C]\in R_{\mathrm{or}}$. This is the proof of \cite[Lemma on the basis of $C^g_*$]{Yu1} without the division by $2$. Hence $\partial(R_{\mathrm{or}}+G^{\mathrm{orb}})\subseteq R_{\mathrm{or}}+G^{\mathrm{orb}}$. In $C^m_*(Y;\Z)$, $p$ is a chain map, since the identity $\partial[\Delta^n,\sigma]=\sum_j(-1)^j[\Delta^{n-1},\sigma d_j]$ of Section~\ref{sec_setting} is integral, and $C^m_*(Y;\Z)\to C^{\mathrm{orb}}_*(Y;\Z)$ commutes with $\partial$.
\end{proof}

\begin{theorem}\label{thm_Z_orb}
For every compact polyhedron $Y\subseteq\R^K$ and every smooth manifold $Y$ without boundary, $p_*:H_*(Y;\Z)\to H_*(\overline{GC}{}^{G}_{\mathrm{orb}}(Y))$ is an isomorphism in all degrees.
\end{theorem}

\begin{corollary}\label{cor_Z_orb}
\begin{enumerate}[(a)]
\item Theorem~\ref{thm_Z_orb} implies Theorem~\ref{thm_orb}, with a proof that uses neither de Rham's theorem nor Stokes' theorem on orbifolds with corners.
\item The maps $\GCG(Y;\Z)\to\overline{GC}{}^{G}_{\mathrm{orb}}(Y)$ induced by $C^m_*(Y;\Z)\subseteq C^{\mathrm{orb}}_*(Y;\Z)$ and $\cG\subseteq\cG^{\mathrm{orb}}$ are isomorphisms on homology.
\end{enumerate}
\end{corollary}

\begin{proof}
(a) By Lemma~\ref{lem_Z_orb_complex}(a), right exactness and flatness of $\Q$, $\overline{GC}{}^{G}_{\mathrm{orb}}(Y)\otimes\Q=GC^G_{\mathrm{orb}}(Y)$ and homology commutes with $\otimes\Q$; the proof of Theorem~\ref{thm_Z_orb} below uses neither theorem, see the list of places where (A5) enters. (b) The maps commute with $p$, and $p_*$ is an isomorphism on both sides by Theorems~\ref{thm_Z} and \ref{thm_Z_orb}.
\end{proof}

\emph{Transversality.} For orbifold sources, regular values and $c$-transversality are read in charts, as in the proof of Proposition~\ref{prop_abstract_orb}. Likewise $[Q,f]$, $Q\in\cS^{\mathrm{orb}}$, is \emph{$K$-transverse} if for every chart $(\tilde U,H)$ the lift $\tilde f$ is transverse, on every open corner stratum of $\tilde U$, to every open simplex of $K'$. This does not depend on the charts, since embeddings of charts map corner strata locally diffeomorphically onto corner strata, and by the chart form of (J2) the boundary of a $K$-transverse class is $K$-transverse, so the $K$-transverse classes span a direct summand of $C^{\mathrm{orb}}_*(X;\Z)$ that is stable under $\partial$. $K$-transverse gluing data and $\cR^\pitchfork$ are defined as in Section~\ref{sec_int}. Since $K$-transversality does not depend on the charts, it may be checked in linear charts centred at the points; at depth $0$ these exist by averaging a metric (Bochner linearization), and at depth one by \cite[Doubling lemma]{Yu1}(1).

\begin{lemma}[rank]\label{lem_rank}
Let $Q\in\cS^{\mathrm{orb}}$, $\dim Q=k$, let $(\tilde U,H)$ be a linear chart centred at a point $y$ of depth $0$ or $1$ with $H\ne1$, $\tilde f$ an $H$-invariant smooth map, and $S$ the open corner stratum through $y$. Then $\operatorname{rank}d(\tilde f|_S)_y\le\dim S-2$. In particular, if $\tilde f$ is transverse on $S$ to a submanifold $E$ of codimension $\ge\dim S-1$, then $y\notin\tilde f^{-1}(E)$.
\end{lemma}

\begin{proof}
For $g\in H$, $\tilde f\circ g=\tilde f$ gives $d\tilde f_y\circ g=d\tilde f_y$ on $T_yS$, which $g$ preserves, so $d(\tilde f|_S)_y$ vanishes on $(g-1)T_yS$ and its rank is at most $\dim\operatorname{Fix}_{T_yS}(H)$. If $y$ has depth $0$, $H$ acts effectively and preserving orientation, and a non-trivial orientation preserving orthogonal map has a fixed space of codimension $\ge2$. If $y$ has depth one, by \cite[Doubling lemma]{Yu1}(1) $g(u,l)=(u,A_gl)$ with $\det A_g=\det g=1$, and $H$ acts effectively on $T_yS=\{0\}\times\R^{k-1}$, since $g$ fixes $u$, so again the fixed space has codimension $\ge2$ in $T_yS$. (At depth $2$ the bound is only $\dim S-1$, e.g.\ for $g(u_1,u_2,l_1,l_2)=(u_2,u_1,l_1,-l_2)$; from depth $3$ on, $H$ may act trivially on $T_yS$, e.g.\ by $g(u_1,u_2,u_3,l)=(u_2,u_3,u_1,l)$. Points of depth $\ge2$ are excluded by dimension below.) Transversality at $y$ to $E$ needs rank $\ge\operatorname{codim}E$.
\end{proof}

\begin{proposition}\label{prop_Z_orb_inj}
For every smooth manifold $X$ without boundary and every compact polyhedron $Y\subseteq\R^K$, $p_*:H_k(\,\cdot\,;\Z)\to H_k(\overline{GC}{}^{G}_{\mathrm{orb}}(\,\cdot\,))$ is injective.
\end{proposition}

\begin{proof}
Define $I$ on the $K$-transverse classes of $C^{\mathrm{orb}}_*(X;\Z)$ as in Section~\ref{sec_int}. For $\dim Q=\dim\tau=k$, $f^{-1}(D(\tau))$ lies in the depth-$0$ strata of charts, since strata of depth $\ge1$ have dimension $<k$, and by Lemma~\ref{lem_rank}, applied in charts centred at the points, it avoids points with non-trivial isotropy. So it is a finite set of interior points at which $|Q|$ is a manifold, and the count $\langle f,D(\tau)\rangle\in\Z$ is defined. We check Lemma~\ref{lem_dual} for orbifold sources.

(a) As before; $I$ takes values in the torsion-free group $C^{\mathrm{simp}}_*(K;\Z)$, so it vanishes on $R_{\mathrm{or}}$.

(b) For $\dim\rho=k-1$, by Lemma~\ref{lem_rank} applied on open strata and on strata of depth one, $f^{-1}(D(\rho))$ avoids the points with non-trivial isotropy, and it avoids strata of depth $\ge2$ by dimension. So it lies in the part of $|Q|$ that is a manifold with boundary near it, with ends at points of depth one, each with a single local boundary component, and the argument of Lemma~\ref{lem_dual}(b), which is local, applies once its ends are matched with the points counted in $\langle f\circ i_Q,D(\rho)\rangle$. The points of $(f\circ i_Q)^{-1}(D(\rho))$ are, by the first paragraph applied to the $K$-transverse class $\partial[Q,f]$, interior points of $\partial Q$ with trivial isotropy; by the chart form of (J2) a depth-one point of $Q$ has the same local group as its preimage in $\partial Q$, on which that group acts effectively since $g(u,l)=(u,A_gl)$. So these points are exactly the ends of the arcs.

(c) Let $(P;(Q_i,\iota_i))\in\fG^{\mathrm{orb}}$ be $K$-transverse with map $f$. Each $f\iota_i$ is $K$-transverse on the boundary strata of $Q_i$, which have dimension $<k$, so $f\iota_i(B_{Q_i})\cap D(\tau)=\emptyset$. Let $x\in f^{-1}(D(\tau))$; it is an interior point of $P$ with trivial isotropy. By (G2) it lies in $\iota_i(\inter Q_i)$ for exactly one $i$, $x=\iota_i(q)$, and by Lemma~\ref{lem_orb_interior}(b) $q$ has trivial isotropy and $\iota_i$ is an orientation preserving local diffeomorphism at $q$. So $I\rel{\fd}{f}=0$, and $I(\partial\rel{\fd}{f})=0$ by (b).

(d) Unchanged, since $p(c)$ has manifold sources.

Now let $c$ be a simplicial cycle with $p(c)=0$ in $H_k(\overline{GC}{}^{G}_{\mathrm{orb}}(X))$. Lifting to $C^{\mathrm{orb}}_*(X;\Z)$, $p(c)=\partial w+r+\theta$ with $r\in G^{\mathrm{orb}}(X;\Z)$, written with finitely many data, and $\theta\in R_{\mathrm{or}}$. The map $(\mathbf F_\beta)_*$ is a homomorphism, so it preserves $R_{\mathrm{or}}$, and it preserves $G^{\mathrm{orb}}$ by (A1). For each of the finitely many charts and corner strata $S$ involved, $(s,\beta)\mapsto\mathbf F(\tilde f(s),\beta)$ is a submersion for small $\beta$, so by parametric transversality all classes and data involved are $K$-transverse after moving by $\mathbf F_\beta$, for almost all small $\beta$. Then
\[
c=I(p(c)_\beta)=\partial I(w_\beta)+I(r_\beta)+I(\theta_\beta)=\partial I(w_\beta).
\]
The polyhedral case follows as in Proposition~\ref{prop_Z_inj}.
\end{proof}

\begin{proof}[Proof of Theorem~\ref{thm_Z_orb}]
Injectivity is Proposition~\ref{prop_Z_orb_inj}. For surjectivity we run the proofs of Propositions~\ref{prop_abs_poly} and \ref{prop_abs_mfd} in $\bar C$, for the family $\fG^{\mathrm{orb}}$, which satisfies (A1)--(A4) by Proposition~\ref{prop_orb_A15}, with the items of the proof of Proposition~\ref{prop_abstract_orb} other than \emph{Integration}. (A5) enters the proof of Theorem~\ref{thm_abstract} only through Propositions~\ref{prop_stokes} and \ref{prop_poly_inj}, which are replaced by Proposition~\ref{prop_Z_orb_inj}, Step 5 of Lemma~\ref{lem_MV}, and degree $0$ of Lemma~\ref{lem_star}(b); the last two are treated below. Homomorphisms of $C^{\mathrm{orb}}_*(\cdot\,;\Z)$, such as $\psi_*$, the prism operator of Lemma~\ref{lem_prism} and $(\mathbf F_\beta)_*$, preserve torsion and descend to $\bar C$; the identities between them, such as the homotopy formula, hold generator by generator. The remaining points are the following.
\begin{itemize}
\item \emph{Cut operators.} $c$-transversality is a property of isomorphism classes, checked componentwise, so the span $C^\pitchfork_c(Y;\Z)$ of the $c$-transverse classes is a direct summand of $C^{\mathrm{orb}}_*(Y;\Z)$, spanned by a subset of the summands of Lemma~\ref{lem_Z_orb_complex}(a). The formulas for $\cut_\pm$ and $\Lambda$ respect isomorphisms and disjoint unions and satisfy $(-Q)_\pm=-(Q_\pm)$, $L_{-Q}=-L_Q$, so they are well defined on $C^\pitchfork_c(Y;\Z)$; being homomorphisms, they descend to its image $\bar C^\pitchfork_c(Y)$ in $\bar C(Y)$, a subcomplex spanned by a subset of the basis. The first two identities \eqref{eq_cut} hold generator by generator by the chart form of Lemma~\ref{lem_pieces}(b), and the third follows from $\partial^2=0$ in $\bar C$. So Lemmas~\ref{lem_cutops}, \ref{lem_cutrel} and Steps 1--4, 6, 7 of Lemma~\ref{lem_MV} hold in $\bar C$.
\item \emph{Step 5 of Lemma~\ref{lem_MV}.} Injectivity on the halves, which are compact polyhedra, is Proposition~\ref{prop_Z_orb_inj}. This replaces (A5) in Lemma~\ref{lem_MV}.
\item \emph{Degree $0$.} Sources of dimension $0$ are finite sets of points, and those of dimension $1$ are compact $1$-manifolds, since an orientation preserving linear action on $\R$, and any linear action on $[0,\infty)$, is trivial, and local groups act effectively. So the signed count of points, which kills torsion and hence is defined on $\bar C_0$, $\varepsilon:\bar C_0\to\Z$ vanishes on $\partial\bar C_1$ and, by (G1$'$) and (G2), on gluing relations in degree $0$, which gives Lemma~\ref{lem_star}(b) in degree $0$; $p$ is onto in degree $0$ at chain level.
\item \emph{(A6).} For $Q\in\cS^{\mathrm{orb}}$ of dimension $k\ge1$, Steps 2 and 3 of the proof of Proposition~\ref{prop_orb_point} are bisections and fan decompositions, all with coefficients $1$, so $[Q,c]-\sum_{a,\sigma}[Q_{a,\sigma},c]\in\operatorname{span}_\Z\cG^{\mathrm{orb}}(\mathrm{pt})$. Each $Q_{a,\sigma}$ is a manifold with embedded faces, and by the proof of Corollary~\ref{cor_point}, $[Q_{a,\sigma},c]\equiv\sum_i\varepsilon_i[\Delta^k,c]$ modulo $\operatorname{span}_\Z\cG(\mathrm{pt})$. The transposition of two vertices is an orientation reversing automorphism of $(\Delta^k,c)$, so $[\Delta^k,c]\in R_{\mathrm{or}}$. Hence $\overline{GC}{}^{G}_{\mathrm{orb},k}(\mathrm{pt})=0$ for $k\ge1$ at chain level.
\item \emph{Effectiveness.} All sources are effective, and all constructions stay in $\cS^{\mathrm{orb}}$ (Section~\ref{sec_orb}), so the weights of non-effective sources in \cite{Yu1} do not occur.
\end{itemize}
Proposition~\ref{prop_amb} is integral. This proves the theorem.
\end{proof}

\begin{remark}\label{rem_Z_orb}
For manifold sources the division by $R_{\mathrm{or}}$ is not needed (Theorem~\ref{thm_Z}), and the same proof shows that dividing $\Cm_*(Y;\Z)$ by its torsion does not change the homology either.
\end{remark}

\section{Orbifold targets}\label{sec_orbtarget}
In Sections~\ref{sec_orb} and \ref{sec_int} the targets are polyhedra and manifolds. Here the target is an effective orbifold, possibly with corners and not compact, and maps are liftable maps or orbifold morphisms, as in \cite[Section on liftable maps and orbifold morphisms]{Yu1}, where the comparison is proved with the relation $D$. We show that the gluing relation alone suffices. The proof of Section~\ref{sec_proof} embeds manifolds in $\R^K$ and reduces to polyhedra, which is not available for orbifolds; we follow instead the pattern of \cite[proof of the main theorem]{Yu1} (good covers, Mayer--Vietoris for open sets, compact supports) and replace the three places where $D$ enters there.

\subsection{Setting}
Let $X$ be an effective orbifold with corners. Sources are in $\cS^{\mathrm{orb}}$ (Section~\ref{sec_orb}). For $\star\in\{\mathrm{lift},\mathrm{mor}\}$ let $C^\star_*(X)$ be the complex of classes $[Q,f]$, $Q\in\cS^{\mathrm{orb}}$, with $f:|Q|\to|X|$ liftable ($\star=\mathrm{lift}$), resp.\ of classes $[Q,\hat f]$ with $\hat f:Q\to X$ an orbifold morphism ($\star=\mathrm{mor}$), with the isomorphisms, the boundary and the basis of \cite[Definition of chains given by liftable maps and orbifold morphisms]{Yu1}; let $C^{\mathrm{sm},\star}_*(X)$ be the complex of liftable singular simplices, resp.\ of ($2$-isomorphism classes of) morphisms $\hat\sigma:\Delta^n\to X$, and $p$ the comparison map. For $O\subseteq X$ open, $C^\star_*(O)\to C^\star_*(X)$ and $C^{\mathrm{sm},\star}_*(O)\to C^{\mathrm{sm},\star}_*(X)$ are injective, onto the span of a subset of the basis: an open suborbifold is a full subgroupoid, so morphisms with image in $O$, their $2$-isomorphisms, the isomorphisms of classes and orientation reversing automorphisms are the same in $O$ and in $X$. In the notation of \cite{Yu1}, $\star=\mathrm{lift}$ is $\star=\mathrm{orb}$, our $C^\star_*$ is $C^{g,\star}_*$ and our $C^{\mathrm{sm},\star}_*$ is $C^\star_*$; the complex $C^{\mathrm{orb}}_*(Y)$ of Section~\ref{sec_orb} is different. The map $p_*$ of Theorem~\ref{thm_orbtarget} is $p_*$ on $H_*(C^{\mathrm{sm},\star}_*(X))$ composed with the inverse of the isomorphism $H_*(C^{\mathrm{sm},\star}_*(X))\cong H_*(|X|;\Q)$ \cite[Proposition on chains given by liftable maps and orbifold morphisms, part (A)]{Yu1}.

\begin{lemma}[morphisms of embeddings]\label{lem_emb_mor}
Let $P,Q\in\cS^{\mathrm{orb}}$ have the same dimension, and let $\iota:|Q|\to|P|$ satisfy the injectivity and lifting conditions of \textup{(G1$'$)}. Then there is an orbifold morphism $\hat\iota:Q\to P$ with $|\hat\iota|=\iota$, and every morphism $\hat\varphi$ with $|\hat\varphi|=\iota$ is $2$-isomorphic to $\hat\iota$ by a unique $2$-isomorphism. If $j:|Q'|\to|Q|$ also satisfies these conditions, then $\widehat{\iota\circ j}\cong\hat\iota\circ\hat j$. The existence of a morphism over $\iota\circ j_{L_Q}:|L_Q|\to|P|$, and its uniqueness up to a unique $2$-isomorphism, also hold, where $L_Q=\{h\circ\iota=c\}\subseteq Q$ and $L_P=\{h=c\}\subseteq P$ are level sets of a smooth $h:|P|\to\R$, $(h,c)$ is regular in the chart sense for $h$ on $P$ and $h\circ\iota$ on $Q$, and $j_{L_Q}$, $j_{L_P}$ are the inclusions: $\hat\iota\circ\hat j_{L_Q}\cong\hat j_{L_P}\circ\hat\iota_0$, where $\iota_0:L_Q\to L_P$ is the restriction, which satisfies \textup{(G1$'$)} by the argument for (A3) in Proposition~\ref{prop_orb_A15}, and $\hat\iota_0$ is given by the first part.
\end{lemma}

\begin{proof}
\emph{Local uniqueness.} Let $(\tilde U,H)$ be a connected chart of $Q$, $(\tilde V,G)$ a connected chart of $P$, and $\tilde\iota_1,\tilde\iota_2:\tilde U\to\tilde V$ two smooth lifts of $\iota$, with $\tilde\iota_1$ an embedding. The closed sets $F_\eta=\{\tilde\iota_2=\eta\tilde\iota_1\}$, $\eta\in G$, cover $\tilde U$, since points of $\tilde V$ with the same image in $|P|$ lie in one $G$-orbit; so the union of their interiors $E_\eta$ is dense and the closures $\bar E_\eta$ cover $\tilde U$. Let $w\in\bar E_\eta\cap\bar E_{\eta'}$, $y=\tilde\iota_1(w)$ and $g=\eta^{-1}\eta'$. Passing to the limit, $gy=y$ and $dg_y\circ d\tilde\iota_1(w)=d\tilde\iota_1(w)$, where $d$ denotes differentials of smooth extensions on all of $\R^n$; since $\dim Q=\dim P$ and $\tilde\iota_1$ is an embedding, $d\tilde\iota_1(w)$ is invertible, so $dg_y=\id$. For an element $g$ of finite order $m$ with $gy=y$ and $dg_y=\id$, the average $\phi=\frac1m\sum_kg^k$ of a smooth extension has $d\phi_y=\id$, so it is injective near $y$, and $\phi\circ g=\phi$; hence $g$ is the identity near $y$, also at points of positive depth. Applied at the points of the closure of $A=\inter\operatorname{Fix}(g)\ne\emptyset$, where $dg=\id$ by continuity, this shows that $A$ is closed, so $A=\tilde V$, and $g=1$ by effectiveness. Hence the sets $\bar E_\eta$ are pairwise disjoint, and since $\tilde U$ is connected, $\tilde\iota_2=\eta\tilde\iota_1$ for a single $\eta$. It is unique: if $\eta\tilde\iota_1=\eta'\tilde\iota_1$, then $\eta^{-1}\eta'$ fixes $\tilde\iota_1(\tilde U)$, which contains a non-empty open set, hence is trivial by the same argument. For $\gamma\in H$, $\tilde\iota\circ\gamma$ is another lift, so $\tilde\iota\circ\gamma=\lambda(\gamma)\circ\tilde\iota$ for a unique $\lambda(\gamma)$; $\lambda$ is a homomorphism by uniqueness, and injective since $H$ acts effectively.

\emph{Gluing.} Choose an atlas of $Q$ by connected charts $\tilde U_i$ with embedding lifts $\tilde\iota_i$ into connected charts $\tilde V_{k(i)}$ of $P$. Let $\varphi$ be a transition germ of $Q$ from $u\in\tilde U_i$ to $u'\in\tilde U_j$. (i) There is a transition germ $\psi$ of $P$ with $\psi\circ\tilde\iota_i=\tilde\iota_j\circ\varphi$: take any transition germ $\psi_0$ of $P$ from $\tilde\iota_i(u)$ into $\tilde V_{k(j)}$; then $\psi_0\circ\tilde\iota_i$ and $\tilde\iota_j\circ\varphi$ are embedding lifts of $\iota$ near $u$ into the same chart, so they differ by a unique element $\eta$ of its local group, and $\psi=\eta\psi_0$. (ii) $\psi$ is unique: the quotient of two such germs is a self-transition germ of $P$ at $y=\tilde\iota_i(u)$, hence the germ of an element $g$ of the local group \cite[Lemma on orbifold morphisms, (d)]{Yu1}, which fixes $\tilde\iota_i(\tilde U_i)$ near $y$. This set need not be a neighbourhood of $y$ (for instance at the vertex of a cone piece in Proposition~\ref{prop_orb_point}), but it contains an open set accumulating at $y$, so $g=1$ as above. By (i) and (ii), $u\mapsto\tilde\iota_i(u)$, $\varphi\mapsto\psi$ is a functor between the groupoids of germs of the two atlases, smooth since the arrow map is locally constant, and it defines $\hat\iota$. A morphism $\hat\varphi$ with $|\hat\varphi|=\iota$ is represented near each point by a pair $(\tilde\varphi,\mu)$ \cite[Lemma on orbifold morphisms, (a)]{Yu1}, and $\tilde\varphi=\eta\tilde\iota$, $\mu=\eta\lambda\eta^{-1}$ by local uniqueness; the elements $\eta$, compatible on overlaps by local uniqueness, form a $2$-isomorphism, unique by local uniqueness. For $\iota\circ j$, after composing with embeddings of charts as in \cite[proof of the Lemma on orbifold morphisms, (b)]{Yu1}, the lifts compose to embeddings, and $\hat\iota\circ\hat j$ is a morphism over $\iota\circ j$.

\emph{Level sets.} Both $\hat\iota\circ\hat j_{L_Q}$ and $\hat j_{L_P}\circ\hat\iota_0$ have lifts which are codimension one embeddings into the invariant hypersurface $\tilde L=\{\tilde h=c\}$ of a chart of $P$, over the same map. In the argument above, $g$ now fixes $y$ and $dg_y=\id$ on $T_y\tilde L$; since $\tilde h$ is invariant, $dg_y$ preserves $d\tilde h_y\ne0$, so $dg_y\nu-\nu\in T_y\tilde L$ for a normal vector $\nu$, and $dg_y$ is unipotent of finite order, hence the identity, as in the proof of the item \emph{Pieces} of Proposition~\ref{prop_abstract_orb}. The rest of the argument applies, with open subsets of $\tilde L$ in place of open subsets of $\tilde V$ and the same covector argument for uniqueness.
\end{proof}

So the morphisms of the pieces of an orbifold gluing datum are determined by the datum. The boundary inclusions and the inclusions of the halves and of the level set of a cut, used in the boundary formulas of Lemma~\ref{lem_cutops} and in the homotopy formula, are strict in the induced atlases (inclusions of hyperplanes and half-spaces, with the identity on groups) \cite[Lemma on orbifold morphisms, (c)]{Yu1}, so their compositions with each other and with the restrictions of morphisms agree on the atlases, and the identities $\hat j\circ\hat i\cong\hat i\circ\hat j$ needed there hold. For $\star=\mathrm{mor}$ we put
\[
\rel{\fd}{\hat f}:=[P,\hat f]-\sum_i[Q_i,\hat f\circ\hat\iota_i],
\]
and for $\star=\mathrm{lift}$, $\rel{\fd}{f}=[P,f]-\sum_i[Q_i,f\circ\iota_i]$ as before. Let $\cG^\star(X)$ be the set of these relations for $\fd\in\fG^{\mathrm{orb}}$ and all $f$, resp.\ $\hat f$, $G^\star(X)=\operatorname{span}\cG^\star(X)+\partial\operatorname{span}\cG^\star(X)$, and $GC^G_\star(X):=C^\star_*(X)/G^\star(X)$. By Lemma~\ref{lem_emb_mor}, the morphisms of the pieces of the cut data in (A3) and of the prism data in (A4) are the restrictions, resp.\ products, of the given ones, up to unique $2$-isomorphism; and those of the cone pieces of Proposition~\ref{prop_orb_point} are induced by the inclusions $(\sigma\cap\tilde V_a,1)\to G_a\ltimes\tilde V_a$.

\begin{theorem}[orbifold targets]\label{thm_orbtarget}
For every effective orbifold with corners $X$ and $\star\in\{\mathrm{lift},\mathrm{mor}\}$, $p_*:H_*(|X|;\Q)\to H_*(GC^G_\star(X))$ is an isomorphism in all degrees.
\end{theorem}

\subsection{Integration}
For an $n$-form $\omega$ on an open set $O\subseteq X$ and $[Q,f]\in C^\star_n(O)$ put $\int_{[Q,f]}\omega:=\int_Qf^*\omega$, with the weights $1/|G|$ of the orbifold integral. Here $f^*\omega$ is defined for liftable maps and is independent of the lift \cite[Remark on injectivity via Stokes]{Yu1}; for morphisms it depends only on the $2$-isomorphism class. It is invariant under isomorphisms and changes sign with the orientation, so it is well defined, and Stokes' theorem on orbifolds with corners holds \cite[Remark on injectivity via Stokes]{Yu1}. Since $(f\circ\iota_i)^*\omega=\iota_i^*f^*\omega$, the proof of (A5) in Proposition~\ref{prop_orb_A15}, applied to the $n$-form $f^*\omega$ on $P$, shows that integration vanishes on $\cG^\star(O)$, hence by Stokes on $G^\star(O)$.

Fix a good cover $\cU$ of an effective orbifold $X$ without boundary as constructed in the proof of \cite[Proposition on good covers of orbifolds]{Yu1}: every non-empty finite intersection of its members is a finite disjoint union of sets $W/H$, where $(W,H)$ is a chart, $W$ is strongly convex for an invariant metric, $H$ fixes a point $c\in W$, and $W$ lies in a ball around $c$ of radius smaller than the injectivity radius at $c$, so that $\exp_c^{-1}$ is smooth on $W$. Let $\cU^\cap$ be the set of these intersections.

\begin{lemma}[de Rham pairing]\label{lem_dR_pair}
Let $O$ be a finite union of members of $\cU^\cap$. Integration $\Omega^*(O)\to\operatorname{Hom}(C^{\mathrm{sm},\star}_*(O),\R)$, with the coboundary $\varphi\mapsto\varphi\circ\partial$, is a quasi-isomorphism. Hence a cycle $s\in C^{\mathrm{sm},\star}_n(O)$ with $\int_s\omega=0$ for all closed $\omega$ is a boundary in $C^{\mathrm{sm},\star}_*(O)$.
\end{lemma}

\begin{proof}
On a simplex, $\sigma^*\omega$ is defined locally by a lift, independently of it, and Stokes' theorem on $\Delta^n$ shows that integration is a cochain map, natural under inclusions of open sets. For $O=W/H$, the geodesic contraction $K(x,s)=\exp_c((1-s)\exp_c^{-1}x)$ is smooth and $H$-equivariant, so the homotopy operator $\omega\mapsto\int_0^1\iota_{\partial_s}K^*\omega\,ds$, which satisfies $dh+hd=K_1^*-K_0^*$, preserves invariant forms, and $H^*_{dR}(W/H)=\R$ in degree $0$. By \cite[proof of the Proposition on chains given by liftable maps and orbifold morphisms, items (3) and (6)]{Yu1}, $H_*(C^{\mathrm{sm},\star}_*(W/H))=\Q$ in degree $0$, so the dual complex has cohomology $\R$ in degree $0$, and integration maps the constant $1$ to the cocycle with value $1$ on points; this is an isomorphism, also for finite disjoint unions. For $O=A\cup B$ the sequence $0\to\Omega(O)\to\Omega(A)\oplus\Omega(B)\to\Omega(A\cap B)\to0$ is exact by invariant partitions of unity ($O$ is paracompact, being metrizable). The sequence $0\to C^{\mathrm{sm},\star}(A\cap B)\to C^{\mathrm{sm},\star}(A)\oplus C^{\mathrm{sm},\star}(B)\to C^{\mathrm{sm},\star,\{A,B\}}(O)\to0$ of free modules is exact by the injectivity statement of the setting, and stays exact after dualising; the inclusion of $\{A,B\}$-small simplices is a quasi-isomorphism by \cite[proof of the Proposition on chains given by liftable maps and orbifold morphisms, item (7)]{Yu1}. Integration maps the first sequence to the dual of the second, so by the five lemma, if the statement holds for $A$, $B$ and $A\cap B$, it holds for $O$. For $O=V_1\cup\dots\cup V_m$ with $V_i\in\cU^\cap$ put $A=V_1\cup\dots\cup V_{m-1}$ and $B=V_m$; then $A\cap B=\bigcup_{j<m}(V_j\cap V_m)$ with $V_j\cap V_m\in\cU^\cap\cup\{\emptyset\}$, and induction on $m$ proves the first statement. The second follows from the universal coefficient theorem over $\R$ and the injectivity of $H_*(C;\Q)\to H_*(C;\Q)\otimes\R$.
\end{proof}

\subsection{Proof of Theorem~\ref{thm_orbtarget}}
\emph{Corners.} By \cite[Corollary on orbifolds with corners, and proof of the Proposition on chains given by liftable maps and orbifold morphisms, item (9)]{Yu1}, the flow of an invariant inward vector field gives $r:X\to X^\circ$ and homotopies $\id_X\simeq i\circ r$, $\id_{X^\circ}\simeq r\circ i$ which are induced by strict morphisms of groupoids, in particular liftable. Composition with such a homotopy is a prism operator as in Lemma~\ref{lem_prism} (homotopy formula as in \cite[same proof, item (3)]{Yu1}), which maps $G^\star$ into $G^\star$ by (A4) and Lemma~\ref{lem_emb_mor}. So $i_*$ is an isomorphism on $H_*(GC^G_\star)$, and also on singular homology, and we may assume that $X$ has no boundary.

\emph{Injectivity.} Let $U\subseteq X$ be open and $s\in C^{\mathrm{sm},\star}_n(U)$ a cycle with $p(s)=\partial w+e$, $e\in G^\star(U)$. Choose a good cover of $U$ as above. Fix an expression of $e$ as an element of $\operatorname{span}\cG^\star(U)+\partial\operatorname{span}\cG^\star(U)$. The images of the finitely many sources and pieces involved in $s$, $w$ and this expression are compact, so they lie in a finite union $O$ of members of $\cU^\cap$; the relations of the expression have maps into $O$, so they lie in $\cG^\star(O)$ and $e\in G^\star(O)$. The injectivity statement of the setting is used only to regard $s$, $w$, $e$ as chains of $O$. For closed $\omega$ on $O$, $\int_s\omega=\int_wd\omega+\int_e\omega=0$, so $s$ is a boundary in $C^{\mathrm{sm},\star}(O)$ by Lemma~\ref{lem_dR_pair}, hence in $C^{\mathrm{sm},\star}(U)$, and $H_*(C^{\mathrm{sm},\star}(U))\cong H_*(|U|;\Q)$ by \cite[Proposition on chains given by liftable maps and orbifold morphisms, part (A)]{Yu1}.

\emph{Contractible pieces.} Let $W/H$ be as above and $\hat K:[W/H]\times I\to[W/H]$ the strict morphism induced by the geodesic contraction to $c$, with $\hat K_0=\id$. Let $z$ be a cycle of $GC^G_\star(W/H)$ of degree $k\ge1$. By the prism operator, $z\equiv(\hat K_1)_*z$ modulo boundaries and $G^\star$. Let $[Q,\hat f]$ be a term of $z$. Near each point $q$ of $Q$ choose a linear chart $(\tilde U,G)$ which is convex and $G$-invariant, on which $\hat f$ is represented by a pair $(\tilde f,\lambda)$ into $(W,H)$: by \cite[Lemma on orbifold morphisms, (a)]{Yu1} into a small chart $(B,H_{\tilde y})$ centred at a lift $\tilde y$ of $|\hat f|(q)$, composed with the embedding $B\hookrightarrow W$, so $\lambda:G\to H_{\tilde y}\subseteq H$; then $\hat K_1\circ\hat f$ is represented by $(\mathrm{const}_c,\lambda)$. Run Steps 2 and 3 of the proof of Proposition~\ref{prop_orb_point} with finitely many such charts. The relations obtained do not depend on the map, so $[Q,\hat K_1\hat f]\equiv\sum_{a,\sigma}[Q_{a,\sigma},\hat K_1\hat f\hat\iota_{a,\sigma}]$ modulo $\operatorname{span}\cG^\star(W/H)$. On a cone piece the local group is trivial, so $\hat K_1\hat f\hat\iota_{a,\sigma}$ is represented by $(\mathrm{const}_c,1)$, i.e.\ it is $\hat c\circ\pi$, where $\pi:Q_{a,\sigma}\to\mathrm{pt}$ and $\hat c:\mathrm{pt}\to[W/H]$ is the strict morphism at the fixed point $c$. By Corollary~\ref{cor_point}, $[Q_{a,\sigma},\pi]\in\operatorname{span}\cG(\mathrm{pt})$, and $\hat c_*$ maps it into $\operatorname{span}\cG^\star(W/H)$. Hence $(\hat K_1)_*z\in\operatorname{span}\cG^\star(W/H)$ and $H_k(GC^G_\star(W/H))=0$ for $k\ge1$. For $\star=\mathrm{lift}$, $K_1\circ f$ is constant, and $[Q,K_1f]\in\operatorname{span}\cG^{\mathrm{lift}}(W/H)$ by Proposition~\ref{prop_orb_point}, pushed forward along the constant map. In degree $0$, sources are finite sets of points and those of dimension $1$ are compact $1$-manifolds, since a finite group acting effectively and preserving orientation on an interval or a circle chart is trivial, morphisms from a point to $X$ are determined up to $2$-isomorphism by their image, so the signed count of points vanishes on $\partial C^\star_1$ and on gluing relations of degree $0$, and $p_*$ is an isomorphism for $W/H$. Finite disjoint unions: a relation splits over the components of $P$, since the partner faces of (G3$'$) have the same image and lie over the same component.

\emph{Mayer--Vietoris.} Let $O=U\cup V$ with $U$, $V$ open, and assume that $p_*$ is onto in degree $k$ for $U$ and $V$ and in degree $k-1$ for $U\cap V$. Choose a smooth $\phi:O\to[0,1]$ with $\phi=0$ on $O\setminus V$ and $\phi=1$ on $O\setminus U$, and, by Sard's theorem, $t\in(0,1)$ such that the finitely many generators involved in a cycle $x$ of $GC^G_\star(O)$ and in a chosen expression of $\partial x$ are transverse to $t$ for $\phi\circ f$ in the chart sense of Proposition~\ref{prop_abstract_orb}. The cut operators $\cut_\pm$ and $\Lambda$ along $\phi$ satisfy Lemmas~\ref{lem_cutops} and \ref{lem_cutrel}, read in charts as in Proposition~\ref{prop_abstract_orb}, with the morphisms of the pieces given by Lemma~\ref{lem_emb_mor}; they map into $C^\star(\{\phi\le t\})\subseteq C^\star(U)$, $C^\star(\{\phi\ge t\})\subseteq C^\star(V)$ and $C^\star(\{\phi=t\})\subseteq C^\star(U\cap V)$. Steps 1--7 of the proof of Lemma~\ref{lem_MV} now apply with $U$, $V$, $U\cap V$ in place of $Y_-$, $Y_+$, $Y_0$: Step 4 uses surjectivity for $U\cap V$ (if $U\cap V=\emptyset$, then $z=0$), Step 5 uses injectivity for $U$ and $V$, and Step 7 uses surjectivity for $U$ and $V$, followed by the push forward into $O$. Hence $p_*$ is onto in degree $k$ for $O$.

\emph{Conclusion.} By induction on the degree and on the number of members, $p_*$ is onto for finite unions of members of $\cU^\cap$, since $(V_1\cup\dots\cup V_{m-1})\cap V_m$ is a union of $m-1$ members. A cycle $x$ of $GC^G_\star(X)$, with an expression of $\partial x$, involves finitely many generators, whose images lie in such a finite union $O$; there $x\equiv p(s)$, and pushing into $X$ gives surjectivity. Injectivity was shown above.
\qed

\begin{proposition}[finite coverings]\label{prop_cover}
Let $P,P'\in\cS^{\mathrm{orb}}$ and let $\pi:P'\to P$ be an orientation preserving orbifold covering of degree $d$: every point of $|P|$ has a linear chart $(\tilde U,G)$ with image $U_0$ such that $\pi^{-1}(U_0)$ is isomorphic over $U_0$ to $\bigsqcup_k G_k\ltimes\tilde U$, with $G_k\subseteq G$, $\pi$ induced by $\id_{\tilde U}$, and $\sum_k[G:G_k]=d$ (the terms with $G_k=G$ are trivial sheets). Then for every orbifold morphism $\hat f:P\to X$ to an effective orbifold, $[P',\hat f\circ\hat\pi]-d\,[P,\hat f]\in\operatorname{span}\cG^{\mathrm{mor}}(X)$; the same holds for liftable maps, for smooth maps to polyhedra and manifolds, and, for the latter, with integer coefficients in the complexes of Section~\ref{sec_int}.
\end{proposition}

\begin{proof}
Run Step 2 of the proof of Proposition~\ref{prop_orb_point} on $P$ with finitely many evenly covered linear charts, shrunk by \cite[Lemma on orbifold morphisms, (a)]{Yu1} so that $\hat f$ is represented on each by a single pair (a $G$-invariant neighbourhood of the centre of an evenly covered chart is evenly covered), choosing the values $t_j$ regular also for the lifts to the charts of $P'$, which are the same functions on the same $\tilde U_j$; and on $P'$ with the functions $\chi_j\circ\pi$. These are relations of $\cG^\star$, and $\pi^{-1}(Q_a)=\bigsqcup_kG_{a,k}\ltimes\tilde V_a$. In Step 3 use on each $G_{a,k}\ltimes\tilde V_a$ the $G_a$-invariant fan $\Sigma_a$, which is $G_{a,k}$-invariant, and the properties of Step 1 hold for the subgroup $G_{a,k}$. As representatives of the $G_{a,k}$-orbits of $n$-cones take the cones $g\sigma$, where $\sigma$ runs over the chosen representatives of the $G_a$-orbits and $g$ over representatives of $G_{a,k}\backslash G_a$; so each $\sigma$ gives $[G_a:G_{a,k}]$ pieces $g\sigma\cap\tilde V_a$ in this sheet, and $d$ pieces in all. The linear map $g^{-1}$ is an orientation preserving isomorphism $g\sigma\cap\tilde V_a\to\sigma\cap\tilde V_a$. If $\hat f$ is represented on $\tilde U_a$ by $(\tilde f,\lambda)$, then on the cone pieces, whose local groups are trivial, $\hat f\circ\hat\pi\circ\hat\iota'_{g\sigma}$ and $\hat f\circ\hat\iota_\sigma\circ\widehat{g^{-1}}$ are represented by $(\tilde f|_{g\sigma},1)$ and $(\tilde f\circ g^{-1}|_{g\sigma},1)=(\lambda(g)^{-1}\tilde f|_{g\sigma},1)$, which are $2$-isomorphic. Hence the $d$ pieces over $\sigma$ have the class $[Q_{a,\sigma},\hat f\hat\iota_{a,\sigma}]$, and $[P',\hat f\hat\pi]\equiv d\sum_{a,\sigma}[Q_{a,\sigma},\hat f\hat\iota_{a,\sigma}]\equiv d\,[P,\hat f]$. For liftable or smooth maps, $f\circ g=f$ on $|P|$. Only relations with coefficients $\pm1$ occur.
\end{proof}

\begin{remark}[orbifold targets]\label{rem_orbtarget}
(1) The relation $D$ enters \cite[proof of the main theorem]{Yu1} at three places: for the contractible pieces, where for orbifold morphisms the thin factorizations must be allowed through non-effective orbifolds \cite[Remark on what changes]{Yu1}; in the Mayer--Vietoris argument, through thin interfaces; and in the colimit, through the locality of $D$. For $D=0$ the last two would be trivial, but not for $G$, whose locality is not known (Remark~\ref{rem_locality}). Here they are replaced by the cone pieces of Proposition~\ref{prop_orb_point}, by cutting an explicit expression of the boundary as in Lemma~\ref{lem_MV}, and by compact supports; injectivity comes from integration instead of Mayer--Vietoris. For instance, a circle mapped to $\R^2/\mathbb Z_3$ through a non-trivial principal $\mathbb Z_3$-bundle and the singular point, which is not thin for effective factorizations, lies in $\operatorname{span}\cG^{\mathrm{mor}}$ at chain level: cut $S^1$ into arcs over which the bundle is trivial; on each arc the morphism is $\hat c\circ\pi$, and $[I,\hat c\circ\pi]=0$ in $C^{\mathrm{mor}}$, the flip of $I$ being an orientation reversing automorphism.

(2) Under $D$, Proposition~\ref{prop_cover} fails at chain level in $C/D$, for instance for $[S^1,z^2]-2[S^1,z]$, neither term being thin. We expect it to be the chain level identity needed for the usual connected form of the master equation for string vertices \cite{YuV}.

(3) We do not treat pullbacks along orbifold bundles over orbifold targets, which are defined on $C^{\mathrm{mor}}$. The proof of Proposition~\ref{prop_orb_ops}(b) does not extend to $\cG^{\mathrm{mor}}$: a piece may have a smaller local group than the glued source at a boundary point, and its pullback need not be injective there, so the pulled back datum need not be an orbifold gluing datum. (For the cone pieces of $D/\mu_3$ and $E=(D\times S^1)/\mu_3$, with $\mu_3$ rotating $S^1$, the pulled back relation still lies in $\operatorname{span}\cG^{\mathrm{mor}}$, by cutting $S^1$ into $3m$ arcs permuted by $\mu_3$.) We do not know whether $\pi^*$ maps $G^{\mathrm{mor}}$ into $G^{\mathrm{mor}}$ in general.
\end{remark}

\section{Relative homology, products, pullbacks and Morse complexes}\label{sec_relprod}
We record targets with corners for manifold sources, the relative version of the comparison theorems, the cross product, pullbacks along fibre bundles, and the chain map from the Morse complex. The relative part is formal; it applies, with the modifications stated there, to the complexes $GC^G_{\mathrm{orb}}$ of Section~\ref{sec_orb} and $GC^G_\star$ of Section~\ref{sec_orbtarget}.

\subsection{Targets with corners}
\begin{proposition}[boundary and corners]\label{prop_corners}
Let $X$ be a manifold with corners, let $\Cm(X)$ be spanned by the classes $[Q,f]$, $Q\in\cS$, with $f:Q\to X$ smooth (Convention~\ref{conv_smooth}), let $G(X)$ be defined by the gluing data as in Definition~\ref{def_rel}, and let $\Csm(X)$ be the smooth singular chains. Then $H_*(\Csm(X))\cong H_*(X;\Q)$, and $p_*:H_*(X;\Q)\to H_*(\GCG(X))$ is an isomorphism.
\end{proposition}

\begin{proof}
As in the step \emph{Corners} of the proof of Theorem~\ref{thm_orbtarget}, with trivial groups, the flow of an inward vector field gives $r:X\to\inter X$ and homotopies $H:\id_X\simeq j\circ r$ in $X$ and $\id_{\inter X}\simeq r\circ j$ in $\inter X$, where $j$ is the inclusion, weakly smooth on $X\times I$ \cite[Corollary on orbifolds with corners]{Yu1}. For smooth $f:Q\to X$ the map $H\circ(f\times\id)$ is smooth, so $\cH[Q,f]=(-1)^n[Q\times I,H\circ(f\times\id)]$ is defined and satisfies the homotopy formula \cite[Proposition on homotopy invariance]{Yu1}; the proof of Lemma~\ref{lem_prism}, stated for homotopies on open subsets of $\R^K$, applies verbatim. It maps $\cG$ into $\cG$ by Lemma~\ref{lem_products}, since gluing data do not involve the maps, and hence $G$ into $G$ as in the proof of Lemma~\ref{lem_prism}. The same holds for smooth singular chains. So $j$ induces isomorphisms on $H_*(\GCG)$ and on $H_*(\Csm)$, compatibly with $p$. The manifold $\inter X$ has no boundary; regarded as a closed submanifold of some $\R^K$, its intrinsic complexes are those of Section~\ref{sec_setting}, $\Csm(\inter X)$ computes $H_*(\inter X;\Q)\cong H_*(X;\Q)$, and Theorem~\ref{thm_A_G} applies. For orbifold sources and orbifold targets with corners see Theorem~\ref{thm_orbtarget}.
\end{proof}

\subsection{Relative homology}
For $W\subseteq Y\subseteq\R^K$ the inclusion induces a chain map $i:\GCG(W)\to\GCG(Y)$, since $G(W)\subseteq G(Y)$. Put
\[
\GCG(Y,W):=\operatorname{Cone}(i),\qquad\GCG_n(Y,W)=\GCG_{n-1}(W)\oplus\GCG_n(Y),
\]
with $\partial(a,b)=(-\partial a,\partial b+ia)$, and let $\Csm(Y,W)=\Csm(Y)/\Csm(W)$; since $\Csm(W)\to\Csm(Y)$ is injective, $\Csm(Y,W)$ is quasi-isomorphic to the cone of the inclusion. The map $(a,b)\mapsto(pa,pb)$ on cones gives $p:H_*(\Csm(Y,W))\to H_*(\GCG(Y,W))$, compatible with the long exact sequences of the pairs.

\begin{proposition}[relative comparison]\label{prop_rel}
\begin{enumerate}[(a)]
\item A map $\psi$ as in Section~\ref{sec_setting} with $\psi(Y)\subseteq Y'$ and $\psi(W)\subseteq W'$ induces a chain map $\GCG(Y,W)\to\GCG(Y',W')$. If $H$ is a homotopy as in Lemma~\ref{lem_prism} with $H(Y\times I)\subseteq Y'$ and $H(W\times I)\subseteq W'$, then $H_0$ and $H_1$ induce chain homotopic maps.
\item If $p_*$ is an isomorphism for $W$ and for $Y$, then $p_*:H_*(\Csm(Y,W))\to H_*(\GCG(Y,W))$ is an isomorphism. This holds when each of $W$ and $Y$ is a compact polyhedron or a smooth manifold without boundary, for instance for a compact polyhedron and a closed subpolyhedron, and for a manifold without boundary and an open subset or a closed submanifold; here a manifold may be any submanifold of $\R^K$ without boundary that is closed in an open subset of $\R^K$, see the proof. Then $H_*(\Csm(Y,W))\cong H_*(Y,W;\Q)$.
\item \emph{Excision.} If $(Y,W)\to(Y',W')$ is an inclusion of pairs as in \textup{(b)} inducing an isomorphism of singular homology, it induces an isomorphism $H_*(\GCG(Y,W))\to H_*(\GCG(Y',W'))$. For example, if $Y$ is a manifold without boundary, $W\subseteq Y$ open and $Z\subseteq W$ closed in $Y$, then $H_*(\GCG(Y\setminus Z,W\setminus Z))\cong H_*(\GCG(Y,W))$.
\end{enumerate}
The same holds for $GC^G_{\mathrm{orb}}$ and the targets of Theorem~\ref{thm_orb}, and for $GC^G_\star$, $\star\in\{\mathrm{lift},\mathrm{mor}\}$, with $Y$ an effective orbifold with corners and $W$ an open subset; there the maps in (a) are liftable maps, resp.\ morphisms, and the homotopies are liftable maps, resp.\ morphisms, $Y\times I\to Y'$, whose prism operators map $G^\star$ into $G^\star$ by (A4) and Lemma~\ref{lem_emb_mor}, as in the step \emph{Corners} of the proof of Theorem~\ref{thm_orbtarget}.
\end{proposition}

\begin{proof}
(a) $\psi_*$ commutes with the inclusions, so $(a,b)\mapsto(\psi_*a,\psi_*b)$ is a chain map of cones. With $\cH$ the operators of Lemma~\ref{lem_prism} for $W$ and $Y$, which commute with the inclusions, $(a,b)\mapsto(-\cH a,\cH b)$ is a chain homotopy between the maps induced by $H_0$ and $H_1$: $\partial(-\cH a,\cH b)+(-\cH(-\partial a),\cH(\partial b+ia))=(\partial\cH a+\cH\partial a,\partial\cH b+\cH\partial b+i\cH a-i\cH a)$, which is $(H_{1*}-H_{0*})(a,b)$.
(b) Five lemma on the long exact sequences of the cones. The singular chains compute $H_*(W;\Q)$ and $H_*(Y;\Q)$ by Proposition~\ref{prop_amb} and \cite{Lee}, and $p_*$ is an isomorphism for them by Theorem~\ref{thm_A_G}. Theorem~\ref{thm_A_G} is stated for closed submanifolds of $\R^K$, but for a submanifold $X$ without boundary that is closed in an open set $O'\subseteq\R^K$ (such as an open subset of a closed submanifold) the proofs of Propositions~\ref{prop_stokes} and \ref{prop_abs_mfd} apply unchanged: they use de Rham's theorem on $X$, a tubular neighbourhood $O\subseteq O'$ and a smooth retraction $\varrho:O\to X$, which exist for such $X$ \cite{Lee}.
(c) Naturality of $p$ and singular excision. For orbifolds use Theorem~\ref{thm_orbtarget} for the open sets and the singular homology of $|\cdot|$.
\end{proof}

\begin{remark}[the quotient complex]\label{rem_rel_quotient}
(1) In \cite{Yu1} the relative complex is the quotient of $C^{g,\infty}_*(X)$ by $C^{g,\infty}_{*,A}(X)+D^\infty_*(X)$, and excision is proved by cutting, for arbitrary subsets $A$ of smooth spaces, with no good cover assumed; the proof uses the locality of $D$ \cite[Lemma on locality of degenerate chains]{Yu1}. For $G$ the quotient is $\GCG(Y)/i(\GCG(W))=\Cm(Y)/(\Cm(W)+G(Y))$, and the natural map from $\operatorname{Cone}(i)$ to this quotient is a quasi-isomorphism if and only if $i:\GCG(W)\to i(\GCG(W))$ is, i.e.\ if and only if $(\Cm(W)\cap G(Y))/G(W)$ is acyclic. This is a form of the chain level locality of $G$, which is not known (Remark~\ref{rem_locality}); so we use the cone, and excision is obtained from the comparison theorem rather than by cutting.

(2) We do not treat the smooth spaces of \cite{Yu1} beyond the targets above. In all our proofs injectivity comes from integration (Propositions~\ref{prop_stokes} and \ref{prop_poly_inj}, Lemma~\ref{lem_dR_pair}) or, over $\Z$, from intersection with the dual cells of a smooth triangulation (Section~\ref{sec_int}); neither is available on a general smooth space, and the localisation of \cite[proof of the main theorem]{Yu1} uses $D$.
\end{remark}

\subsection{Products}
For $Y\subseteq\R^K$ and $Z\subseteq\R^{K'}$ put $[Q,f]\times[M,g]:=[Q\times M,f\times g]\in\Cm(Y\times Z)$, $Y\times Z\subseteq\R^{K+K'}$, with the product orientation. Products of sources are sources, the formula respects isomorphisms, orientation reversal and disjoint unions, and $\partial(x\times y)=\partial x\times y+(-1)^{|x|}x\times\partial y$, since $\partial(Q\times M)=\partial Q\times M\sqcup(-1)^{\dim Q}Q\times\partial M$ as oriented manifolds with corners \cite[Section on operations on geometric chains]{Yu1}.

\begin{lemma}\label{lem_prod_glue}
If $(P;(Q_i,\iota_i)_i)$ is a gluing datum and $M$ a source, then $(P\times M;(Q_i\times M,\iota_i\times\id)_i)$ and $(M\times P;(M\times Q_i,\id\times\iota_i)_i)$ are gluing data.
\end{lemma}

\begin{proof}
We treat $P\times M$. (G0)--(G2) are clear: products of orientation preserving embeddings are orientation preserving embeddings, and the images of the interiors are $\iota_i(\inter Q_i)\times\inter M$. The faces of $Q_i\times M$ are the products $F\times M_0$ of faces $F$ of $Q_i$ with components $M_0$ of $M$, and the products $Q_{i,0}\times F_M$ of components $Q_{i,0}$ of $Q_i$ with faces $F_M$ of $M$. The latter map into $P\times i_M(F_M)\subseteq B_{P\times M}$ and are outer. If $F$ is outer, so is $F\times M_0$. If $F$ is glued to $F'$ by $\psi$, then $F\times M_0$ is glued to $F'\times M_0$ by $\psi\times\id$, a diffeomorphism of manifolds with corners.
\end{proof}

\begin{proposition}[cross product]\label{prop_cross}
The product induces a chain map $\times:\GCG(Y)\otimes\GCG(Z)\to\GCG(Y\times Z)$. For singular simplices $\sigma,\upsilon$, $p(\sigma)\times p(\upsilon)-p(\mathrm{EZ}(\sigma\otimes\upsilon))\in\operatorname{span}\cG(Y\times Z)$, where $\mathrm{EZ}$ is the shuffle map. Hence, if $p_*$ is an isomorphism for $Y$, $Z$ and $Y\times Z$ (for instance if $Y$ and $Z$ are both compact polyhedra, or both manifolds without boundary), then under $p_*$ the product $\times$ is the singular cross product.
\end{proposition}

\begin{proof}
By Lemma~\ref{lem_prod_glue}, $r\times c\in\operatorname{span}\cG(Y\times Z)$ for $r\in\operatorname{span}\cG(Y)$ and $c\in\Cm(Z)$, and by the Leibniz rule $(\partial r)\times c=\partial(r\times c)-(-1)^{|r|}r\times\partial c\in G(Y\times Z)$; likewise in the second factor. So $\times$ descends, and it is a chain map by the Leibniz rule. The staircase triangulation $t$ of $\Delta^n\times\Delta^m$ is admissible, since each open simplex lies in a product of open faces (as for $t\times I$ in Section~\ref{ssec_adm}), and its simplices are affine embeddings, so it is Whitehead. Hence $[\Delta^n\times\Delta^m,\sigma\times\upsilon]-p\,\tau(\Delta^n\times\Delta^m,\sigma\times\upsilon,t)\in\cG(Y\times Z)$ by Lemma~\ref{lem_glue}(e). The chain $\tau(\Delta^n\times\Delta^m,\sigma\times\upsilon,t)$ is $\mathrm{EZ}(\sigma\otimes\upsilon)$, the sign of each shuffle simplex being its local degree for the product orientation. The last claim follows since $\mathrm{EZ}$ induces the cross product on smooth singular chains.
\end{proof}

\subsection{Pullbacks}
Let $X$ be a smooth manifold without boundary, $\pi:E\to X$ a smooth fibre bundle whose fibre $F$ is a closed oriented manifold, with orientation preserving structure group, and regard $X\subseteq\R^K$ and $E\subseteq\R^{K'}$ as closed submanifolds. For a geometric simplex $(P,x)$ in $X$ the fibre product $R=P\times_XE$ is locally $P\times F$, hence a compact manifold with corners, oriented base first, and we put $\pi^*[P,x]:=[R,\mathrm{pr}_E]$. The faces of $R$ are the components of $\partial P\times_XE$; a component lies over one face of $P$, so it is embedded, and $R\in\cS$.

\begin{lemma}\label{lem_pull_glue}
If $(P;(Q_i,\iota_i)_i)$ is a gluing datum and $x:P\to X$ is smooth, then $(P\times_XE;(Q_i\times_XE,\iota_i\times\id)_i)$, with $Q_i\times_XE$ formed with $x\iota_i$, is a gluing datum.
\end{lemma}

\begin{proof}
Over a trivialising open set $V\subseteq X$ the datum restricts to the product of the restriction of $(P;(Q_i,\iota_i))$ to $x^{-1}(V)$ with $F$, so (G0)--(G2) hold as in Lemma~\ref{lem_prod_glue}. The faces of $Q_i\times_XE$ are the components $C$ of $F_Q\times_XE$ for faces $F_Q$ of $Q_i$. If $F_Q$ is outer, $C$ maps into $B_P\times_XE=B_R$. If $F_Q$ is glued to $F'$ by $\psi$, then $\psi\times\id:F_Q\times_XE\to F'\times_XE$ is defined since $x\iota_j\psi=x\iota_i$, it is a diffeomorphism, and it glues $C$ to the component $(\psi\times\id)(C)$.
\end{proof}

\begin{proposition}[pullback]\label{prop_pull}
The map $\pi^*$ induces a chain map from $\GCG_*(X)$ to $\GCG_{*+\dim F}(E)$.
\end{proposition}

\begin{proof}
The construction respects isomorphisms, orientation reversal and disjoint unions, so it is defined on $\Cm(X)$. Since $\partial F=\emptyset$, $\partial R=\partial P\times_XE$ with the base first orientation, so $\partial\pi^*=\pi^*\partial$. By Lemma~\ref{lem_pull_glue}, $\pi^*\rel{\fd}{x}=\rel{\fd\times_XE}{\mathrm{pr}_E}\in\cG(E)$, since $\mathrm{pr}_E\circ(\iota_i\times\id)$ is the map of the piece $Q_i\times_XE$. Hence $\pi^*G(X)\subseteq G(E)$.
\end{proof}

\begin{remark}\label{rem_pull}
(1) If $\partial F\ne\emptyset$, $F\in\cS$ and the structure group maps each face of $F$ to itself, the same construction gives the boundary formula $\partial\pi^*x=\pi^*\partial x+(-1)^{|x|}(i_\partial)_*\pi_\partial^*x$ of \cite[Section on operations on geometric chains]{Yu1}, where $\pi_\partial:\partial_vE\to X$ is the bundle of faces of the fibres and $i_\partial:\partial_vE\to E$ its map, with $(i_\partial)_*\pi_\partial^*[P,x]:=[P\times_X\partial_vE,i_\partial\circ\mathrm{pr}]$ as in the Subsection on orbifold sources, and $\pi^*$ preserves $G$, with $E$ embedded in some $\R^{K'}$: by Lemma~\ref{lem_pull_glue} for $E$ and for $\partial_vE$ (a bundle whose fibre has corners, to which the proof applies verbatim), the relations of both pulled back data lie in $\cG(E)$, provided their faces are embedded. They are: the structure group preserves each face $F_k$ of $F$, so a face of $F_k$, a component of some $F_k\cap F_l$, is moved by the structure group only within $F_k\cap F_l$, whose components are disjoint, so $i$ is injective on each face of $P\times_XE$ and of $P\times_X\partial_vE$. Without the condition on faces, the faces of $P\times_XE$ need not be embedded: for $F$ a square, $X=S^1$ and monodromy the rotation by a quarter turn, the four edges of the fibre form one face of $E$, which meets itself at the vertices. The condition is sufficient, not necessary.
For orbifold sources no condition on faces is needed, and $F$ need not have embedded faces (Proposition~\ref{prop_orb_ops}).
\end{remark}

\subsection{Orbifold sources}
For orbifold sources faces need not be embedded, which removes the condition on faces in Remark~\ref{rem_pull}. Let $F$ be a compact oriented manifold with corners, $\pi:E\to X$ a smooth fibre bundle with fibre $F$ and orientation preserving structure group over a manifold $X$ without boundary, and $i_\partial:\partial_vE\to E$ the bundle of the boundaries $\partial F$ of the fibres, a bundle over $X$ with fibre $\partial F$ and orientation preserving structure group. Regard $X\subseteq\R^K$ as a closed submanifold and $E\subseteq\R^{K'}$ by a smooth embedding (the proof of Whitney's theorem \cite{Lee} applies to manifolds with corners, whose charts are restrictions of smooth maps on open subsets of $\R^n$); smooth maps into $E$ are then the smooth maps into $\R^{K'}$ with image in $E$. The fibre product $P\times_X\partial_vE$ below is formed in the same way, with fibre $\partial F$, and for $[P,x]\in C^{\mathrm{orb}}(X)$ we write $(i_\partial)_*\pi_\partial^*[P,x]:=[P\times_X\partial_vE,i_\partial\circ\mathrm{pr}]\in C^{\mathrm{orb}}(E)$; no complex of chains in $\partial_vE$ is used, since $i_\partial$ need not be injective. For $P\in\cS^{\mathrm{orb}}$ and $x:|P|\to X$ smooth, the fibre product $R=P\times_XE$ is the orbifold with corners with charts $(\tilde U\times W,H)$, where $(\tilde U,H)$ is a chart of $P$ with image $U$ such that $E$ is trivial over a neighbourhood of $x(U)$, $W\subseteq F$ is a chart, and $H$ acts on $\tilde U$ only. It is compact, oriented base first and effective, so $R\in\cS^{\mathrm{orb}}$, and $\partial R=\partial P\times_XE\sqcup(-1)^{\dim P}P\times_X\partial_vE$ as oriented orbifolds, the sign coming from $\partial(\tilde U\times W)$, as in the manifold case of \cite[Section on operations on geometric chains]{Yu1}. Put $\pi^*[P,x]:=[R,\mathrm{pr}_E]$; this respects isomorphisms, orientation reversal and disjoint unions, an orientation reversing automorphism $\phi$ of $(P,x)$ giving one, $\phi\times\id$, of $(R,\mathrm{pr}_E)$, so $\pi^*$ is defined on $C^{\mathrm{orb}}(X)$.

\begin{lemma}\label{lem_orb_ops_glue}
Let $(P;(Q_i,\iota_i)_i)\in\fG^{\mathrm{orb}}$.
\begin{enumerate}[(a)]
\item For $M\in\cS^{\mathrm{orb}}$, $(P\times M;(Q_i\times M,\iota_i\times\id)_i)$ and $(M\times P;(M\times Q_i,\id\times\iota_i)_i)$ are in $\fG^{\mathrm{orb}}$.
\item For $x:|P|\to X$ smooth, $(P\times_XE;(Q_i\times_XE,\iota_i\times\id)_i)$, with $Q_i\times_XE$ formed with $x\iota_i$, is in $\fG^{\mathrm{orb}}$.
\end{enumerate}
\end{lemma}

\begin{proof}
(a) Products of effective compact oriented orbifolds with corners are such. On underlying spaces $\iota_i\times\id$ is injective, and its lifts are products of embedding lifts with identities, so (G1$'$) holds; (G2) holds since the images cover $|P|\times|M|$ and $|\inter(Q_i\times M)|=|\inter Q_i|\times|\inter M|$. The faces of $Q_i\times M$ are the components of $\partial Q_i\times M$ and of $Q_i\times\partial M$. A component $C$ of $F\times M$, $F$ a face of $Q_i$, is outer if $F$ is, since $\iota_ii_{Q_i}(F)\times|M|\subseteq B_P\times|M|\subseteq B_{P\times M}$; if $F$ satisfies (G3$'$)(b) with $\psi:F\to F'$, then $\psi\times\id$ maps $C$ diffeomorphically onto a component of $F'\times M$ and $(\iota_j\times\id)(i_{Q_j}\times\id)(\psi\times\id)=(\iota_i\times\id)(i_{Q_i}\times\id)$ on $C$. The components of $Q_i\times\partial M$ map into $|P|\times B_M\subseteq B_{P\times M}$ and are outer.

(b) Over a trivialising open set $V\subseteq X$ the datum is the product of the restriction to $x^{-1}(V)$ with $F$, so (G0$'$)--(G2) hold as in (a). The faces of $Q_i\times_XE$ are the components of $\partial Q_i\times_XE$, treated as in (a) with $\psi\times\id$, which is defined since $x\iota_ji_{Q_j}\psi=x\iota_ii_{Q_i}$, and the components of $Q_i\times_X\partial_vE$, which map into the image of $\id\times i_\partial$, consisting of points of positive depth in the fibre, and are outer.
\end{proof}

\begin{proposition}\label{prop_orb_ops}
\begin{enumerate}[(a)]
\item The product induces a chain map $\times:GC^G_{\mathrm{orb}}(Y)\otimes GC^G_{\mathrm{orb}}(Z)\to GC^G_{\mathrm{orb}}(Y\times Z)$, and if $Y$ and $Z$ are both compact polyhedra or both manifolds without boundary, then under $p_*$ it is the singular cross product.
\item $\pi^*$ maps $G^{\mathrm{orb}}(X)$ into $G^{\mathrm{orb}}(E)$, and $\partial\pi^*x=\pi^*\partial x+(-1)^{|x|}(i_\partial)_*\pi_\partial^*x$, where $\pi_\partial:\partial_vE\to X$ and $i_\partial:\partial_vE\to E$. If $\partial F=\emptyset$, $\pi^*:GC^G_{\mathrm{orb},*}(X)\to GC^G_{\mathrm{orb},*+\dim F}(E)$ is a chain map.
\end{enumerate}
\end{proposition}

\begin{proof}
(a) As for Proposition~\ref{prop_cross}, with Lemma~\ref{lem_orb_ops_glue}(a) and the Leibniz rule for orbifolds \cite[Section on operations on geometric chains]{Yu1}; the staircase relations lie in $\cG\subseteq\cG^{\mathrm{orb}}$, and $p_*$ is an isomorphism by Theorem~\ref{thm_orb}. (b) The boundary formula follows from the formula for $\partial R$ and the definition of $(i_\partial)_*\pi_\partial^*$. By Lemma~\ref{lem_orb_ops_glue}(b), applied to the bundles $E$ and $\partial_vE$ over $X$, the data $\fd\times_XE$ and $\fd\times_X\partial_vE$ are in $\fG^{\mathrm{orb}}$ for $\fd\in\fG^{\mathrm{orb}}$; since a gluing relation does not depend on the map, the relations of the second with the map $i_\partial\circ\mathrm{pr}$ lie in $\cG^{\mathrm{orb}}(E)$. So $\pi^*$ and $(i_\partial)_*\pi_\partial^*$ map $\operatorname{span}\cG^{\mathrm{orb}}(X)$ into $\operatorname{span}\cG^{\mathrm{orb}}(E)$, and $\pi^*\partial r=\partial\pi^*r-(-1)^{|r|}(i_\partial)_*\pi_\partial^*r\in G^{\mathrm{orb}}(E)$.
\end{proof}

\begin{remark}\label{rem_orb_ops}
(1) For manifold sources, without a condition on faces the pullback can leave $\cS$ (the square example of Remark~\ref{rem_pull}); for orbifold sources the boundary of \cite[Convention on the boundary]{Yu1} handles non-embedded faces, and $\cS\subseteq\cS^{\mathrm{orb}}$, so the pullback of a manifold chain along a bundle as in the square example is defined in $GC^G_{\mathrm{orb}}$.
(2) For orbifold targets, the product of liftable maps is liftable and that of morphisms is a morphism, and by Lemma~\ref{lem_emb_mor} the morphisms of the pieces of the product data are the products of those of the factors, up to unique $2$-isomorphism; so Lemma~\ref{lem_orb_ops_glue}(a) gives a chain map $\times$ on $GC^G_\star$, the Leibniz rule using $\hat i_{Q\times M}|_{\partial Q\times M}\cong\hat i_Q\times\id$, by the uniqueness in Lemma~\ref{lem_emb_mor}. No comparison with the singular product is claimed there. Pullbacks along orbifold bundles are not treated (Remark~\ref{rem_orbtarget}(3)).
\end{remark}

\subsection{Morse complexes}
The relation $G$ kills a class of thin chains at chain level.

\begin{lemma}[product chains]\label{lem_prod_thin}
Let $M,N\in\cS$ with $\dim M\ge1$, and let $g:N\to\R^K$ be smooth with $g(N)\subseteq Y$. Then $[M\times N,g\circ\mathrm{pr}_2]\in\operatorname{span}\cG(Y)$.
\end{lemma}

\begin{proof}
Let $j=\dim M$ and $t\in\cT^{C^1}(M)$ (Theorem~\ref{thm_whitehead}). By Lemma~\ref{lem_glue}(e) and Lemma~\ref{lem_prod_glue}, $(M\times N;(\varepsilon_i\Delta^j\times N,\sigma_i\times\id)_i)$ is a gluing datum, and the maps of its pieces are $g\circ\mathrm{pr}_2$. So $[M\times N,g\circ\mathrm{pr}_2]\equiv\sum_i\varepsilon_i[\Delta^j\times N,g\circ\mathrm{pr}_2]$ modulo $\operatorname{span}\cG(Y)$. Since $j\ge1$, the affine map of $\Delta^j$ exchanging two vertices, times $\id_N$, is an orientation reversing automorphism of $(\Delta^j\times N,g\circ\mathrm{pr}_2)$, so each term is $0$ in $\Cm(Y)$.
\end{proof}

Let $M$ be a closed manifold with a Morse function $h$ and a Morse--Smale, locally trivial metric, and let $CM_*(h;\Q)$, the compactified descending manifolds $\overline U(a)$, the evaluation maps $e_a$ and the spaces $\overline{\mathcal M}(a,b)$ be as in \cite[Further applications, Morse theory]{Yu1}, after \cite{qin}.

\begin{proposition}[Morse complex]\label{prop_morse_G}
$\Phi(a):=[\overline U(a),e_a]$ defines a chain map $\Phi:CM_*(h;\Q)\to\GCG_*(M)$, depending only on $h$, the metric and the orientations of the descending manifolds. If $M$ is oriented, $\Phi$ is a quasi-isomorphism.
\end{proposition}

\begin{proof}
$\overline U(a)$ is a compact manifold with faces \cite[Theorems 3.3, 3.4]{qin}, oriented by the chosen orientation of $U(a)$, with embedded faces by step (i) of the proof of \cite[Proposition on the Morse complex]{Yu1}, so it is a source, and $e_a$ is smooth. Let $k=\mathrm{ind}\,a$. Each component of $\overline{\mathcal M}(a,b)\times\overline U(b)$ is diffeomorphic to a face of $\overline U(a)$, so it is a source, and so is $\overline{\mathcal M}(a,b)$, since a non-embedded face of $\overline{\mathcal M}(a,b)$ would give a non-embedded face of the product ($\overline U(b)\ne\emptyset$). By step (i), which takes place in the complex of classes before any relation is imposed, and since $\Cm(M)$ is spanned by a subset of the basis of that complex, in $\Cm_{k-1}(M)$
\[
\partial\Phi(a)=\Phi(da)+T_a,\qquad T_a=\sum_{\mathrm{ind}\,b\le k-2}[\overline{\mathcal M}(a,b)\times\overline U(b),e_b\circ\mathrm{pr}_2].
\]
For $\mathrm{ind}\,b\le k-2$, $\dim\overline{\mathcal M}(a,b)=k-\mathrm{ind}\,b-1\ge1$, so $T_a\in\operatorname{span}\cG(M)$ by Lemma~\ref{lem_prod_thin}, and $\partial\Phi=\Phi d$ in $\GCG(M)$. The independence of choices is as in \cite{Yu1}.

For the quasi-isomorphism, integration $I[Q,f](\omega)=\int_Qf^*\omega$ is a chain map from $\GCG(M)\otimes\R$ to currents, by Propositions~\ref{prop_stokes} and \ref{prop_additive} and Stokes' theorem. $I\circ p$ is integration over smooth singular chains, a quasi-isomorphism by de Rham's theorem, and $p\otimes\R$ is one by Theorem~\ref{thm_A_G}, $\R$ being flat over $\Q$; hence so is $I$. The rest is step (v) of the proof in \cite{Yu1}, which uses $I\circ\Phi$, the description of the faces in step (i), \cite[Lemma on Morse moduli spaces]{Yu1} (after \cite{qin}), de Rham's theorem for currents and the results of Harvey and Lawson \cite{hl,hl2}; since $\R$ is faithfully flat over $\Q$, $\Phi$ is then a quasi-isomorphism.
\end{proof}

In \cite{Yu1} the terms $T_a$ are degenerate and $\Phi$ is a chain map modulo $D$; here they lie in $\operatorname{span}\cG$, by Lemma~\ref{lem_prod_thin}.

\section{Continuous maps: outlook}\label{sec_cont}

When the maps are only continuous, $\Cm_{\rm top}(X)$ is spanned by continuous maps from compact oriented manifolds with corners, with faces not necessarily embedded, into a topological space $X$, and the relation $R_{\rm top}$, the quotient $GC^{R_{\rm top}}(X)$ and the kernel complex $K_{\rm top}(X)$ are defined as $R_\Delta$, $\GCRD$ and $K$ in Section~\ref{sec_R}, with topological triangulations that are admissible and depth-subordinate. Such triangulations exist by \cite[Theorem~1.1(a)]{Yu2}. The structure theorem holds verbatim, with the proof of Proposition~\ref{prop_structure}: $p$ is onto at chain level with kernel $K_{\rm top}$, so $p_*$ is an isomorphism in all degrees if and only if $K_{\rm top}(X)$ is acyclic. We do not study $K_{\rm top}$ here and only record some questions.

\subsection*{Open problems}
\begin{enumerate}[(1)]
\item Is $K_{\rm top}(X)$ acyclic for every topological space, or every metrizable absolute neighbourhood retract $X$? Already: does every $3$-cycle of $K_{\rm top}(S^3)$ have degree $0$?
\item Interpolation: given admissible topological triangulations $t,t'$ of $Q$ and one of $\partial Q\times I$ restricting to $t_\partial,t'_\partial$, is there an admissible triangulation of $Q\times I$ extending them? For $\dim Q\ge4$ topological triangulations can be non-combinatorial.
\item Locality of $K_{\rm top}$: is $K(X)\cap(\Calttop(U)+\Calttop(V))=K(U)+K(V)$ for open covers $X=U\cup V$? We expect that, by subdivision invariance, a Mayer--Vietoris argument reduces to this.
\end{enumerate}
The methods of this paper do not apply: they cut along regular levels, which continuous (or Lipschitz) maps do not have (Remark~\ref{rem_lipschitz}), and Theorems~\ref{thmA}--\ref{thmC} concern only smooth maps; only the injectivity argument and the low-degree arguments extend to Lipschitz maps (Theorem~\ref{thm_lip_inj}, Proposition~\ref{prop_lip_low}).

\section*{Declaration of generative AI and AI-assisted technologies}
During the preparation of this work the author used AI to check proofs and, under the author's direction, to draft some proofs and passages of text. The gluing relation and the overall approach are the author's. The author reviewed and edited all content and takes full responsibility for the paper.

\bigskip
\noindent Hao Yu, {\it Independent researcher, Beijing, China.}\\
\indent {\it Previously: School of Mathematics, University of Minnesota, Minneapolis, MN 55455, USA.}\\
{\tt dustless2014@163.com}

\end{document}